%% file: main.tex
\PassOptionsToPackage{dvipsnames}{xcolor}
\documentclass[onefignum,onetabnum,twoside]{siamart171218}

\usepackage{amsmath}
\usepackage{amssymb}
\usepackage{booktabs}
\usepackage{color}
\usepackage{comment}
\usepackage{gnuplot-lua-tikz}
\usepackage{mathtools}
\usepackage{subfig}
\usepackage{subfig}
\usepackage{tikz}

\usepackage{microtype}
\microtypesetup{
  tracking=true,
  kerning=true,
  spacing=true
}
\microtypecontext{spacing=nonfrench}

\newcommand{\Z}{\mathbb{Z}}
\newcommand{\R}{\mathbb{R}}

\newcommand{\e}{\varepsilon}
\newcommand{\eff}{\varepsilon^\text{eff}}

\renewcommand{\vec}[1]{\boldsymbol{#1}}
\newcommand{\vH}{\vec H}
\newcommand{\vE}{\vec E}

\newcommand{\vJ}{\vec J}
\newcommand{\vJa}{\vec J_{a}}

\newcommand{\vJSd}{\vec J_{\Sigma^d}}
\newcommand{\vME}{\vec{\mathcal{E}}}
\newcommand{\vMH}{\vec{\mathcal{H}}}
\newcommand{\vn}{\vec \nu}

\newcommand{\vx}{\vec x}
\newcommand{\vy}{\vec y}

\newcommand{\dx} {\,{\mathrm d}x}
\newcommand{\dox}{\,{\mathrm d}o_x}
\newcommand{\dy} {\,{\mathrm d}y}
\newcommand{\doy}{\,{\mathrm d}o_y}

\newcommand{\js}[1]{\left[#1\right]_{\Sigma}}
\newcommand{\je}[1]{\left\{#1\right\}_{\partial\Sigma}}
\newcommand{\jsd}[1]{\left[#1\right]_{\Sigma^d}}
\newcommand{\jse}[1]{\left\{#1\right\}_{\partial\Sigma^d}}

\headers
  {Homogenization of plasmonic crystals}
  {Maier, Mattheakis, Kaxiras, Luskin, Margetis}

\title{Homogenization of plasmonic crystals: Seeking the %
  epsilon-near-zero effect}

\author{
  Matthias Maier%
  \thanks{%
    Department of Mathematics, Texas A\&M University, College %
    Station, TX 77843, USA.}
  \and Marios Mattheakis%
  \thanks{%
    John A. Paulson School of Engineering and Applied Sciences, Harvard %
    University, Cambridge, MA 02138, USA.}
  \and Efthimios Kaxiras%
  \footnotemark[2]
  \thanks{%
    Department of Physics, Harvard University, Cambridge, MA 02138, USA.}
  \and Mitchell~Luskin%
  \thanks{%
    School of Mathematics, University of Minnesota, %
    Minneapolis, MN 55455, USA.}
  \and Dionisios Margetis%
  \thanks{%
    Department of Mathematics, and Institute for Physical Science and %
    Technology, and Center for Scientific Computation and Mathematical %
    Modeling, University of Maryland, College Park, MD 20910, USA.}
}

\begin{document}

\maketitle

\begin{abstract}
  By using an asymptotic analysis and numerical simulations, we derive and
  investigate a system of homogenized Maxwell's equations for conducting
  material sheets that are periodically arranged and embedded in a
  heterogeneous and anisotropic dielectric host. This structure is
  motivated by the need to design plasmonic crystals that enable the
  propagation of electromagnetic waves with no phase delay
  (epsilon-near-zero effect). Our microscopic model incorporates the
  surface conductivity of the two-dimensional (2D)
  material of each sheet and a corresponding line charge density through a
  line conductivity along possible edges of the sheets. Our analysis
  generalizes averaging principles inherent in previous Bloch-wave
  approaches. We investigate physical implications of our findings. In
  particular, we emphasize the role of  the vector-valued corrector field,
  which expresses microscopic modes of surface waves on the 2D material.
  \textcolor{black}{
    We demonstrate how our homogenization procedure may set the foundation for
    computational investigations of: effective optical responses of reasonably general geometries, and
   complicated design problems in the plasmonics of 2D materials.
  }
\end{abstract}

\section{Introduction}
\label{sec:intro}

The advent of two-dimensional (2D) materials with controllable electronic
structures has opened up an unprecedented wealth of optical phenomena that
challenge the classical diffraction limit of electromagnetic waves. Notable
examples of related applications range from optical
holography~\cite{wintz2015}, tunable metamaterials~\cite{nemilentsau2016},
and cloaking~\cite{alu2005}, to subwavelength focusing
lenses~\cite{cheng2014}. A striking feature of many of these applications
is the possible emergence of an unusual \color{black} parameter regime
\color{black} with no refraction, referred to as epsilon-near-zero (ENZ)
effect~\cite{maier2018a,mattheakis2016,huang2011,moitra2013,silveirinha2007}.

This effect calls for designing novel \emph{plasmonic crystals} made of
stacked metallic or semi-metallic 2D material structures arranged
periodically with subwavelength spacing, and embedded in a dielectric
host~\cite{maier2018a,mattheakis2016,maier2018b,caloz2005}. The
growing need to describe, engineer, and tune the optical properties of such
plasmonic crystals motivates the present paper. To gain insight into their
effective optical behavior, we utilize a homogenization procedure that
systematically illuminates how such macroscopic properties emerge from the
plasmonic microstructure.

\textcolor{black}{%
In this paper, our key objectives are: (i) to elucidate how effective
parameter values of plasmonic crystals are derived from two-scale
asymptotics; (ii) to demonstrate how this homogenization procedure can be
used to compute effective optical responses of complicated (periodic)
microstructures numerically; and (iii) to investigate physical implications
of the resulting, effective description for a few prototypical geometries.}
The homogenized system can be controlled by tuning (microscopic) geometry,
(periodic) spacing, frequency, and conductivity of the 2D material. This
description \color{black} implies \color{black} precise conditions for
the ENZ effect. Our work extends previous homogenization results for
time-harmonic Maxwell's
equations~\cite{maier2018b,wellander2003,amirat2011,amirat2016,guenneau2007}.
In particular, we \textcolor{black}{develop the following aspects of the
homogenization of plasmonic crystals:}

\begin{itemize}
  \item
    \textcolor{black}{We introduce a \emph{general} homogenization
    result for} a microscopic model that incorporates the surface
    conductivity of arbitrarily curved 2D material sheets \emph{and} a
    corresponding line charge density along possible edges of the sheets
    (see Section~\ref{sec:micro-formulation}). The former may give rise to
    a \emph{2D surface plasmon-polariton} (SPP) \cite{samaier07,bludov13}
    whereas the latter \color{black} may influence the appearance of
    \color{black} an \emph{edge plasmon-polariton} (EPP)
    \cite{krenn2004,fei2015}. These waves are special fine-scale surface
    \textcolor{black}{and edge} modes, \color{black} respectively.
    \color{black} The microscopic model describes a large class of
    plasmonic crystals consisting of periodic inclusions of metallic
    finite-size and curved 2D materials in a heterogeneous and anisotropic
    dielectric host (Sections~\ref{sec:micro-formulation},
    and~\ref{sec:asymptotics}).
  \item
    \textcolor{black}{We demonstrate analytically} how the combination
    of the complex-valued surface conductivity of the material sheets and
    the line conductivity along their edges can yield an ENZ effect
    (Section~\ref{sec:enz}). \textcolor{black}{In this framework, we
    derive precise ENZ conditions for planar sheets that extend results
    previously extracted from Bloch wave
    theory~\cite{mattheakis2016,maier2018a}.} In addition, by numerical
    simulations \color{black} based on the finite element method,
    \color{black} we discuss the ENZ effect occurring in a few prototypical
    geometries consisting of nanoscale 2D structures.
  \item
    \textcolor{black}{%
    We extend the homogenization procedure to include line charge models in
    correspondence to the line charge density of our microscopic description.
    This consideration accounts for possible edge effects.}
    \textcolor{black}{We discuss how this edge contribution extends known
    ENZ conditions (Section~\ref{sec:enz}), and discuss further implications of
    this possibility in the \textcolor{black}{conclusion}
    (Section~\ref{sec:conclusion}).}
  \item
    \textcolor{black}{We introduce a \emph{computational platform}
    based on our homogenization procedure that can serve as a foundation
    for investigating the effective optical response of (reasonably general)
    microscopic geometries (Section~\ref{sec:numerics}).} \color{black}
    This framework enables the systematic computational investigation of
    complicated design problems in the plasmonics of 2D materials
    \cite{molesky2018,liberal2017,miller2012} (Sections~\ref{sec:numerics}
    and~\ref{sec:conclusion}). \color{black}
\end{itemize}

We discuss several physical implications of our findings. For instance, we
identify the physical role of the corrector \textcolor{black}{field} in
our formulation: this field \color{black} contains microscopic wave
modes on the 2D material. We introduce a line charge density in the
modeling and homogenization procedure, and
\color{black} show how the line charge density
introduced in our microscopic model alters the ENZ effect of plasmonic
crystals described previously~\cite{maier2018a, mattheakis2016}.
\textcolor{black}{We recover a Lorentzian function} for the effective
dielectric permittivity tensor \textcolor{black}{of selected
prototypical geometries} that partly validates our results since its real
and imaginary parts automatically satisfy the Kramers-Kronig relations.

\subsection{Motivation: Epsilon-near-zero effect}
\label{subsec:intro:motivation}

Recently, 2D materials such as graphene and black phosphorus have been the
subject of extensive experimental and theoretical studies. From the
viewpoint of Maxwell's equations, the dielectric permittivity of a
conducting 2D material may have a negative real part. As a result,
SPPs of transverse-magnetic (TM) polarization may exist on the conducting
sheet with a dispersion relation that allows for a transmitted wavenumber,
$k_{\text{SPP}}$, much larger in magnitude than the free-space wavenumber,
$k_0$~\cite{low2017}.

Specifically, for an infinite, flat conducting sheet in an isotropic and
homogeneous ambient space, the SPP dispersion relation is~\cite{samaier07}
\begin{equation*}
  \sqrt{k_0^2-k^2_{\text{SPP}}} =
  -\left(\frac{2k_0}{\omega\mu_0\sigma}\right)\,k_0,
\end{equation*}
where $\sigma$ is the surface conductivity of the sheet, $\mu_0$ is the
magnetic permeability of the ambient space, and $\omega$ denotes the
angular frequency. Note that $|k_{\text{SPP}}|\gg k_0$ if
$\text{Im}\,\sigma >0$ and $|\omega\mu_{0}\sigma|\ll k_0$ for an assumed
$e^{-i\omega t}$ time dependence and lossless surrounding medium. Hence,
the wavelength of the TM-polarized SPP scales linearly with $\sigma$ if
dissipation is relatively small ($0<\text{Re}\sigma\ll \text{Im}\sigma$).

In general, the dispersion relation of transmitted waves through given 2D
materials can be altered dramatically by introducing different geometries
of the sheet, or different arrangements of sheets in a dielectric host. In
particular, plasmonic crystals are structures that consist of stacked,
periodically aligned metallic layers. When the period is of the order of
the SPP wavelength, unusual optical phenomena may occur, such as the ENZ
effect and negative
refraction~\cite{mattheakis2016,maier2017,silveirinha2007, moitra2013}.
These properties can be precisely controlled by tuning the electronic
structure of the 2D material, through chemical doping or other means, and
the operating frequency, $\omega$~\cite{xia2014,shirodkar2018}.

The ENZ effect implies that at least one eigenvalue of the \emph{effective}
dielectric permittivity of the plasmonic structure is close to zero. This
effect causes surprising optical features, which have not been obtained by
traditional photonic systems. These features include decoupling of spatial
and temporal field variations, tunneling through very narrow channels,
constant phase transmission, strong field confinement, diffraction-free
propagation, and ultrafast phase transitions~\cite{silveirinha2007,
niu2018}. Many novel functional devices based on plasmonic crystals have
been proposed, indicating the broad prospects of photonics
based on the ENZ effect~\cite{mattheakis2016,maier2017,silveirinha2007,
moitra2013, niu2018, xia2014}.

Motivated by this perspective, in this paper we develop a general
homogenization procedure for plasmonic crystals. In addition, we
investigate the possibility for emergence of the ENZ effect in prototypical
geometries \color{black} with, e.g., graphene layers, nanoribbons, and nanotubes. \color{black}

\subsection{Microscopic model and geometry}
\label{subsec:intro:micro}
\textcolor{black}{The geometry is shown in Figure~\ref{fig:geometry}. It
consists of periodically  stacked, possibly curved sheets, $\Sigma^d$, of a
2D material with surface conductivity $\sigma^d(\vx)$. We assume that a
charge accumulation may occur via a \emph{line conductivity},
$\lambda^d(\vx)$, on the edges $\partial\Sigma^d$ of the sheets. The sheets
are embedded in a dielectric with \emph{heterogeneous} permittivity
$\e^d(\vx)$.}
At the microscale, we invoke time-harmonic Maxwell's equations for the
electromagnetic field $(\vE^d,\vH^d)$ in domain $\Omega$; see
\eqref{eq:maxwell}. The conductivities $\sigma$ and $\lambda$ are
responsible for the induced current density
\begin{align*}
  \vJSd=\vec\delta_{\Sigma^d}\sigma^d\vE^d+\vec\delta_{\partial\Sigma^d}\lambda^d\vE^d.
\end{align*}
\textcolor{black}{Here, $\delta_{\Sigma^d}$ denotes the Dirac delta
function associated with the (possibly curved) surface of $\Sigma^d$, and
$\delta_{\partial\Sigma^d}$ is the Dirac delta function associated
with the boundaries of $\Sigma^d$. The induced current density $\vJSd$ will
lead to jump conditions of electromagnetic field components over $\Sigma^d$
and $\partial\Sigma^d$. Detailed discussions and derivations of the
governing equations and jump conditions are given in
Section~\ref{sec:micro-formulation}.}

\subsection{Homogenized theory}
\label{subsec:intro:effective}

We will demonstrate by a formal \emph{asymptotic analysis} that in the case
of scale separation, meaning if $d$ is sufficiently small compared to the
free-space wavelength, $2\pi/k_0$, the above problem manifesting the
microstructure can be expressed by the following \emph{homogenized system}:
\textcolor{black}{%
\begin{align}
  \label{eq:intro:maxwellhomogenized}
    \nabla\times\vME = i\omega\mu_0\vMH,
    \qquad
    \nabla\times\vMH = -i\omega\eff\vME+\vJa,
\end{align}}
in which the dependence on (microscale) spacing $d$ is eliminated. Here,
$\vJa$ is a \textcolor{black}{current-carrying} source, $(\vME, \vMH)$
describes an effective electromagnetic field, \textcolor{black}{and
$\eff$ is an \emph{effective permittivity tensor}.}

\textcolor{black}{We assume that the functions $\e$, $\sigma$ and $\lambda$
depend on a slow scale and are periodic and rapidly oscillating on a
(fast) scale proportional to the \emph{scaling parameter} $d$; typically,
$d\ll 2\pi/k_0$:
\begin{align*}
  \e^d(\vx) = \e(\vx,\vx/d),
  \qquad
  \sigma^d\,=\,d\,\sigma(\vx,\vx/d),
  \qquad
  \lambda^d=d^2\,\lambda(\vx,\vx/d).
\end{align*}
This scaling assumption leads to an effective permittivity tensor, viz.,}
\begin{multline*}
  \eff_{ij}:=
  \int_Y\e(\vx,\vy)\big(\vec e_j+\nabla_y\chi_j(\vx,\vy)\big)
  \cdot\vec e_i \dy
  -\frac1{i\omega}
  \int_\Sigma
  \sigma(\vx, \vy)\big(\vec e_j+\nabla_y\chi_j(\vx,\vy)\big) \cdot \vec
  e_i\doy
  \\
  -\frac1{i\omega}
  \int_{\partial\Sigma}
  \lambda(\vx, \vy)\big(\vec e_j+\nabla_y\chi_j(\vx,\vy)\big) \cdot \vec
  e_i\,\text{d} s.
\end{multline*}
In the above equation, $Y$ denotes the unit cell with embedded boundary
$\Sigma$ whose edge is $\partial\Sigma$ (see Figure~\ref{fig:geometry});
$\vec e_i$ is the unit vector in the $i$-th direction; and
$\nabla_y\vec\chi$ is the Jacobian of the \emph{corrector}, $\vec\chi$,
that solves the associated cell problem which encodes the microscopic
details, cf.~\eqref{eq:cell_problem}. A crucial property of $\eff$ is that
it manifests an interplay, and possible mutual cancellation, of three
distinct averages. Thus, by tuning geometry, \color{black} (periodic) spacing, \color{black} frequency and
conductivities of the 2D material, we can force one or more eigenvalues of
$\eff$ to be close to zero. This ENZ effect is detailed in
Section~\ref{sec:enz}.

\subsection{Past works}
\label{subsec:intro:past-lit}

The derivation presented in this paper is based on a formal asymptotic
analysis in the spirit of~\cite{Bensoussan1978}. Note that
in~\cite{maier2018b} a rigorous approach invoking two-scale convergence is
applied to plasmonic crystals without a line charge density along edges. In
the framework of homogenization theory for time-harmonic Maxwell's
equations, we should also mention a number of other related, rigorous or
formal, results~\cite{Nevard97, wellander2003, guenneau2007,
amirat2011, amirat2016}. In particular, in \cite{guenneau2007} a formal
asymptotic analysis is applied to \emph{finite} photonic crystals; and in
\cite{Nevard97} the authors homogenize Maxwell's equations in the presence
of rough boundaries (see also~\cite{kristensson2005}).

Broadly speaking, the design of  structures with unusual optical properties
is a highly active direction of research. For recent advances in photonics,
we refer the reader to~\cite{Jahani2016,Zheludev2016}. In fact, the
computation of effective material parameters in Maxwell's equations has a
long history in physics and engineering. It is impossible to exhaustively
list the related bibliography here.  In regard to the homogenization of
periodic and heterogeneous systems, we mention as
examples~\cite{kristensson2005, shelukhin2009, caloz2005, lalanne1996}.

Notably, homogenization results for Maxwell's equations can be related to
approaches based on \emph{Bloch theory} for waves in periodic structures.
In Section~\ref{sec:enz}, we demonstrate  that our homogenization result
readily generalizes an averaging principle that was previously found in
particular Bloch-wave solutions constructed for simple
settings~\cite{maier2018a, mattheakis2016}. In the context of layered
structures usually only a few Bloch waves effectively contribute to the
macroscopic field \cite{sjoeberg2005}. In contrast, our approach relies
entirely on periodic upscaling principles and is thus \emph{independent} of
the choice of particular solutions. In fact, we derive effective equations
and material parameters that are valid for a wide range of geometries
without choosing any particular solutions \emph{a priori}.

\subsection{Limitations}
\label{subsec:intro:limitations}

Although the homogenization results presented in this paper are fairly general, our analysis bears limitations. We should mention the following
issues:
\begin{itemize}
  \item
    Our analysis leaves out \emph{boundary effects} in the homogenization
    procedure due to the interaction of the microstructure (conducting
    sheet) with boundaries of the domain. This simplification restricts the
    homogenization result either to layered microstructures immersed in a
    scattering domain, or to domains with periodic boundary conditions.
  \item
    Our asymptotic analysis relies on a strong periodicity assumption for
    the microstructure. Even though we allow material parameters to also
    depend on a slow scale, \color{black} we do not account for a
    slowly varying geometry of the microstructure.
  \item
    We assume a scale separation between the free-space wavenumber, $k_0$,
    which we treat as of the order of unity ($k_0\sim 1$), and the
    \emph{SPP wavenumber}, $k_{\text{SPP}}$, with $k_{\text{SPP}}\sim 1/d$.
    This assumption rules out resonance effects, if the period (spacing),
    $d$, is close to the free-space wavelength. In this vein, we do not
    discuss the case with resonant scaling, $k_{\text{SPP}}\sim 1$.
\end{itemize}

\subsection{Paper organization}
\label{subsec:intro:organiz}

The remainder of the paper is organized as follows. In
Section~\ref{sec:micro-formulation}, we introduce the microscale model and
scaling assumptions. In Section~\ref{sec:asymptotics}, we derive a system
of homogenized equations and the corresponding cell problem, which
determine the effective material parameters. Section~\ref{sec:enz} provides
a detailed discussion of the ENZ effect from our homogenization. In
Section~\ref{sec:numerics}, we complement our analytical findings by
\textcolor{black}{a demonstration via examples of how the
homogenization results can serve as a platform for computing the effective
optical response of various microscopic geometries.}
Section~\ref{sec:conclusion} concludes the paper with a summary of the key
results and an outline of open problems. In the appendix \color{black}
we provide detailed analytical derivations needed in the main text, for the
sake of brevity. \textcolor{black}{The supplementary material contains
some additional analytical derivations that are not essential for our
paper.} \color{black} We use the  $e^{-i\omega t}$ time dependence
throughout. \color{black}

\section{Microscopic theory and geometry}
\label{sec:micro-formulation}
In this section, we give a detailed description of the (full) microscale
model and geometry. We also describe our scaling assumptions and discuss
the rationale underlying them.
\begin{figure}
  \begin{center}
    \subfloat[]{
      \begin{tikzpicture}[scale=0.90]
        \path[thick,draw,->] (-2.5,-0.6) -- (-2.0,-0.6);
        \path[thick,draw,->] (-2.5,-0.6) -- (-2.5,-0.1);
        \path[thick,draw,->] (-2.5,-0.6) -- (-2.8,-0.9);
          \node at (-2.5, +0.1) {\small $y_2$};
          \node at (-1.8, -0.6) {\small $y_1$};
          \node at (-3.0, -0.9) {\small $y_3$};
      \end{tikzpicture}
      \hspace{-3em}
      \begin{tikzpicture}[scale=0.90]
        \path [thick, draw]
          (-1, -1) -- (1, -1) -- (1, 1) -- (-1, 1) -- cycle;
        \path [thick, draw] (1,-1) -- (1.5, -0.5) -- (1.5, 1.5) -- (1, 1);
        \path [thick, draw] (-1,1) -- (-0.5, 1.5) -- (1.5, 1.5);
        \path [thick, draw, dashed]
          (-0.5, 1.5) -- (-0.5, -0.5) -- (1.5, -0.5);
        \path [thick, draw, dashed] (-0.5, -0.5) -- (-1.0, -1.0);
        \path [fill=black, fill opacity=0.2]
          (-1, 0) to[out=-30,in=170] (1, 0) to[out=-30,in=170] (1.5, 0.5)
          to[out=170,in=-30] (-0.5, 0.5) to[out=170,in=-30] (-1.0, 0.0);

        \path [thick, draw] (-1, 0) to[out=-30,in=170] (1, 0);
        \path [thick, draw] (1, 0) to[out=-30,in=170] (1.5, 0.5);
        \path [thick, draw, dashed] (-0.5, 0.5) to[out=-30,in=170] (1.5, 0.5);
        \path [thick, draw, dashed] (-1.0, 0.0) to[out=-30,in=170] (-0.5, 0.5);
        \path [thick, draw] (-1.4, 0.5) -- (0.0, 0.1);
        \node at (-1.6, 0.5) {\small $\Sigma$};
      \end{tikzpicture}}
    \hspace{3em}
    \subfloat[]{
      \begin{tikzpicture}[scale=0.90]
          \path[thick,draw,->] (-2.5,-0.6) -- (-2.0,-0.6);
          \path[thick,draw,->] (-2.5,-0.6) -- (-2.5,-0.1);
          \path[thick,draw,->] (-2.5,-0.6) -- (-2.8,-0.9);
          \node at (-2.5, +0.1) {\small $x_2$};
          \node at (-1.8, -0.6) {\small $x_1$};
          \node at (-3.0, -0.9) {\small $x_3$};
        \foreach \x in {0,...,6} {
          \path [fill=white]
            (-1.5, -1+0.3*\x) to[out=-30,in=170] (-1.0, -1+0.3*\x)
                               to[out=-30,in=170] (-0.5, -1+0.3*\x)
                               to[out=-30,in=170]  (0.0, -1+0.3*\x)
                               to[out=-30,in=170]  (0.5, -1+0.3*\x)
                               to[out=-30,in=170]  (1.0, -1+0.3*\x)
                               to[out=-30,in=170]  (1.5, -1+0.3*\x)
                               to[out=-30,in=170] (1.75, -0.75+0.3*\x)
                               to[out=-30,in=170] (2.00, -0.5+0.3*\x)
                               to[out=170,in=-30] (1.50, -0.5+0.3*\x)
                               to[out=170,in=-30] (1.00, -0.5+0.3*\x)
                               to[out=170,in=-30] (0.50, -0.5+0.3*\x)
                               to[out=170,in=-30] (0.00, -0.5+0.3*\x)
                               to[out=170,in=-30] (-0.5, -0.5+0.3*\x)
                               to[out=170,in=-30] (-1.0, -0.5+0.3*\x)
                               to[out=170,in=-30] (-1.25, -0.75+0.3*\x)
                               to[out=170,in=-30] cycle;
          \path [thick, draw, fill=black, fill opacity=0.2]
            (-1.5, -1+0.3*\x) to[out=-30,in=170] (-1.0, -1+0.3*\x)
                               to[out=-30,in=170] (-0.5, -1+0.3*\x)
                               to[out=-30,in=170]  (0.0, -1+0.3*\x)
                               to[out=-30,in=170]  (0.5, -1+0.3*\x)
                               to[out=-30,in=170]  (1.0, -1+0.3*\x)
                               to[out=-30,in=170]  (1.5, -1+0.3*\x)
                               to[out=-30,in=170] (1.75, -0.75+0.3*\x)
                               to[out=-30,in=170] (2.00, -0.5+0.3*\x)
                               to[out=170,in=-30] (1.50, -0.5+0.3*\x)
                               to[out=170,in=-30] (1.00, -0.5+0.3*\x)
                               to[out=170,in=-30] (0.50, -0.5+0.3*\x)
                               to[out=170,in=-30] (0.00, -0.5+0.3*\x)
                               to[out=170,in=-30] (-0.5, -0.5+0.3*\x)
                               to[out=170,in=-30] (-1.0, -0.5+0.3*\x)
                               to[out=170,in=-30] (-1.25, -0.75+0.3*\x)
                               to[out=170,in=-30] cycle;
        }
        \node at (-2.25, +1.0) {\small $\Omega$};
        \node at (-2.25, +0.6) {\small $\varepsilon^d(\vx)$};
        \node at (2.5, +0.8) {\small $\Sigma^d$};
        \path [thick, draw] (2.3, 0.75) -- (1.8, 0.55);
        \node at (2.5,  0.4) {\small $\sigma^d(\vx)$};
      \end{tikzpicture}}
  \end{center}
  \caption{Schematic of geometry. (a) The unit cell, $Y=[0,1]^3$, with
    microstructure $\Sigma$, a conducting sheet. (b) Computational domain
    $\Omega$ with rescaled periodic layers $\Sigma^d$ and spatially
    dependent surface conductivity $\sigma^d(\vx)$. The ambient medium has
    a heterogeneous permittivity, $\e^d(\vx)$.}
  \label{fig:geometry}
\end{figure}
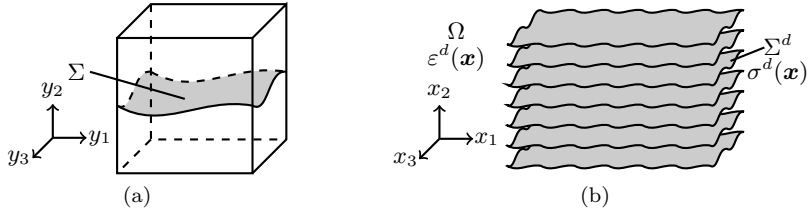

Our objective is to study Maxwell's equations in the (not necessarily
bounded) domain $\Omega\subset\R^3$ that contains microstructures with
periodically aligned, plasmonic hypersurfaces $\Sigma^d$;
see~Figure~\ref{fig:geometry}. These hypersurfaces are generated from the
inscribed hypersurface $\Sigma$ via scaling with a length scale, $d$, and
periodically repeating a \emph{unit cell}, $Y=[0,1]^3$. More precisely, we
define
\begin{align*}
  \Sigma^d = \big\{d\,\vec z + d\,\vec \varsigma\;\;:\;\;\vec z \in \Z^3,\;
  \vec \varsigma\in\Sigma\big\},
\end{align*}
where $\Z$ denotes the set of integers. Note that we make no specific
assumption about $\Sigma$, except that it is smooth (in $Y$ and when
periodically continued) so that it admits a well-defined surface normal.
This property does neither imply that $\Sigma$ has to be a connected
manifold, nor that $\Sigma$ has to have an intersection with the boundary
of $Y$; see Section~\ref{sec:numerics} for examples.

The relevant material parameters can be heterogeneous and tensor-valued. We
assume that the ambient medium, which is contained in
$\Omega\setminus\Sigma^d$, is described by a dielectric permittivity,
$\e^d(\vx)$. Furthermore, the plasmonic sheets $\Sigma^d$ have the
associated surface conductivity $\sigma^d(\vx)$. In addition, we allow for
a \emph{line conductivity} $\lambda^d(\vx)$ on the edges,
$\partial\Sigma^d$, of the sheets.

To derive a physically appealing homogenization limit, we assume that the
spatial dependence of the material parameters can be separated into a
slowly oscillating macroscale behavior and a (locally) periodic microscale
behavior with structural period $d$ \cite{Bensoussan1978,pavliotis2007},
viz.,
\begin{align*}
  \e^d(\vx) = \e\big(\vx,\vx/d\big),\qquad
  \sigma^d(\vx) = d\,\sigma\big(\vx,\vx/d\big),
  \qquad
  \lambda^d(\vx) = d^2\,\lambda\big(\vx,\vx/d\big)~.
\end{align*}
Here, $\e(\vx,.)$, $\sigma(\vx,.)$, and $\lambda(\vx,.)$ are $Y$-periodic
and tensor-valued, i.e., $\e\big(\vx,\vy+\vec e_i\big) =
\e\big(\vx,\vy\big)$, $\sigma\big(\vx,\vy+\vec e_i\big) =
\sigma\big(\vx,\vy\big)$, and $\lambda\big(\vx,\vy+\vec e_i\big) =
\lambda\big(\vx,\vy\big)$, for any unit vector $\vec e_i$ and
$\vx\in\Omega$, $\vy\in\R^3$.

\textcolor{black}{%
The surface conductivity $\sigma(\vx,\vy)$ and line conductivity
$\lambda(\vx,\vy)$ deserve particular attention. The former describes the
linear optical response of the (possibly curved) hypersurface $\Sigma$, or induced
surface current along $\Sigma$, due to an incident electric field.
Similarly, $\lambda$ describes the linear response of the
edges, $\partial\Sigma$, of $\Sigma$, or induced line current along
$\partial\Sigma$, due to an incident electric field. These considerations imply that both
tensors,  $\sigma(\vx,\vy)$ and
$\lambda(\vx,\vy)$, are rank deficient with one, or two zero eigenvalues, respectively.
More precisely, by fixing $(\vx,\vy)$ and expressing the
tensors $\sigma(\vx,\vy)$ and $\lambda(\vx,\vy)$ in local coordinates
$(\vec \tau, \vec n, \vn)$, where $\vec \tau$ is a vector field parallel to
the edge $\partial\Sigma$, $\vn$ is the (interior) normal on $\Sigma^d$,
and $\vec n$ denotes the outward-pointing unit vector orthogonal to
$\vec\tau$ and $\vn$, we have} \textcolor{black}{
\vspace{-1em}
\begin{align*}
  \sigma^d\;=\;
  \begin{pmatrix}
    \sigma_{11} & \sigma_{12} & 0 \\
    \sigma_{21} & \sigma_{22} & 0 \\
    0 & 0 & 0
  \end{pmatrix},
  \qquad
  \lambda^d\;=\;
  \begin{pmatrix}
    \lambda_{11} & 0 & 0 \\
    0 & 0 & 0 \\
    0 & 0 & 0
  \end{pmatrix}.
\end{align*}
}
Since we will apply a formal asymptotic analysis, our assumptions on the
\color{black} related parameters \color{black} are not too restrictive.
Consequently, we use arbitrary tensor-valued functions $\e$ and $\sigma$. A
mathematically rigorous convergence result typically requires more
restrictive assumptions \cite{maier2018b}.

\subsection{On the scaling of conductivities}
\label{subsec:scaling}

Our particular choice of scalings of the surface and line conductivities
with $d$, viz.,
\begin{align*}
  \sigma^d\sim d~,\qquad \lambda^d\sim d^2~,
\end{align*}
deserves some explanation. \textcolor{black}{We recall} that the
wavenumber, $k_{\text{SPP}}$, of the desired, fine-scale SPP on an infinite
conducting sheet scales inversely proportional to the surface conductivity
(Section~\ref{sec:intro}). Hence, our choice of scaling of $\sigma^d$
distinctly separates two scales: one related to the free-space wavenumber
$k_0$ ($k_0\sim 1$), determined by the average of $\e^d$, and another for
the SPP wavenumber, $k_{\text{SPP}}\sim1/\sigma^d\sim\,1/d$, on the
conducting sheets~\cite{maier2017}. By our scaling, the \emph{interaction
range} of the SPP on each sheet is of the order of $d$. More precisely, in
the limit $d\to 0$ the strength with which SPPs on one sheet influence
neighboring sheets remains constant. In addition, our choice of scaling
implies that the total surface current on the sheets remains \emph{finite}.
In a similar vein we scale the line conductivity, $\lambda^d$, with $d^2$.
In the limit $d\to 0$, this assumption yields a \emph{finite} total line
current on the edges of the conducting sheets.

It is worth mentioning that other scaling scenarios may lead to different
homogenization results. For example, the assumption $\sigma^d\sim 1$
corresponds to a resonance which in turn yields an effective
\emph{permeability} tensor~\cite{schweizer2010}. Such a configuration can
exhibit a respective \emph{mu-near-zero} effect. A variety of structures,
e.g., nanorings, that exhibit a mu-near-zero effect are studied
in~\cite{papadakis2018}.

\subsection{Heterogeneous Maxwell's equations}
\label{subsec:maxwell}

The time-harmonic Maxwell equations for the electromagnetic field $(\vE^d,
\vH^d)$ in $\Omega\setminus\Sigma^d$ are:
\begin{align}
  \begin{cases}
    \begin{aligned}
      \nabla\times\vE^d &= i\omega\mu_0\vH^d~,
      \quad
      \nabla\times\vH^d = -i\omega\e^d\vE^d+\vJa~,
      \\[0.3em]
      \nabla\cdot(\e^d\vE^d) &= \frac{1}{i\omega}{\nabla\cdot\vJa}~,
      \quad
      \nabla\cdot\vH^d = 0~.
    \end{aligned}
  \end{cases}
  \label{eq:maxwell}
\end{align}
The surface conductivity, $\sigma^d$, is responsible for the appearance of
the current density $\vJSd=\vec\delta_{\Sigma^d}\sigma^d\vE^d$, on
$\Sigma^d$. Accordingly, we must impose the following jump conditions
\textcolor{black}{on $\Sigma^d$ (away from the boundaries
$\partial\Sigma^d$):}
\begin{align}
  \begin{cases}
    \begin{aligned}
      \jsd{\vn\times\vE^d} &= 0~,
      &\,
      \jsd{\vn\times\vH^d} &= \sigma^d\vE^d~,
      \\[0.3em]
      \jsd{\vn\cdot(\e^d\vE^d)}  &=
      \frac{1}{i\omega}\nabla\cdot(\sigma^d\vE^d)~,
      &\,
      \jsd{\vn\cdot\vH^d} &= 0~.
    \end{aligned}
  \end{cases}
  \label{eq:jump}
\end{align}
Here, $[\,.\,]_{\Sigma^d}$ denotes the jump over $\Sigma^d$ with respect to
a chosen normal $\vn$, viz.,
\begin{align*}
  \jsd{\vec F}(\vx)
  \,:=\,
  \lim_{\alpha\searrow0}
  \Big(\vec F(\vx+\alpha\vn) - \vec F(\vx-\alpha\vn)\Big)\qquad
  \vx\in\Sigma^d.
\end{align*}
Equations \eqref{eq:maxwell} and \eqref{eq:jump} are supplemented  with the
following internal boundary conditions on the edges of the plasmonic
sheets, $\Sigma^d$, due to the line conductivity $\lambda^d$:
\begin{align}
  \begin{cases}
    \begin{aligned}
      \jse{\vec n\times\vH^d}\;&=\;\lambda^d\vE^d,
      \\[0.3em]
      \vec n \cdot\big(\sigma^d\vE^d\big)\;&=\;
      \nabla\cdot\big(\lambda^d\vE^d\big) \quad\text{on
      }\partial\Sigma^d,
    \end{aligned}
  \end{cases}
  \label{eq:compatibility_condition}
\end{align}
where $\vec n$ denotes the outward-pointing unit vector tangential to the
2D sheet $\Sigma^d$ and normal to curve
$\partial\Sigma^d$. In~\eqref{eq:compatibility_condition}, the symbol
$\jse{\,.\,}$ denotes a \emph{singular} jump over the edge,
$\partial\Sigma^d$, viz.,
\begin{align}
  \jse{\vec F}(\vx)
  \,:=\,
  \lim_{\alpha\searrow0}\,
  \int_{-\alpha}^{\;\alpha}
  \big(\vec F(\vx+\alpha^2\vec n+\zeta\vn)-\vec F(\vx-\alpha^2\vec
  n+\zeta\vn)\big)\,
  \text{d}\zeta\qquad
  \vx\in\partial\Sigma^d.
  \label{eq:singular_jump}
\end{align}
The jump condition in \eqref{eq:compatibility_condition} stems from the
property that, by virtue of the Amp\`ere-Maxwell law, the line current
$\lambda^d\vE^d$ creates a \emph{singular} magnetic field, $\vH\sim
(1/r)\,\vec e_\varphi$ as $r\to 0$, where $r$ denotes the distance to the
edge, $\partial\Sigma^d$, and $\vec e_\varphi$ is the unit vector in the
azimuthal direction of the local cylindrical coordinate system. In this
sense, jump \eqref{eq:singular_jump} measures the strength of the
singularity of the magnetic field $\vH$ at any point
$\vx\in\partial\Sigma^d$. We refer the reader to
Section~\ref{sec:gauss_law} of the supplementary material for a detailed
discussion including a derivation of the jump $\jse{\,.\,}$.

The second compatibility condition in \eqref{eq:compatibility_condition} is
a direct consequence of charge conservation along the edges,
$\partial\Sigma^d$, in view of the fact that $\sigma^d$ vanishes
identically outside the sheets. Let us explain. The \emph{line charge
accumulation} at the edges is related to the jump of the electric
displacement field in the $\vec n$-direction over the edge,
$\partial\Sigma^d$. This jump is in turn equal to $(i\omega)^{-1} \vec n
\cdot\big(\sigma^d\vE^d\big)$, since $\sigma^d\equiv 0$ outside $\Sigma^d$.
The corresponding line conductivity must be balanced by the divergence of
the line current $\lambda^d \vE^d$. Alternatively, compatibility conditions
\eqref{eq:compatibility_condition} can be derived from a variational
formulation; see Section~\ref{sec:weak_formulation} of the supplementary
material.

Equations \eqref{eq:maxwell}--\eqref{eq:compatibility_condition} form a
closed system if they are complemented by suitable boundary conditions on
$\partial\Omega$. We can impose either \emph{absorbing} conditions, or
\emph{periodic} conditions, or, in the case with an unbounded domain,
$\Omega$, Silver-M\"uller radiation conditions at
infinity~\cite{muller69,samaier07}. We make the important assumption that
the domain, $\Omega$, and the chosen boundary conditions on
$\partial\Omega$ are compatible with the microstructure in the sense of the
following simultaneous requirements:
\textcolor{black}{(i) the microstructure is bounded, that is
$\sigma(\vx,.)=0$ and $\e(\vx,.)=\e_0$ for $|\vx|\ge C$ for some $C>0$;
and (ii) the microstructure only touches those parts of the boundary
$\partial\Omega$ with periodic boundary conditions.}
These requirements ensure that the formal homogenization result derived in
this paper remains valid everywhere. If any one of these requirements is
not satisfied, a \emph{boundary layer} might occur, which necessitates a
\color{black} subtle homogenization procedure \color{black} near
$\partial\Omega$~\cite{pavliotis2007}.

\section{Two-scale expansions and asymptotics}
\label{sec:asymptotics}

In this section, we provide a derivation of effective system
\eqref{eq:intro:maxwellhomogenized} by means of a formal asymptotic
analysis. For the sake of brevity, we summarize the derivation here and
refer the reader to Appendix~\ref{sec:detailed_asymptotics} for details.

\textcolor{black}{We start with the microscale description of
\eqref{eq:maxwell}, \eqref{eq:jump} and \eqref{eq:compatibility_condition}.
Consider the two-scale expansions~\cite{Bensoussan1978, pavliotis2007}}
\begin{align*}
  \vE^d & \to \vE^{(0)}(\vx,\vy) + d \vE^{(1)}(\vx,\vy) + d^2
  \vE^{(2)}(\vx,\vy) + \ldots,
  \\
  \vH^d & \to \vH^{(0)}(\vx,\vy) + d \vH^{(1)}(\vx,\vy) + d^2
  \vH^{(2)}(\vx,\vy) + \ldots,
\end{align*}
with the corresponding substitutions
\begin{align*}
  \nabla \to \nabla_x+\frac1d\nabla_y,
  \qquad
  \e^d \to \e(\vx,\vy),
  \qquad
  \sigma^d \to d\sigma(\vx,\vy),
  \qquad
  \lambda^d \to d^2\lambda(\vx,\vy).
\end{align*}
Here, $\vy\in Y$ denotes the fast variable, a now \emph{independent}
coordinate for  the microscale. The small parameter, $d$, is treated as
non-dimensional since it is scaled with $1/k_0$ which is set to unity.

\textcolor{black}{The application of the above formal expansions to
\eqref{eq:maxwell} and \eqref{eq:jump} and the subsequent dominant balance
need to be carried out to the two lowest orders in the small parameter,
$d$. This procedure results in eight equations and eight boundary
conditions (\eqref{eq:orderdm1}--\eqref{eq:order0b} in Appendix~\ref{sec:detailed_asymptotics}), as well as two interface
conditions on $\partial\Sigma$
(\eqref{eq:compatibility_rescaled1}--\eqref{eq:compatibility_rescaled} in Appendix~\ref{sec:detailed_asymptotics}). We now use the expanded
equations to describe $\vE^{(0)}$ and $\vH^{(0)}$ (see \eqref{eq:ansatz1}--\eqref{eq:ansatz2} in Appendix~\ref{sec:detailed_asymptotics}; cf.~\cite{maier2018b}). Thus, we find
\begin{equation}
  \label{eq:ansatz}
  \begin{cases}
    \begin{aligned}
      \vE^{(0)}(\vx,\vy) &= \vME(\vx) + \nabla_y\varphi(\vx,\vy)~,
      \quad \varphi(\vx,\vy) = \sum_{j=1}^3\chi_j(\vx,\vy)\mathcal{E}_j(\vx)~,
      \\[-1em]
      \vH^{(0)}(\vx,\vy) &= \vMH(\vx)~,
    \end{aligned}
  \end{cases}
\end{equation}
where $\chi_j(\vx,\,.\,)$ are the correctors and assumed to be $Y$-periodic
($j=1,2,3$).}

\textcolor{black}{%
Next, we focus on the functions $\vME(\vx)$, $\vMH(\vx)$ and $\chi_j(\vx,\vy)$
($j=1,2,3$). For this purpose, we resort to the corresponding
\emph{cell problem} (see Section
\ref{sec:detailed_asymptotics}\ref{subsec:cell_problem} in Appendix~\ref{sec:detailed_asymptotics}). Each $\chi_j(\vx,\vy)$ satisfies}
\begin{align}
  \label{eq:cell_problem}
  \begin{cases}
    \begin{aligned}
      \nabla_y\cdot\Big(\e(\vx,\vy)\big(\vec
      e_j+\nabla_y\chi_j(\vx,\vy)\big)\Big)
      &= 0
      && \text{in }Y\setminus\Sigma~,
      \\[0.5em]
      \js{\vn\times (\vec e_j+\nabla_y\chi_j(\vx,\vy))}
      &= 0
      && \text{on }\Sigma~,
      \\[0.5em]
      \js{\vn\cdot\Big(\e(\vx,\vy)\big(\vec
      e_j+\nabla_y\chi_j(\vx,\vy)\big)\Big)}
      &=
      \frac{1}{i\omega}\nabla_y\cdot\Big(\sigma(\vx,\vy) \big(\vec
      e_j+\nabla_y\chi_j(\vx,\vy)\big)\Big)
      && \text{on }\Sigma~,
      \\[0.5em]
      \vec n\cdot\Big(\sigma(\vx,\vy) \big(\vec
      e_j+\nabla_y\chi_j(\vx,\vy)\big)\Big)
      &=
      \nabla_y\cdot\Big(\lambda(\vx,\vy) \big(\vec
      e_j+\nabla_y\chi_j(\vx,\vy)\big)\Big)
      && \text{on }\partial\Sigma\setminus\partial Y~.
    \end{aligned}
  \end{cases}\hspace{-2.75em}
\end{align}
Equations~\eqref{eq:cell_problem} along with the condition that
$\chi_j(\vx,\,.\,)$ be $Y$-periodic form the desired, closed cell problem.
Notice that $\vx$ plays the role of a parameter. Hence, the above cell
problem uniquely describes $\chi_j(\vx,\,.\,)$ for any given (macroscopic)
point $\vx$.

A remark on the two jump conditions in \eqref{eq:cell_problem} is in order.
The first jump condition ensures that the tangential part of $\vec
e_j+\nabla_y\chi_j(\vx,\vy)$ on $\Sigma$ is single valued. In regard to the
second jump condition, recall that $\sigma(\vx,\vy)$ only acts on the
tangential part of a vector field (to return a tangential vector). Thus,
the right-hand side of the respective jump condition consists of the
divergence of tangential field components, which in turn implies that this
term is also single valued.

\textcolor{black}{Similarly, we recover the \emph{homogenized system}
that fully describes $\vME$ and $\vMH$, viz.,}
\begin{align}
  \begin{cases}
    \begin{aligned}
      \nabla\times\vME &= i\omega\mu_0\vMH~,
      &\quad
      \nabla\times\vMH &= -i\omega\eff\vME+\vJa~;
      \\[0.3em]
      \nabla\cdot(\eff\vME) &= \frac{1}{i\omega}{\nabla\cdot\vJa}~,
      &\quad
      \nabla\cdot\vMH &= 0~,
    \end{aligned}
  \end{cases}
  \label{eq:maxwellhomogenized}
\end{align}
with the \emph{effective permittivity tensor} $\eff$ given by
\begin{multline}
  \eff_{ij}(\vx):=
  \int_Y\e(\vx,\vy)\big(\vec e_j+\nabla_y\chi_j(\vx,\vy)\big)
  \cdot\vec e_i \dy
  -\frac1{i\omega}
  \int_\Sigma
  \sigma(\vx, \vy)\big(\vec e_j+\nabla_y\chi_j(\vx,\vy)\big) \cdot \vec
  e_i\doy
  \\
  -\frac1{i\omega}
  \int_{\partial\Sigma\setminus\partial Y}
  \lambda(\vx, \vy)\big(\vec e_j+\nabla_y\chi_j(\vx,\vy)\big) \cdot \vec
  e_i\,\text{d} s~.
  \label{eq:effectiveperm}
\end{multline}
\textcolor{black}{For the sake of brevity, we give the derivation of
this formula to Section
\ref{sec:detailed_asymptotics}\ref{subsec:homogenized_problem} of
Appendix~\ref{sec:detailed_asymptotics}.}

\textcolor{black}{Equations~\eqref{eq:ansatz}--\eqref{eq:effectiveperm}
summarize our main homogenization results. In conclusion,
system~\eqref{eq:maxwellhomogenized} describes the large-scale optical
response of the periodic medium in terms of the macroscopic electromagnetic
field $(\vME, \vMH)$ via the effective permittivity tensor $\eff$ given
by~\eqref{eq:effectiveperm} with~\eqref{eq:cell_problem}.}

\section{Epsilon-near-zero effect}
\label{sec:enz}

A crucial feature of the averaging for the effective permittivity,
$\eff(\vx)$, as implied by~\eqref{eq:effectiveperm}, is the
\emph{interplay} of three distinct averages: one for the bulk permittivity,
$\varepsilon$; another for the surface conductivity, $\sigma$; and a third
one for the line conductivity, $\lambda$. Each of these averages has a
positive real part. This interplay can be exploited as follows. By tuning
the microscopic geometry, periodic spacing, frequency, or surface
or line conductivity, we can in principle force at least one eigenvalue of
$\eff(\vx)$ to be close to zero. In the case with an $\vx$-independent
$\eff$, this condition amounts to the ENZ effect in the direction of
propagation determined by the respective
eigenvector(s)~\cite{maier2018a,mattheakis2016}.

In the remainder of this section, we discuss the simplified case with a
vanishing corrector, $\vec\chi$. In this case, we derive an explicit
formula for the special spacing, $d=d_c$, called the critical spacing, that
implies the ENZ effect in the presence of the line conductivity, $\lambda$.
We show how our present averaging procedure is related to results of a
Bloch-wave approach~\cite{maier2018a}. Lastly, we discuss in some detail
the character and physical role of the corrector, $\vec\chi$.

\subsection{Case with vanishing corrector $\vec\chi$: Formalism}
\label{subsec:vanishing}
\color{black} In principle, \color{black} the corrector, $\vec\chi$,
implicitly depends on the geometry and material parameters such as
$\sigma$ and $\lambda$. Accordingly, the derivation of closed, analytical
formulas for $\vec\chi$ may be possible only in a limited number of
situations. We now restrict attention to the particular yet physically
appealing case with a vanishing corrector, $\vec\chi$. Our goal is to
better understand the emerging ENZ effect.

Technically speaking, the solution, $\vec\chi$, of cell
problem~\eqref{eq:cell_problem} vanishes whenever the forcing
term in this problem is identically zero. This term may not be immediately
obvious by inspection of~\eqref{eq:cell_problem}. In order to gain some
insight into its character, we invoke the weak formulation of the cell
problem. To this end, we multiply the first equation
in~\eqref{eq:cell_problem} with a smooth and $Y$-periodic test function,
$\psi$, and integrate over $Y$ (see also Section~\ref{sec:weak_formulation}
of the supplementary material). Two subsequent integrations by parts, and
use of the requisite jump and boundary conditions, yield
\begin{multline*}
  0 \;=\;
  \int_Y \e(\vx,\vy)\big(\vec e_j+\nabla_y\chi_j(\vx,\vy)
  \big)\cdot\nabla_y\psi(\vy)\,\dy
  \\
  -
  \frac{1}{i\omega} \int_{\Sigma}
  \sigma(\vx,\vy)\big(\vec
  e_j+\nabla_y\chi_j(\vx,\vy)\big)\cdot\nabla_y\psi(\vy)\,\doy
  \\
  -
  \frac{1}{i\omega} \int_{\partial\Sigma\setminus\partial Y}
  \lambda(\vx,\vy)\big(\vec
  e_j+\nabla_y\chi_j(\vx,\vy)\big)\cdot\nabla_y\psi(\vy)\,\text{d}s~.
\end{multline*}
The forcing term corresponds to the following contribution:
\begin{align}
  \int_Y \e(\vx,\vy)\,\vec e_j\cdot\nabla_y\psi(\vy)\,\dy
  -
  \frac{1}{i\omega} \int_{\Sigma}
  \sigma(\vx,\vy)\,\vec e_j\cdot\nabla_y\psi(\vy)\,\doy
  -
  \frac{1}{i\omega} \int_{\partial\Sigma\setminus\partial Y}
  \lambda(\vx,\vy)\,\vec e_j\cdot\nabla_y\psi(\vy)\,\text{d}s~.
  \label{eq:forcing}
\end{align}
Thus, the requirement of having a vanishing corrector means
that~\eqref{eq:forcing} must vanish identically for every possible choice
of $\psi$. One more integration by parts converts
expression~\eqref{eq:forcing} to
\enlargethispage{1em}
\begin{multline*}
  i\omega\int_Y \nabla_y\cdot\big(\e(\vx,\vy)\vec e_j\big)\psi(\vy)\,\dy
  + \int_{\partial\Sigma} \vec n\cdot\big(\sigma(\vx,\vy)\vec
  e_j\big)\psi(\vy)\,\text{d}s
  \\
  -
  \int_{\Sigma}
  \nabla_y\cdot\big(\sigma(\vx,\vy)\vec e_j\big)\psi(\vy)\,\doy
  -
  \int_{\partial\Sigma\setminus\partial Y}
  \nabla_y\cdot\big(\lambda(\vx,\vy)\vec e_j\big)\psi(\vy)\,\text{d}s~.
\end{multline*}
Thus, the forcing is necessarily zero whenever the following conditions
hold simultaneously:
\begin{align}
  \label{eq:divergence_free}
  \begin{cases}
    \begin{aligned}
      \nabla_y\cdot\e(\vx,\vy)\;&=\;0
      \quad\text{in }Y~,
      \\[0.1em]
      \nabla_y\cdot\sigma(\vx,\vy)\;&=\;0
      \quad\text{on }\Sigma~,
      \\[0.1em]
      \vec n\cdot\sigma(\vx,\vy)
      -\nabla_y\cdot\lambda(\vx,\vy)
      \;&=\;0
      \quad\text{on }\partial\Sigma\setminus\partial Y~.
    \end{aligned}
  \end{cases}
\end{align}
Note that $\sigma$ and $\lambda$ encode information about the curvature of
$\Sigma$ and $\partial\Sigma$, respectively. This implies that for a
non-flat surface $\Sigma$, or for a surface with (arbitrarily shaped)
internal edges, $\sigma(\vx,\vy)$ and $\lambda(\vx,\vy)$ are in principle
$\vy$-dependent with non-vanishing divergence. Under the assumption
that~\eqref{eq:divergence_free} holds, formula~\eqref{eq:effectiveperm} for
the effective permittivity simply reduces to
\begin{align}
  \label{eq:simplified_eff}
  \eff= \int_Y\e(\vx,\vy)\dy - \frac1{i\omega} \int_\Sigma \sigma(\vx,
  \vy)\doy
  - \frac1{i\omega} \int_{\partial\Sigma\setminus\partial Y} \lambda(\vx,
    \vy)\,\text{d}s~.
\end{align}
By this simplified (geometric) average, one may directly influence the
effective permittivity tensor by either adjusting the operating frequency,
$\omega$, or by tuning the parameters $\sigma(\vx,\vy)$ and
$\lambda(\vx,\vy)$.

We add a remark on the above formalism. A natural question at this point
concerns the existence of suitable configurations that
obey~\eqref{eq:divergence_free}. More precisely, it is of interest to
specify configurations that allow for a vanishing corrector while two or
all three of the distinct averages in \eqref{eq:simplified_eff} are
nonzero. \color{black} We give three related examples in this section.
\color{black}  First, consider the geometry with parallel, planar sheets of
constant $\sigma$ and no edges, $\partial\Sigma=\emptyset$, embedded in a
homogeneous dielectric host with uniform $\varepsilon$; see
Figure~\ref{fig:cell_problems}a. In this case, $\vec\chi\equiv 0$, and we
can have a nonzero contribution from the average of $\sigma$
in~\eqref{eq:simplified_eff}. The second geometry consists of parallel and
periodically aligned nanoribbons embedded in a homogeneous, dielectric
host; see Figure~\ref{fig:cell_problems}b. In this case, we can achieve
$\vec\chi\equiv 0$ through \eqref{eq:divergence_free} and a nonzero
contribution of the line average of $\lambda$ in \eqref{eq:simplified_eff}
by setting $\sigma^d=0$ and choosing a constant $\lambda^d$. Further, in
the same geometry, by choosing a constant $\sigma$ and the $y_3$-dependent
line conductivity $\lambda=\sigma\,y_3$, all three averages
in~\eqref{eq:simplified_eff} are nonzero while condition
\eqref{eq:divergence_free} holds; thus, $\vec\chi\equiv 0$. An example with
a geometry that always has a nonzero corrector (in the presence of a
nonzero conductivity $\sigma$) consists of parallel and periodically
aligned nanotubes; see Figure~\ref{fig:cell_problems}c. \color{black}
For this third case, the corrector $\vec\chi$ is obtained numerically. (For
further discussion on the computational framework for $\vec\chi$, the
interested reader is \textcolor{black}{r}eferred to
Section~\ref{sec:numerics}.)
 \color{black}

\subsection{Vanishing corrector: A notion of critical spacing}
In this subsection, we further simplify average~\eqref{eq:simplified_eff}
in order to derive an explicit formula for the critical spacing, $d_c$.
Recall the ansatz
\begin{align*}
  \sigma^d(\vx) = d\,\sigma(\vx,\vx/d),
  \quad\lambda^d(\vx) = d^2\,\lambda(\vx,\vx/d)~;
\end{align*}
see Section~\ref{sec:intro}\,\ref{subsec:scaling}. Although this choice of
scaling is convenient for the asymptotic analysis
(Section~\ref{sec:asymptotics}), it introduces two parameters, $\sigma$ and
$\lambda$, that couple the physical conductivities, $\sigma^d$ and
$\lambda^d$, with the spacing, $d$. This coupling may obscure the physical
insight possibly gained by the averaging in~\eqref{eq:simplified_eff}. A
reason is that, in realistic settings, the parameters $d$, $\sigma^d$ and
$\lambda^d$ are controlled independently by tuning of the geometry and the
electronic structure of the 2D material~\cite{mattheakis2016}.

As a remedy, consider the following rescaled quantities
\begin{align*}
  \sigma^d(\vx,\vy) := d\,\sigma(\vx,\vy)~,
  \quad\lambda^d(\vx,\vy) := d^2\,\lambda(\vx,\vy)~,
\end{align*}
where the microscale variable, $y$, is singled out and treated as
independent in the actual conductivities, $\sigma^d(\vx)$ and
$\lambda^d(\vx)$. The averages of interest therefore are
\begin{align*}
  \bar\e\,(\vx) = \int_Y\e\,(\vx,\vy)\dy~,
  \qquad
  \bar\sigma^d(\vx) = \int_\Sigma \sigma^d(\vx,\vy)\doy~,
  \qquad
  \bar\lambda^d(\vx) = \int_{\partial\Sigma\setminus\partial
  Y}\lambda^d(\vx,\vy)\,\text{d}s~.
\end{align*}
For slowly varying parameters, these averages simply reduce to the original
parameters. We assume that $\bar\e\,(\vx)$, $\bar\sigma^d(\vx)$ and
$\bar\lambda^d(\vx)$ have the same eigenvectors, or \emph{principal axes},
$\vec r_i$ ($i=1,2,3$). Let $\bar\e_i$, $\bar\sigma^d_i$, and
$\bar\lambda_i^d$, denote the corresponding eigenvalues of the averaged
tensors. (Only two eigenvalues of $\bar\sigma^d$ and one eigenvalue of
$\bar\lambda^d$ are possibly nonzero.) After some
algebra,~\eqref{eq:simplified_eff} gives
\begin{align}
  \frac{\eff_i}{\bar\e_i}
  \;=\; \Big(1-\frac{\xi_{0,i}}d\Big)
  \Big(1+\frac{\bar\lambda^d_i}{i\omega\bar\e_i\,\xi_{0,i}\,d}\Big)~,
  \label{eq:eff_plasmonic}
\end{align}
for all directions $\vec r_i$ with nonzero $\bar\sigma^d_i$, or
$\bar\lambda^d_i$. Here, we define the \emph{generalized plasmonic
thickness}~\cite{maier2018a}
\begin{align}
  \xi_{0,i}\;=\;
  \frac{\bar\sigma^d_i}{2\,i\omega\bar\e_i}
  +
  \sqrt{ \left(\frac{\bar\sigma^d_i}{2\,i\omega\bar\e_i}\right)^2
  + \frac{\bar\lambda^d_i}{i\omega\bar\e_i}}~,
  \label{eq:plasmonic_thickness}
\end{align}
which accounts for the line conductivity. Evidently, $\eff_i\simeq 0$ if
the spacing, $d$, is close to the critical value $d_c= \xi_{0,i}$ for some
$i$ and suitable ranges of values for $\bar\sigma^d_i$ and
$\bar\lambda^d_i$ (e.g., ${\rm Im}\,\bar\sigma^d_i>0, {\rm
Im}\,\bar\lambda^d_i>0$).

We stress that for $|\bar\sigma^d_i|\gg|\bar\lambda^d_i|$
\eqref{eq:eff_plasmonic} with~\eqref{eq:plasmonic_thickness} reduces to the
known result
\begin{align*}
  \frac{\eff_i}{\bar\e_i} \;\sim\; \Big(1-\frac{\xi_{0,i}}d\Big)~,
  \qquad
  \xi_{0,i}\;\approx\; \frac{\bar\sigma^d_i}{i\omega\bar\e_i}~,
\end{align*}
which is also obtained via a Bloch-wave approach for planar
sheets~\cite{maier2018a}. In contrast, if
$|\bar\lambda^d_i|\gg|\bar\sigma^d_i|$, \eqref{eq:eff_plasmonic} reduces to
\begin{align*}
  \frac{\eff_i}{\bar\e_i} \;\sim\; \Big(1-\frac{\xi_{0,i}}d\Big)~,
  \qquad
  \xi_{0,i}\;\approx\; \sqrt{\frac{\bar\lambda^d_i}{i\omega\bar\e_i}}.
\end{align*}
Notice the distinctly different scaling of $\xi_{0,i}$ with the
conductivity parameter in this regime.
%

\subsection{On the physical role of corrector $\vec\chi$}
\label{subsec:corrector}

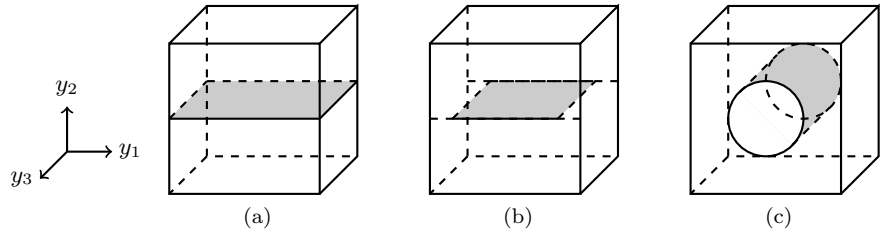
\begin{figure}
  \begin{center}
    \begin{tikzpicture}[scale=1.20]
        \path[thick,draw,->] (-2.5,-0.6) -- (-2.0,-0.6);
        \path[thick,draw,->] (-2.5,-0.6) -- (-2.5,-0.1);
        \path[thick,draw,->] (-2.5,-0.6) -- (-2.8,-0.9);
        \node at (-2.5, +0.1) {\small $y_2$};
        \node at (-1.8, -0.6) {\small $y_1$};
        \node at (-3.0, -0.9) {\small $y_3$};
      \end{tikzpicture}
    \subfloat[]{
      \begin{tikzpicture}[scale=1.0]
        \path [thick, draw]
          (-1, -1) -- (1, -1) -- (1, 1) -- (-1, 1) -- cycle;
        \path [thick, draw] (1,-1) -- (1.5, -0.5) -- (1.5, 1.5) -- (1, 1);
        \path [thick, draw] (-1,1) -- (-0.5, 1.5) -- (1.5, 1.5);
        \path [thick, draw, dashed]
          (-0.5, 1.5) -- (-0.5, -0.5) -- (1.5, -0.5);
        \path [thick, draw, dashed] (-0.5, -0.5) -- (-1.0, -1.0);
        \path [fill=black, fill opacity=0.2]
          (-1, 0) -- (1, 0) -- (1.5, 0.5) -- (-0.5, 0.5) -- (-1.0, 0.0);
        \path [thick, draw] (-1, 0) -- (1, 0);
        \path [thick, draw] (1, 0) -- (1.5, 0.5);
        \path [thick, draw, dashed] (-0.5, 0.5) -- (1.5, 0.5);
        \path [thick, draw, dashed] (-1.0, 0.0) -- (-0.5, 0.5);
      \end{tikzpicture}}
    \hspace{2em}
    \subfloat[]{
      \begin{tikzpicture}[scale=1.0]
        \path [thick, draw]
          (-1, -1) -- (1, -1) -- (1, 1) -- (-1, 1) -- cycle;
        \path [thick, draw] (1,-1) -- (1.5, -0.5) -- (1.5, 1.5) -- (1, 1);
        \path [thick, draw] (-1,1) -- (-0.5, 1.5) -- (1.5, 1.5);
        \path [thick, draw, dashed]
          (-0.5, 1.5) -- (-0.5, -0.5) -- (1.5, -0.5);
        \path [thick, draw, dashed] (-0.5, -0.5) -- (-1.0, -1.0);
        \path [fill=black, fill opacity=0.2]
          (-0.7, 0) -- (0.7, 0) -- (1.2, 0.5) -- (-0.2, 0.5) -- (-0.7, 0.0);
        \path [thick, draw] (-0.7, 0) -- (0.7, 0);
        \path [thick, draw, dashed] (0.7, 0) -- (1.2, 0.5);
        \path [thick, draw, dashed] (-0.2, 0.5) -- (1.2, 0.5);
        \path [thick, draw, dashed] (-0.7, 0.0) -- (-0.2, 0.5);
        \path [thick, draw, dashed] (-1, 0) -- (1, 0);
        \path [thick, draw, dashed] (-0.5, 0.5) -- (1.5, 0.5);
      \end{tikzpicture}}
    \hspace{2em}
    \subfloat[]{
      \begin{tikzpicture}[scale=1.0]
        \path [thick, draw]
          (-1, -1) -- (1, -1) -- (1, 1) -- (-1, 1) -- cycle;
        \path [thick, draw] (1,-1) -- (1.5, -0.5) -- (1.5, 1.5) -- (1, 1);
        \path [thick, draw] (-1,1) -- (-0.5, 1.5) -- (1.5, 1.5);
        \path [thick, draw, dashed]
          (-0.5, 1.5) -- (-0.5, -0.5) -- (1.5, -0.5);
        \path [thick, draw, dashed] (-0.5, -0.5) -- (-1.0, -1.0);

        \path [fill=black, fill opacity=0.2]
        (0.85, 0.15) arc (-45:135:0.5) -- (-0.35, 0.35) -- (0.35, -0.35);
        \path [fill=white] (0.35, -0.35) arc (-45:135:0.5);
        \path [thick, draw, dashed] (0.5, 0.5) circle (0.5);
        \path [thick, draw] (0.0, 0.0) circle (0.5);
        \path [thick, draw, dashed] (0.35, -0.35) -- (0.85, 0.15);
        \path [thick, draw, dashed] (-0.35, 0.35) -- (0.15, 0.85);
      \end{tikzpicture}}
  \end{center}
  \caption{
    Prototypical examples of microscopic geometries with conducting sheets
    for cell problem~\eqref{eq:cell_problem}. (a) Infinite planar sheet,
    with no edges. (b) Planar strip (nanoribbon). (c) Sheet forming
    circular cylinder (nanotube). The corrector, $\vec\chi$, can be
    characterized as follows: (a) $\vec\chi\equiv 0$; (b)
    $\chi_2=\chi_3\equiv 0$ while $\chi_1$ is nontrivial; and (c)
    $\chi_3\equiv 0$ while $\chi_1$ and $\chi_2$ are nontrivial.}
  \label{fig:cell_problems}
\end{figure}

\begin{figure}
  \begin{center}
    \subfloat[$\text{Re}\,\chi_1$, ribbons]{
      \includegraphics[width=2.8cm]{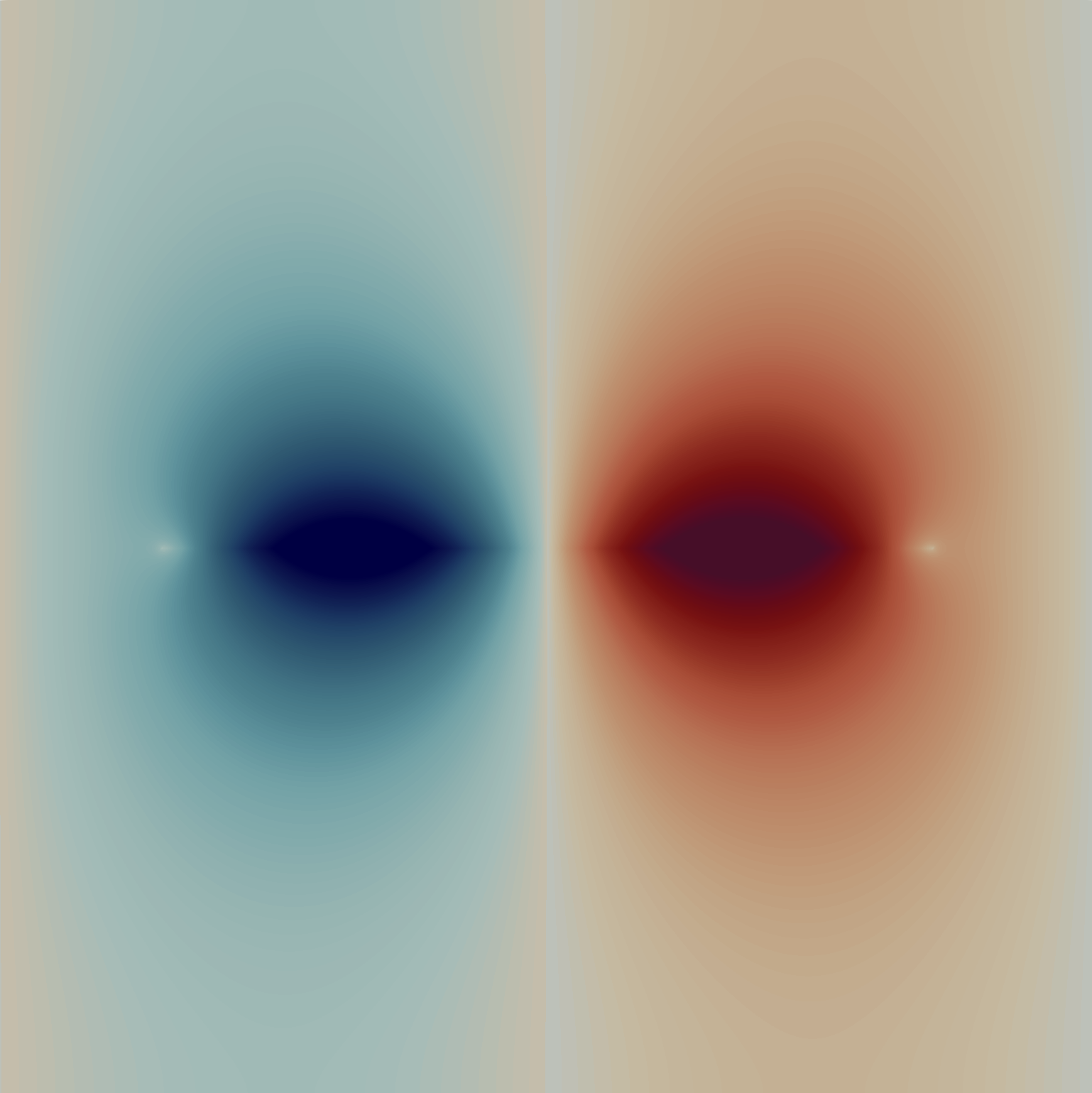}
    }
    \hspace{0.1em}
    \subfloat[$\text{Im}\,\chi_1$, ribbons]{
      \includegraphics[width=2.8cm]{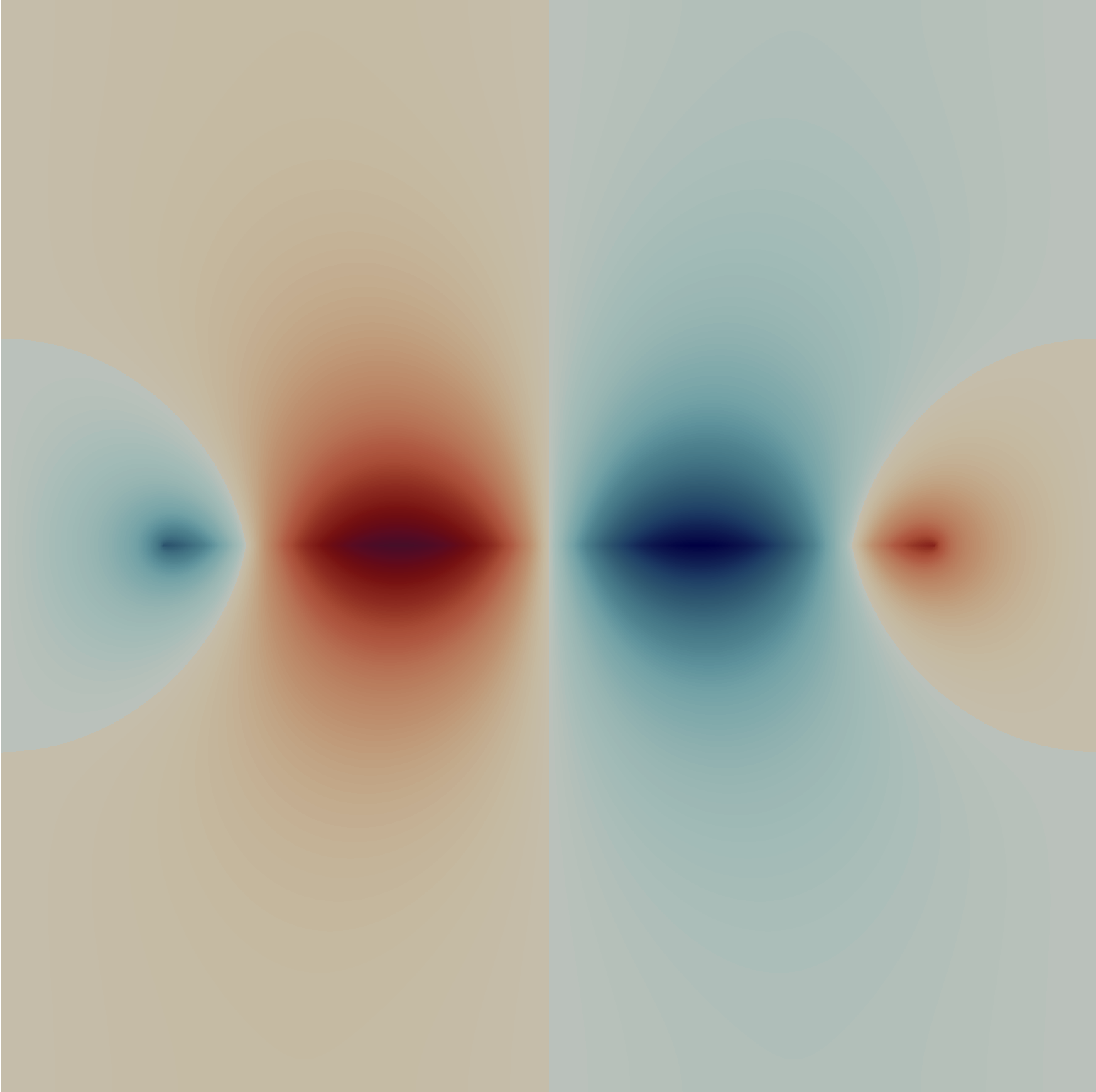}
    }
    \hspace{0.1em}
    \subfloat[$\text{Re}\,\chi_1$, tubes]{
      \includegraphics[width=2.8cm]{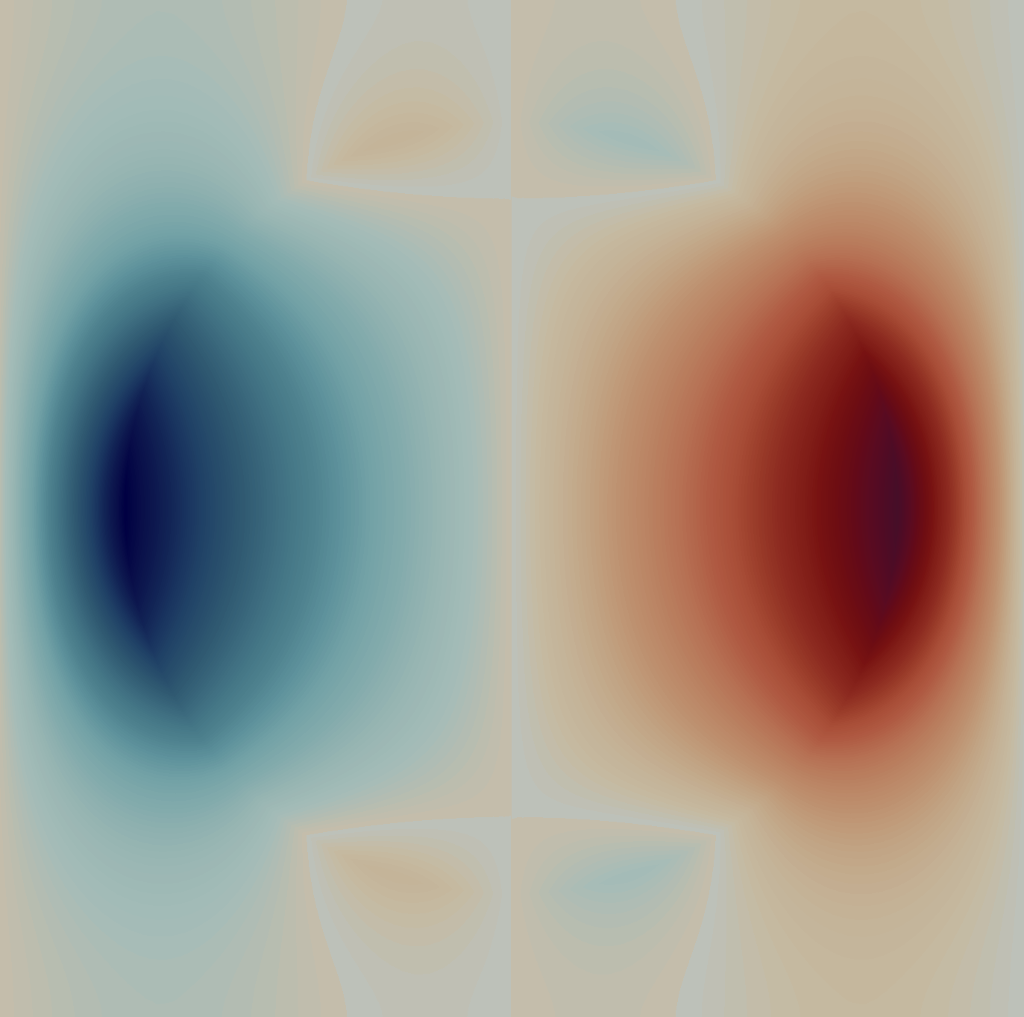}
    }
    \hspace{0.1em}
    \subfloat[$\text{Im}\,\chi_1$, tubes]{
      \includegraphics[width=2.8cm]{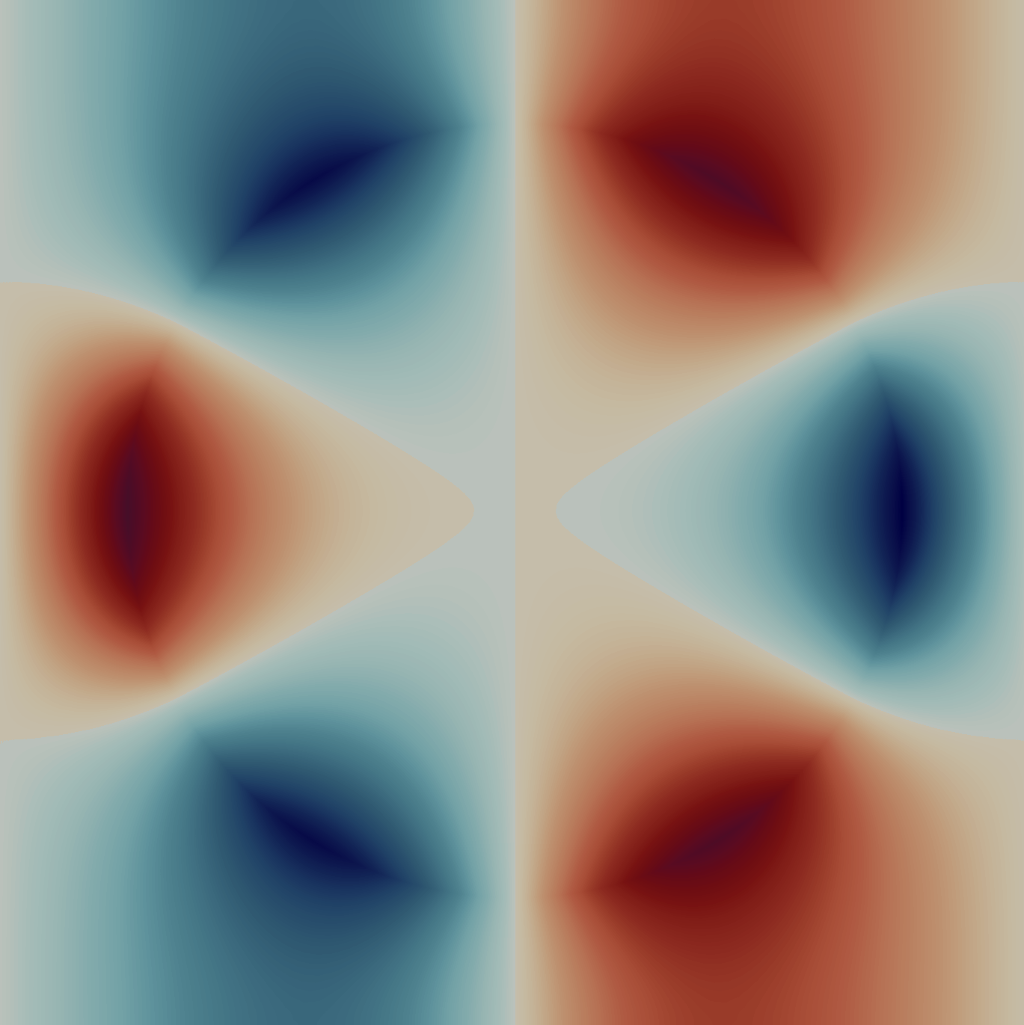}
      \begin{tikzpicture}
        \node at (0.0, 0.0)
          {\includegraphics[width=0.3cm]{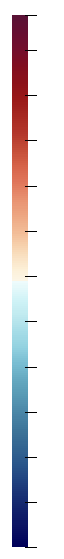}};
        \node at (0.3,  1.25) {\small $+$};
        \node at (0.3, -1.25) {\small $-$};
      \end{tikzpicture}
      \hspace{-3em}
    }
  \end{center}
  \caption{
    Real and imaginary parts of corrector component $\chi_1$ for geometries
    of Figures~\ref{fig:cell_problems}b,c in $y_1y_2$-plane. It is evident
    that internal edges in nanoribbons and curvature in nanotubes create
    SPPs.}
  \label{fig:corrector_response}
\end{figure}

Let us return to cell problem~\eqref{eq:cell_problem}. We now provide a
physical interpretation of its solution, the corrector $\vec\chi$. In
particular, we solve cell problem~\eqref{eq:cell_problem} for the three
prototypical configurations shown in Figure~\ref{fig:cell_problems} (see
also Section~\ref{sec:enz}\,\ref{subsec:vanishing}). In all three
geometries, we set the permittivity and surface conductivity equal to
nonzero constants, $\e(\vx,\vy)=\e$, $\sigma(\vx,\vy)=\sigma$, while we
take $\lambda(\vx,\vy)=0$. By the discussion in
Section~\ref{sec:enz}\,\ref{subsec:vanishing}, we conclude that for the
geometry of planar sheets with no edges (Figure~\ref{fig:cell_problems}(a))
$\vec\chi\equiv 0$. For the geometry of nanoribbons
(Figure~\ref{fig:cell_problems}(b)), $\chi_2\equiv \chi_3\equiv 0$. The
remaining geometry (Figure~\ref{fig:cell_problems}(c)) has only one
vanishing corrector component, $\chi_3\equiv 0$. The real and imaginary
parts of the nontrivial corrector component $\chi_1$ for the last two
geometries are shown in Figure~\ref{fig:corrector_response}. The SPP
excited by the edge discontinuity is evident in these settings;
cf.~\cite{maier2017}.

Motivated by these numerical results, we develop an argument that the
corrector, $\vec\chi$, encodes the microscale response of the system to the
macroscopic electromagnetic field, $(\vME(\vx),\vMH(\vx))$. We start by
noting that our scalings $\sigma^d(\vx)=\,d\,\sigma(\vx,\vx/d)$ and
$\lambda^d(\vx)=\,d^2\,\lambda(\vx,\vx/d)$ imply that the typical length
scales of surface and line waves (such as SPPs and EPPs) scale with $d$ and
$d^2$, respectively. Accordingly, we can identify the forcing in cell
problem \eqref{eq:cell_problem} due to $\vec e_j$. This forcing corresponds to a
(normalized) asymptotically slow planar wave. In this sense, the cell
problem describes the \emph{local response of the microstructure} to all
possible excitations by local plane waves.

This interpretation has an important consequence in light of the discussion
in Section~\ref{sec:enz}\,\ref{subsec:vanishing}: All cases of vanishing
correctors are indeed characterized by conditions~\eqref{eq:divergence_free}.
In fact, conditions \eqref{eq:divergence_free} characterize exactly all
microscale geometries that \emph{do not permit} the excitation of SPPs by
plane waves.

\section{\textcolor{black}{Computational platform and examples of effective permittivity}}
\label{sec:numerics}

\textcolor{black}{%
In this section, we introduce a computational platform that serves as a
foundation for investigating the effective optical response of
(sufficiently general) microscopic geometries.
\textcolor{black}{We implement the computational framework with the
help of the finite element toolkit deal.II~\cite{dealii85}
and demonstrate its effectiveness on} a prototypical and somewhat realistic
example of a plasmonic crystal consisting of corrugated graphene sheets. In
addition, we} present a number of computational results concerning the
frequency response of effective permittivity~\eqref{eq:effectiveperm} by
use of cell problem~\eqref{eq:cell_problem} \color{black} for the
geometries of nanoribbons and nanotubes
(Figure~\ref{fig:cell_problems}b,c). \color{black}

\subsection{\textcolor{black}{Computational framework}}
\label{subsec:framework}

\textcolor{black}{%
Building on computational methods that we developed for plasmonic
problems~\cite{maier2017a, maier2017}, we propose the following
computational approach for present purposes:
\vspace{-1em}
\begin{itemize}
  \item
    We compute approximations for solution $\{\vME,\vMH\}$ of
    homogenized problem \eqref{eq:maxwellhomogenized} by use of an
    adaptive finite element scheme. Typically, for such a scheme to be
    efficient there is a need for advanced computational
    techniques\textcolor{black}{,} such as the construction of a
    \emph{perfectly matched layer} for treatment of absorbing boundary
    conditions\textcolor{black}{, and strategies for} adaptive local
    refinement of meshes. We refer the reader to \cite{maier2017a} for
    details.
  \item
    \textcolor{black}{%
    During the finite element computation of \eqref{eq:maxwellhomogenized}
    we \emph{reconstruct} the effective permittivity tensor $\eff$ when
    needed by first approximating the corrector $\chi_i$ ($i=1,2,3$) that
    is described by cell problem \eqref{eq:cell_problem}, and subsequently
    evaluating \eqref{eq:effectiveperm} with a suitable numerical
    quadrature rule. This is done with the same finite element toolkit
    \cite{dealii85} that we use for effective problem
  \eqref{eq:maxwellhomogenized}.}
\end{itemize}}
\vspace{-1em}
\textcolor{black}{%
We make the remark that using a finite element discretization is
particularly advantageous for approximating cell problem
\eqref{eq:cell_problem} that contains a jump condition over smooth, curved
hypersurfaces. For the sake of brevity we omit algorithmic details but
refer the reader to \cite{maier2018b} for a discussion of the variational
form of equations \eqref{eq:maxwellhomogenized} and
\eqref{eq:cell_problem}, as well as to \cite{dealii85} for algorithmic
details on curved boundary approximations, numerical quadrature, and
numerical linear algebra.}

\begin{figure}[t]
  \centering
  \subfloat[]{%
    \begin{tikzpicture}[xscale=1.00,yscale=0.943]
      \def\nlayers{8}
      \def\nrepeti{19}
      \def\delta{1.5 / 4. / \nlayers};
      \node at (1.5, 1.5) {$\Omega$};
      \foreach \n in {0,...,\nrepeti} {
        \def\x{-1.875 + 4. * \delta * \n};
        \def\yu{-0.75};
        \def\yo{ 0.75};
        \path[thick, draw=black!10, fill=black!10]
                     (\x,               \yu)
          sin        (\x + \delta,      \yu + \delta)
          cos        (\x + 2. * \delta, \yu)
          sin        (\x + 3. * \delta, \yu - \delta)
          cos        (\x + 4. * \delta, \yu)
          --         (\x + 4. * \delta, \yo)
          sin        (\x + 3. * \delta, \yo - \delta)
          cos        (\x + 2. * \delta, \yo)
          sin        (\x + \delta,      \yo + \delta)
          cos        (\x,               \yo)
          --         cycle;
      }
      \foreach \n in {0,...,\nrepeti} {
        \def\x{-1.875 + 4. * \delta * \n};
        \foreach \m in {0,...,\nlayers} {
          \def\y{ -0.75 + 4. * \delta * \m};
          \draw[thick] (\x,         \y)
            sin        (\x + \delta,      \y + \delta)
            cos        (\x + 2. * \delta, \y)
            sin        (\x + 3. * \delta, \y - \delta)
            cos        (\x + 4. * \delta, \y);
        }
      }
      \path [thick, draw]
        (-1.875, -2) -- (1.875, -2) -- (1.875, 2) -- (-1.875, 2) -- cycle;
      \path [line width=0.5mm, ->, draw] (0, -1.40) -- (0, -1.05);
      \node at (0.5, -1.25) {$\vJa$};
    \end{tikzpicture}
  }
  \hspace{0.3em}
  \subfloat[]{\includegraphics[height=3.8cm]{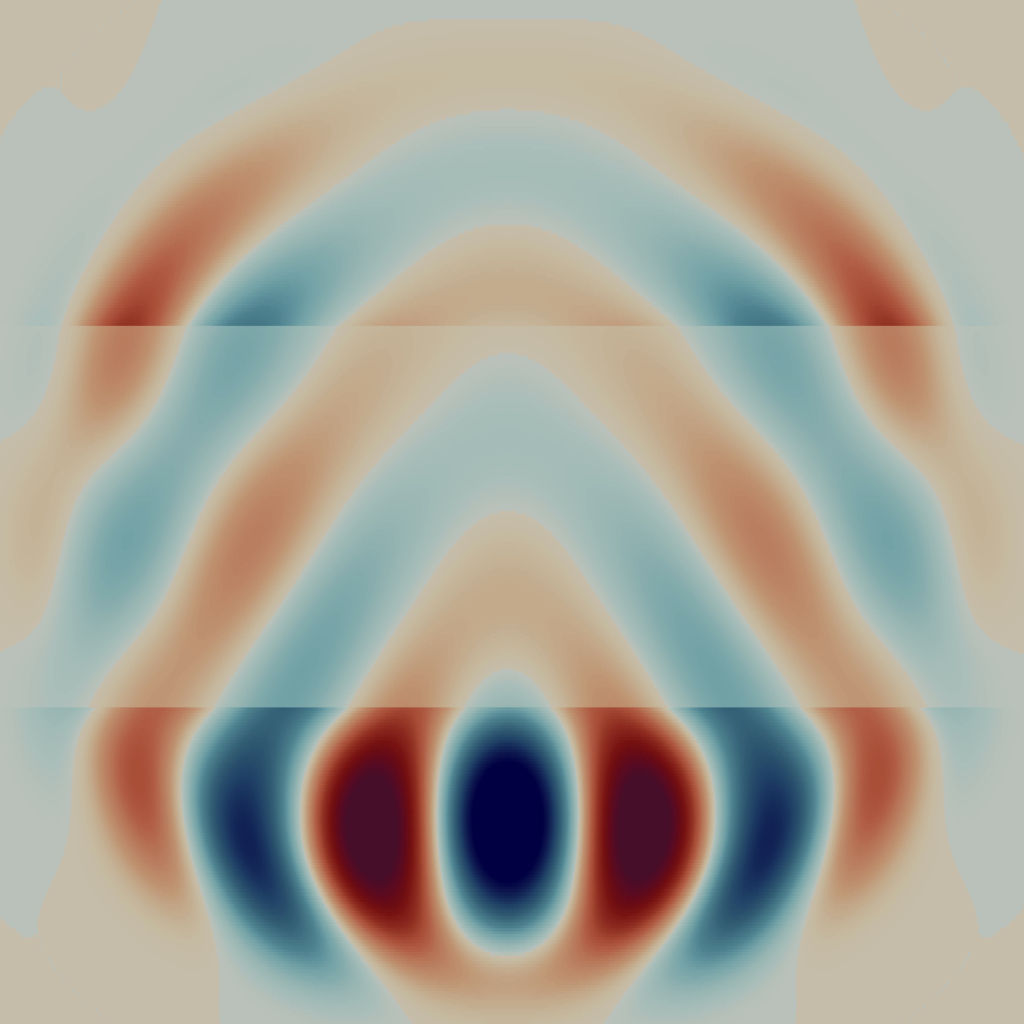}}

  \subfloat[]{\includegraphics[height=3.8cm]{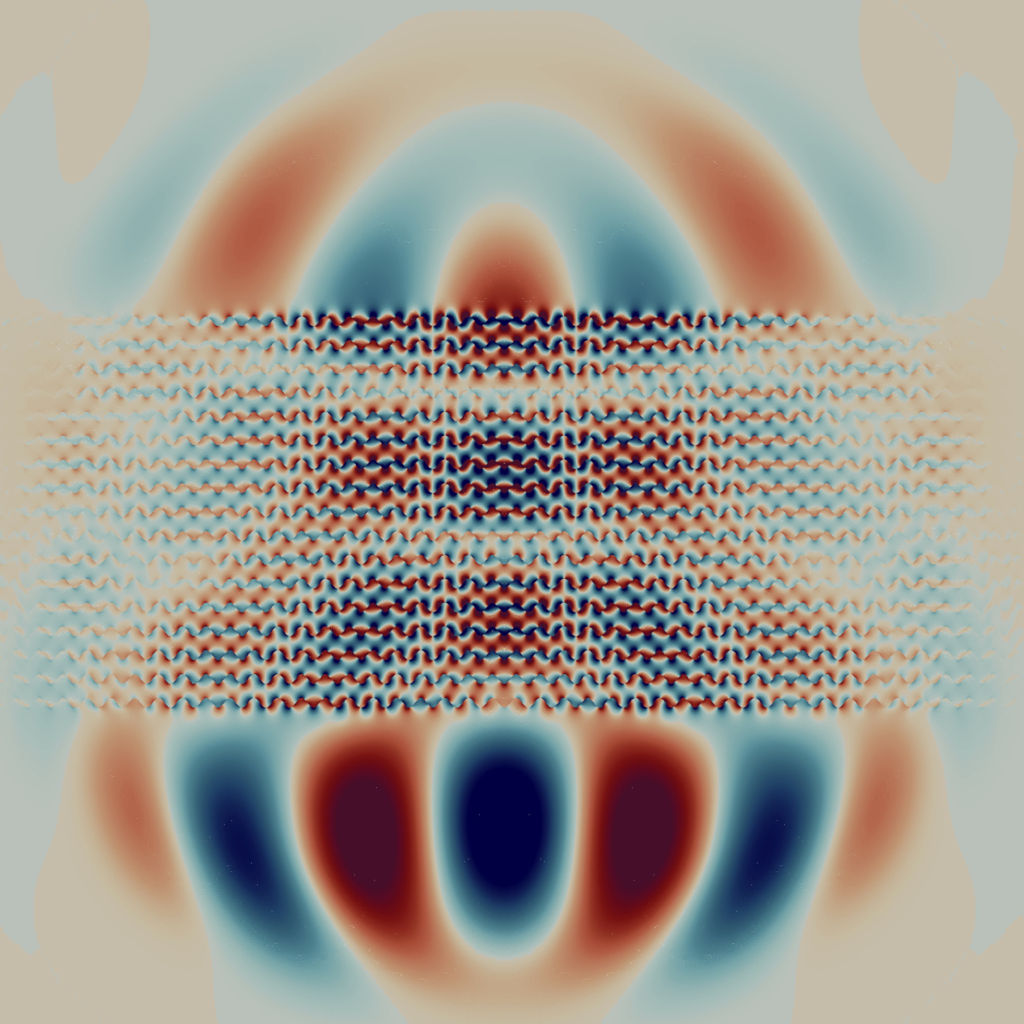}}
  \hspace{0.3em}
  \subfloat[]{\includegraphics[height=3.8cm]{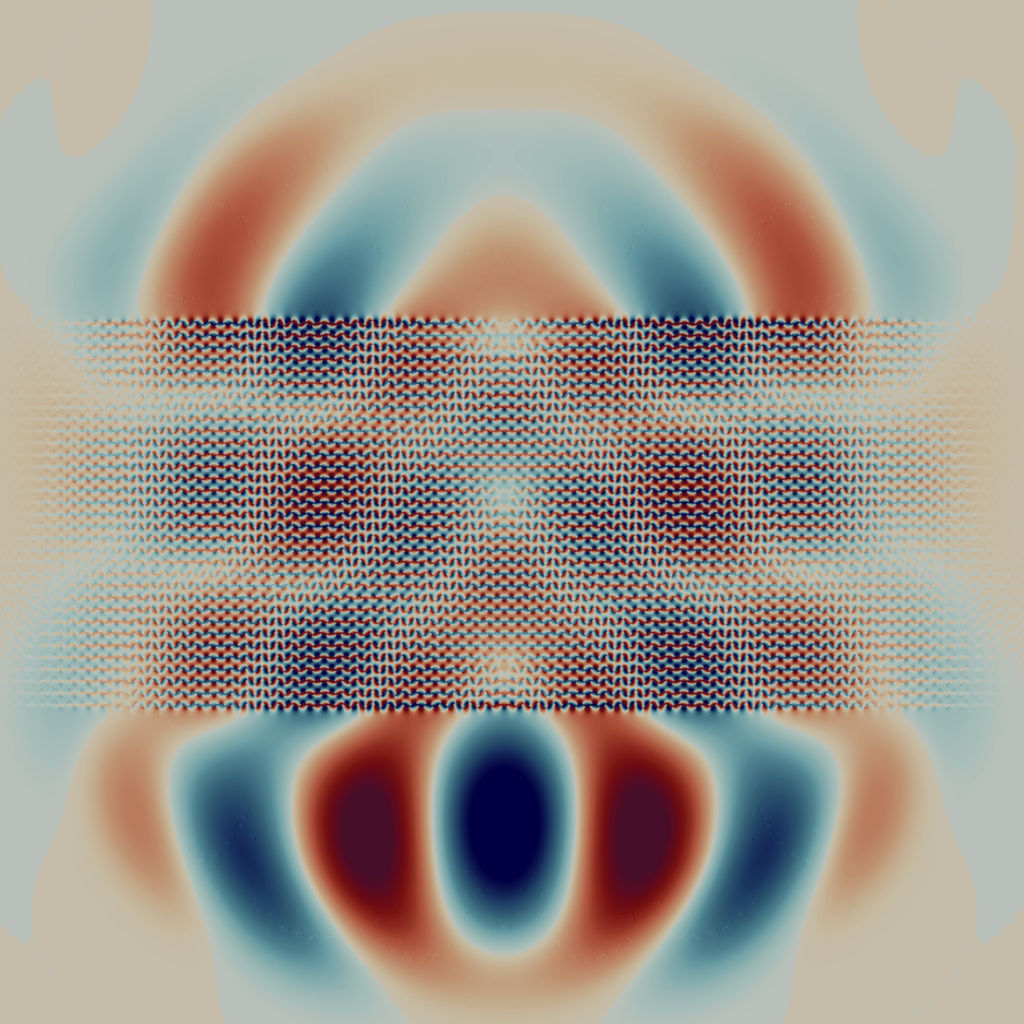}}
  \hspace{0.3em}
  \subfloat[]{\includegraphics[height=3.8cm]{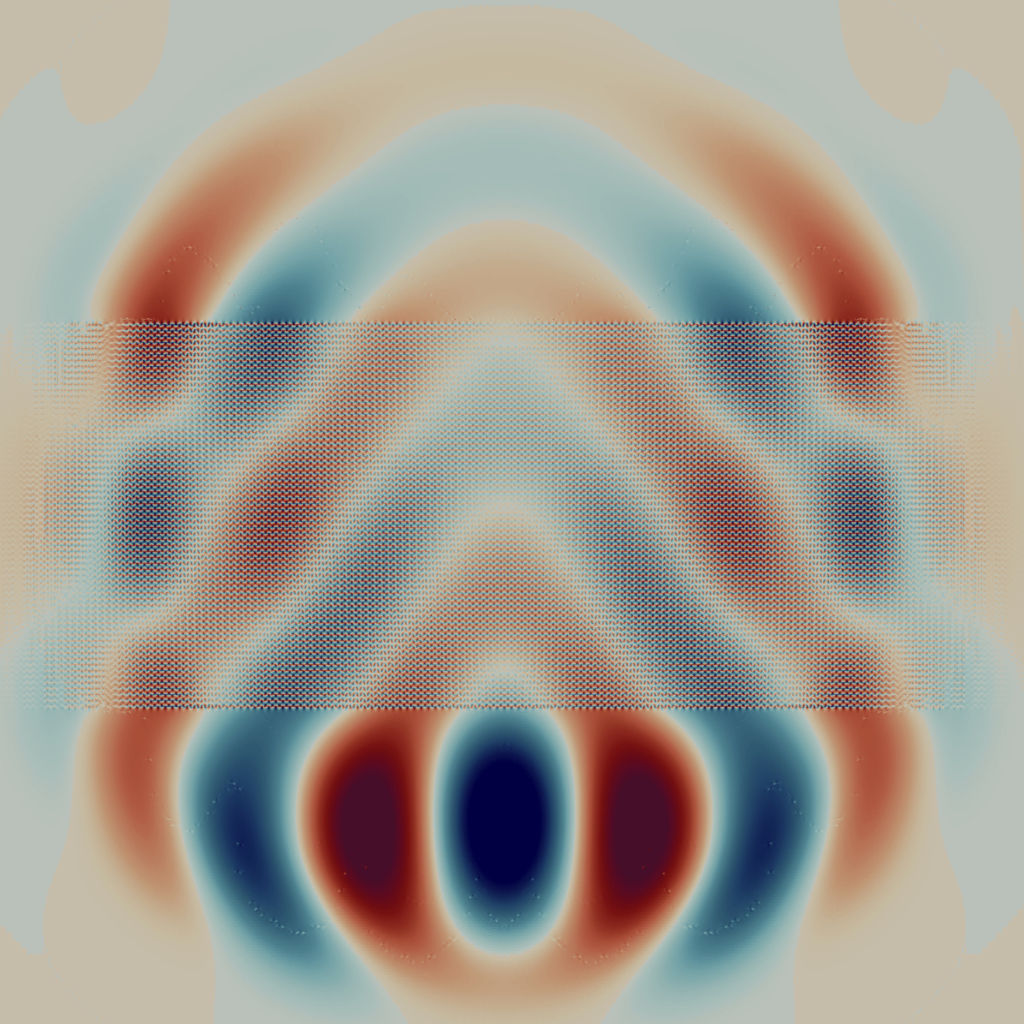}}
  \caption{
    \textcolor{black}{%
    (a) Geometry of a plasmonic-crystal slab of height $1$ with corrugated
    layers of 2D material, where the corrugation is sinusoidal and the
    period is equal to the spacing, $d$ (in the schematic $d\,=\,2^{-3}$).
    (b)} \textcolor{black}{Real part of electric field in the
    $y$-direction in the corresponding homogenization limit (as $d\to 0$),
    via solution of cell problem~\eqref{eq:cell_problem} and computation of
    the homogenized solution given by \eqref{eq:maxwellhomogenized}.
    \textcolor{black}{(c-e)} Real part of  electric field in the
    $y$-direction, based on direct numerical simulations
    of~\eqref{eq:maxwell}--\eqref{eq:compatibility_condition} (see
    \cite{maier2017a}), for decreasing spacing , $d$
    (\textcolor{black}{i.\,e.}, increasing degree of scale
    separation): \textcolor{black}{(c)} $d\,=\,2^{-4}$;
    \textcolor{black}{(d)} $d\,=\,2^{-5}$; \textcolor{black}{(e)}
    $d\,=\,2^{-6}$.}}
  \label{fig:homogenization_limit}
\end{figure}
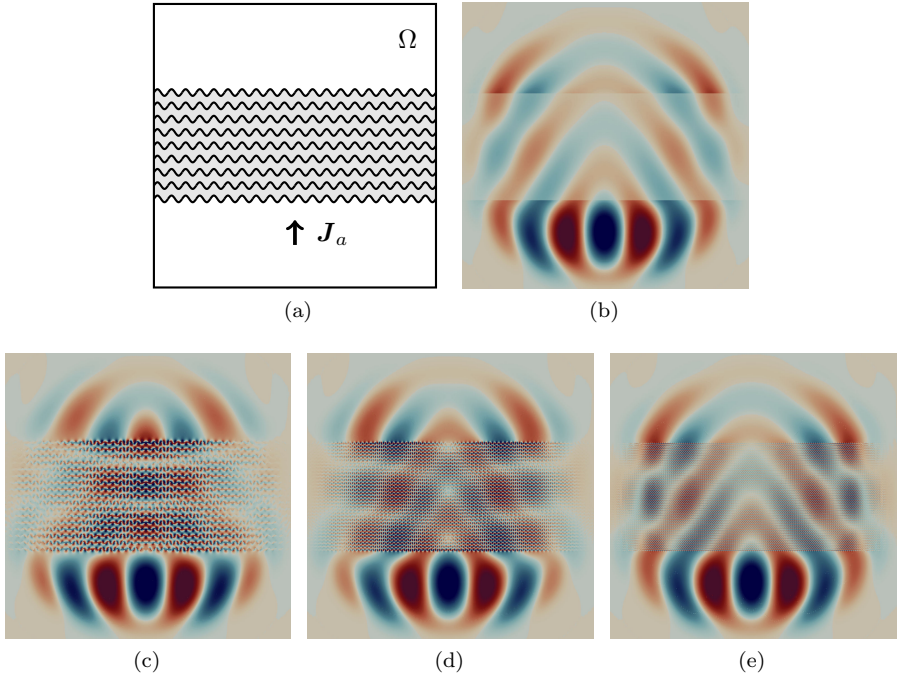
\textcolor{black}{%
We demonstrate the applicability of our computational platform by use of a
prototypical scattering configuration. \textcolor{black}{It} consists
of a dipole situated in close proximity to a  plasmonic crystal (slab) of
height $1$ consisting \textcolor{black}{of many layers of corrugated
2D sheets} at distance $d$; see Figure~\ref{fig:homogenization_limit}. The
sheet corrugation is described by a sine curve of amplitude $d/4$ and
period $d$ (Figure~\ref{fig:homogenization_limit}a). The computational
result of our \textcolor{black}{proposed scheme involving the
homogenized system} is shown in Figure \ref{fig:homogenization_limit}b.
\textcolor{black}{In addition, } Figures
\ref{fig:homogenization_limit}c-e show the results of a direct numerical
simulation \textcolor{black}{of the electric field $\vE$ described by
\eqref{eq:maxwell} and \eqref{eq:jump} for an} increasing level of scale
separation ($d\,=\,2^{-4}$, $2^{-5}$, $2^{-6}$, respectively).}
\textcolor{black}{The direct numerical computations shown in
Figure~\ref{fig:homogenization_limit}b-d require a very fine resolution
that results in linear systems with up to $1.6\,\times10^{7}$ unknowns.
In contrast, our proposed computational framework can be efficiently
implemented by use of a moderate resolution in the approximation scheme of
around $3.2\times10^{4}$ unknowns, a number smaller by a factor of about
$500$ than the one for the direct numerical simulation. This comparison
demonstrates that already for a moderate scale separation of $d=2^{-6}$,
the computational platform leads to a significant saving of computational
resources.
}

\textcolor{black}{%
Our example demonstrates in addition that the predictive quality of the
homogenized solution improves with increasing scale separation: Typically,
we have a good agreement between the electromagnetic field $\{\vE,\vH\}$
that is described system \eqref{eq:maxwell} and \eqref{eq:jump} and the
homogenized field $\{\vME,\vMH\}$ given by \eqref{eq:maxwellhomogenized}
for a scale separation of $d\approx1/64$ and less, even for scattering
configuration with dominant near field character as shown in
Figure~\ref{fig:homogenization_limit}. In general, the upper limit of $d$
for which the homogenized system \eqref{eq:maxwellhomogenized} has some
predictive quality heavily depends on the geometry and the location of
sources.
}

\textcolor{black}{We should add the remark that our homogenization
result remains valid also for \textcolor{black}{current-carrying}
sources, $\vJa$, situated inside the geometry provided the function
$\vJa(\vx)$ is square integrable. This claim would typically require that
dipole sources have to be regularized, as in our computational example. In
the case with (perfect, unregularized) point sources the emerging near
field might not be captured in its entirety by our computational
framework.}

\subsection{Effective permittivity tensors in prototypical geometries}
\label{subsec:effective-perm}

In order to relate numerical results obtained by solving cell problem
\eqref{eq:cell_problem} to the averaging in~\eqref{eq:effectiveperm} and
compute physical quantities as a function of frequency, $\omega$, it is
necessary to use a suitable model for the material parameters. If
$\lambda^d=0$, the only modeling parameters that enter
\eqref{eq:cell_problem} are  $\e(\vx,\vy)$ and
$\frac1{i\omega}\sigma(\vx,\vy)$. We assume that $\e(\vx,\vy)=\e$ is
constant and that the tangential parts of $\sigma(\vx,\vy)=\sigma(\vy)$ are
given by spatially constant values. Thus, the only microscale ($\vy$-)
parameter dependence to be accounted for is the one in
$\frac1{i\omega}\sigma(\vx,\vy)$; \color{black} cf.
Figure~\ref{fig:cell_problems} for geometries of interest. \color{black}

The surface conductivity, $\sigma^d$, of doped graphene can  plausibly be
described by the Kubo formula, which takes into account electronic
excitations and temperature effects~\cite{bludov13}. However,
\color{black} in a suitable parameter regime that includes terahertz
frequencies, \color{black} in which fine-scale SPPs on graphene can
typically be generated, \textcolor{black}{it has been shown that the
Kubo formula reduces to the (much simpler) Drude model~\cite{falkovsky07}.}
By this model, the tangential components of a spatially constant $\sigma^d$
are given by the formula
\begin{align*}
  \sigma^d=\frac{i\,e^2\,E_F}{\e_0\,\pi\hbar^2\,\big(\omega+i/\tau\big)}~.
\end{align*}
Here, $e$ is the electron charge, $\hbar$ denotes the (reduced) Planck
constant, $\e_0$ is the vacuum permittivity, $E_F$ denotes the Fermi
energy, and $\tau$ is the electronic relaxation time. In this context, an
$\vx$-dependence of $\sigma^d$ may arise from spatial variations of the
parameters $E_F$ and $\tau$.

\color{black} We proceed to carry out numerical computations for $\eff_{ij}$. \color{black} Using typical parameter values for graphene~\cite{bludov13}, we set $\tau\;=\;0.5\,10^{-12}\text{s}$ and apply the following rescalings:
$E_F\;=\;\tilde E_F\,10^{-19}\text{J}$, $\omega\;=\;\tilde
\omega\,10^{14}\text{Hz}$, $d\;=\;\tilde d\,10\,\text{nm}$,
with $\;0 \le\tilde E_F \le 1.6\;$, $\;0.5\le\tilde\omega\le4.0\;$, and
$\;0\le\tilde d\le20.0$. The surface average in
\eqref{eq:effectiveperm} of the effective permittivity tensor, $\eff$, has
the constant ($\vy$-independent) prefactor
\begin{align}
  \eta:=\frac{1}{i\omega}\,\sigma \;=\; \frac{1}{i\omega\,d}\,\sigma^d \;=\;
  82.9\, \frac{\tilde E_F}{\tilde
  d\,\tilde\omega\,\big(\tilde\omega+0.02\,i\big)}~.
  \label{eq:drude_rescaled}
\end{align}
Utilizing definition~\eqref{eq:drude_rescaled} we conveniently express the
matrix elements of $\eff$ given in \eqref{eq:effectiveperm} as
\begin{align*}
  \frac{\eff_{ij}}{\e}=
  1\;-\;\frac\eta\e\, \int_\Sigma
  P_t(\vy)\big(\vec e_j+\nabla_y\chi_j(\vx,\vy)\big) \cdot \vec
  e_i\doy~.
\end{align*}
Here, $P_t(\vy)$ denotes the projection onto the tangential space of
$\Sigma^d$ at point $\vy$. This formula for $\eff_{ij}$ uncovers an
important property: Up to a factor of $\e$, computational results for the
cell problem only depend on the ratio $\eta/\e$. Thus, it is sufficient to
compute values for $\eff_{ij}$ by setting $\e=1$.

Let us now recall the discussion about vanishing correctors in
Sections~\ref{sec:enz}\,\ref{subsec:vanishing},\ref{subsec:corrector}. For
the prototypical geometries of Figure~\ref{fig:cell_problems}, the
respective effective permittivities have the matrix forms:

\begin{align}
  \eff_{\text S}
  =\begin{pmatrix}
    \e_{11}^{\text S} & 0 & 0\\
    0 & \e & 0\\
    0 & 0 & \e_{33}^{\text S}
  \end{pmatrix},
  \qquad
  \eff_{\text R}
  =\begin{pmatrix}
    \e_{11}^{\text R} & 0& 0\\
    0 & \e & 0 \\
    0 & 0 & \e_{33}^{\text R}
  \end{pmatrix},
  \qquad
  \eff_{\text T}
  =\begin{pmatrix}
    \e_{11}^{\text T} & 0& 0\\
    0 & \e_{22}^{\text T} & 0\\
    0 & 0 & \e_{33}^{\text T}
  \end{pmatrix}~.
  \label{eq:eff_tensors}
\end{align}
Here, the subscript (S, R or T) for each matrix indicates the type of
geometry: S corresponds to the geometry with planar sheets and no edges
(Figure~\ref{fig:cell_problems}a); R stands for the nanoribbons geometry
(Figure~\ref{fig:cell_problems}b); and T corresponds to the nanotubes
geometry (Figure~\ref{fig:cell_problems}c). Due to vanishing corrector
components, the matrix element $\e_{11}^{\text{S}}$  is given in closed
form by~\eqref{eq:simplified_eff}, viz.,
\begin{align}
  \e_{11}^\text{S} = 1 - \eta~;\quad \e_{11}^{\text{S}}=\e_{33}^{\text{S}}=
  \e_{33}^{\text{R}}=\e_{33}^{\text{T}}~.
  \label{eq:eff_trivial}
\end{align}

In the remaining numerical computations, we focus on the \color{black}
more complicated \color{black} geometries with nanoribbons and nanotubes
(Figure~\ref{fig:cell_problems}b,c). To determine $\e_{11}^{\text{R}}$ and
$\e_{11}^{\text{T}}=\e_{22}^{\text{T}},$ we solve cell problem
\eqref{eq:cell_problem} directly and compute the average
by~\eqref{eq:effectiveperm}. The computations are carried out by the finite
element toolkit deal.II~\cite{dealii85}.
\begin{figure}[tbp]
  \begin{center}
    \subfloat[]{\input{nano-ribbons.tex}}
    \subfloat[]{\input{nano-tubes.tex}}
  \end{center}
  \caption{
    Plots of real and imaginary parts of matrix elements of $\eff$ as a
    function of frequency by~\eqref{eq:eff_tensors} for the geometries of
    Figures~\ref{fig:cell_problems}b, c. (a) $\e^{\text{R}}_{11}$ for
    nanoribbons; and (b) $\e^{\text{T}}_{11}$ or $\e^{\text{T}}_{22}$
    ($\e^{\text{T}}_{11}=\e^{\text{T}}_{22}$) for nanotubes. The shaded
    area indicates the frequency regime for negative real part in each
    case.}
  \label{fig:effective_epsilon}
\end{figure}
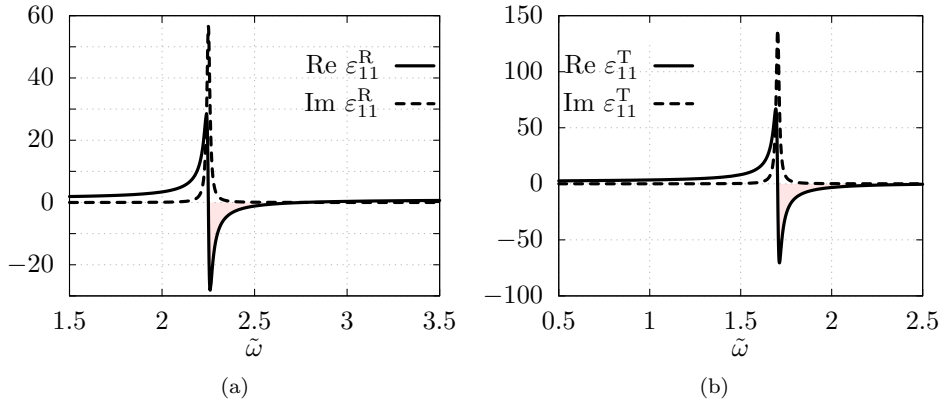
To this end, we sample over the frequency range $0.5\le\tilde\omega\le4.0$
with the choice $\tilde E_F=1$ and $\tilde d = 20.72$. The real and
imaginary parts of $\e^{\text{R}}_{11}$ and $\e^{\text{T}}_{11}$ are
plotted as a function of frequency in
Figures~\ref{fig:effective_epsilon}a, b.
\textcolor{black}{The real and imaginary parts of this $\e$ satisfy the
Kramers-Kronig relations. We observe that each of these functions exhibits
a strong Lorentzian resonance at around
($\e_{11}^{\text{R}}$:) $\tilde\omega_R=2.25002(2)$, and
($\e_{11}^{\text{T}}$:) $\tilde\omega_R=1.70285(3)$.
It is of interest to note that the ENZ effect occurs at the following
(rescaled) frequencies:
($\e_{11}^{\text{R}}$:) $\tilde\omega\approx2.752$,
($\e_{11}^{\text{T}}$:) $\tilde\omega\approx2.755$.
Here, both frequencies are determined by interpolating the computational
results displayed in Figure~\ref{fig:effective_epsilon} (obtained
via the solution of cell problem \eqref{eq:cell_problem} with a finite element
method) with a Lorentzian function. Note that \eqref{eq:eff_trivial}
implies the critical frequency $\tilde\omega\approx2$ for
$\e_{11}^{\text{S}}$ where
$\e_{11}^{\text{S}}=\e_{33}^{\text{S}}=\e_{33}^{\text{R}}=\e_{33}^{\text{T}}$.
}
%

\section{Conclusion}
\label{sec:conclusion}

In this paper, we carried out a formal asymptotic analysis to homogenize
time-harmonic Maxwell's equations in plasmonic crystals made of conducting
sheets with microscale spacing $d$, in the limit as $d\to 0$ via suitable
scalings of material parameters. The homogenized system features an
effective permittivity tensor, $\eff$, given by an averaging procedure that
involves a weighted volume average of the bulk permittivity, as well as a
weighted surface average of the surface conductivity, $\sigma^d$, and a
weighted line average of the line charge density, $\lambda^d$, of the
sheets. The vector-valued corrector field of this procedure solves a closed
cell problem. We showed analytically how the combination of the
complex-valued $\sigma^d$ and $\lambda^d$ yields an ENZ effect; and defined
the related critical spacing which depends on $\sigma^d$ and $\lambda^d$.

\textcolor{black}{%
The introduction of the line charge density, $\lambda^d$, in our microscale
model and homogenization procedure is an aspect that deserves particular
emphasis. We believe that this feature has eluded previous works in
plasmonics. In fact, a general and precise quantitative description of the
influence of the physics at material edges on effective optical parameters
of plasmonic crystals made of 2D materials is a largely open problem. Edge
effects have received significant attention recently because of
observations of EPPs~\cite{krenn2004,fei2015}. Our discussion of a
homogenization procedure involving a line charge density on material edges
contributes to this effort. For example, we demonstrated that the
(generalized) plasmonic thickness has an algebraic dependence on the line
charge density different from that on the surface conductivity of the 2D
material. In addition, our formalism suggests a few  mathematical problems
(regarding the well-posedness and formal proof of homogenization results)
which were not addressed here.
}

\textcolor{black}{%
We also discussed how our homogenization result can be
incorporated into well established computational approaches for time
harmonic Maxwell equations. This procedure involves the computation of
effective material parameters by approximation of the solution of the
corresponding cell problems and averaging. We demonstrated the feasibility
of this approach in a geometry with corrugated sheets; and computed the
Lorentz-type resonance of two prototypical microscopic geometries.
The compuational framework that we introduced paves the road for future,
systematic computational investigations of complicated design problems in
the plasmonics of 2D materials \cite{molesky2018,liberal2017,miller2012}.
}

Our results have a few limitations and point to open problems in asymptotics. For
instance, the asymptotic analysis is based on a strong periodicity
assumption. Further, we do not discuss boundary
layers in the homogenization procedure due to the interaction of the
microstructure with boundaries of the (scattering) domain. It is also worth
mentioning that scaling assumptions different from the ones chosen here may
lead to different homogenization results.

\appendix

\section{Two-scale expansion and asymptotics}
\label{sec:detailed_asymptotics}

\textcolor{black}{In this appendix, we carry out in detail the asymptotic analysis
that was outlined in Section~\ref{sec:asymptotics}. As a
first step, we apply the two-scale asymptotic expansion
to \eqref{eq:maxwell} and \eqref{eq:jump}. By collecting all terms of the
order of $d^{-1}$ in region $\Omega\times Y$ and of the order of $d^0$ on
$\Omega\times\Sigma$, we obtain the equations}
\begin{align}
  \hspace{-1em}
  \begin{cases}
    \begin{aligned}
      \nabla_y\times\vE^{(0)} &= 0~,
      &
      \nabla_y\times\vH^{(0)} &= 0~,
      \\[0.3em]
      \nabla_y\cdot(\e\vE^{(0)}) &= 0~,
      &
      \nabla_y\cdot\vH^{(0)} &= 0~;
      \\[0.3em]
      \js{\vn\times\vE^{(0)}} &= 0~,
      &
      \js{\vn\times\vH^{(0)}} &= 0~,
      \\[0.3em]
      \js{\vn\cdot(\e\vE^{(0)})} &=
      \frac{1}{i\omega}\nabla_y\cdot(\sigma\vE^{(0)})~,
      &
      \js{\vn\cdot\vH^{(0)}} &= 0~.
    \end{aligned}
  \end{cases}
  \hspace{-1em}
  \label{eq:orderdm1}
\end{align}
In a similar vein, a second set of equations is obtained by collecting
all terms of the order of $d^0$ in $\Omega\times Y$ and of the order of
$d^1$ on $\Omega\times\Sigma$, viz.,
\begin{align}
  \begin{cases}
    \begin{aligned}
      \nabla_x\times\vH^{(0)} +
      \nabla_y\times\vH^{(1)} &= -i\omega\e\vE^{(0)}+\vJa~,
      \\[0.3em]
      \nabla_x\times\vE^{(0)} +
      \nabla_y\times\vE^{(1)} &= i\omega\mu_0\vH^{(0)}~,
      \\[0.3em]
      \nabla_x\cdot(\e\vE^{(0)}) +
      \nabla_y\cdot(\e\vE^{(1)}) &= \frac{1}{i\omega}{\nabla_x\cdot\vJa}~,
      \\[0.3em]
      \nabla_x\cdot\vH^{(0)} +
      \nabla_y\cdot\vH^{(1)} &= 0~,
    \end{aligned}
  \end{cases}
  \label{eq:order0a}
\end{align}
and
\begin{align}
  \begin{cases}
    \begin{aligned}
      \js{\vn\times\vE^{(1)}} &= 0~,
      &
      \js{\vn\times\vH^{(1)}} &= \sigma\vE^{(0)}~,
      \\[0.3em]
      \js{\vn\cdot(\e\vE^{(1)})}  &=
      \frac{1}{i\omega}\big(\nabla_x\cdot(\sigma\vE^{(0)})+
      \nabla_y\cdot(\sigma\vE^{(1)})\big)~,
      &
      \js{\vn\cdot\vH^{(1)}} &= 0~.
    \end{aligned}
  \end{cases}
  \label{eq:order0b}
\end{align}
In \eqref{eq:orderdm1}--\eqref{eq:order0b}, the differential equations are
valid for $(\vx,\vy)\in\Omega\times Y$; while the jump conditions hold for
$(\vx,\vy)\in\Omega\times\Sigma$.

So far, we used boundary conditions on the interior of the sheets,
$\Sigma^d\cap \Omega\setminus \partial\Sigma^d$. We now focus
on~\eqref{eq:compatibility_condition}, imposed on the edges,
$\partial\Sigma^d\cap\Omega$.
First, we introduce the rescaled jump $\je{\,.\,}$ in a way analogous to
definition~\eqref{eq:singular_jump} for the jump $\jse{\,.\,}$, viz.,
\begin{align*}
  \je{\vec F}(\vx,\vy)
  \,:=\,
  \lim_{\alpha\searrow0}\,
  \int_{-\alpha}^{\;\alpha}
  \big(\vec F(\vx,\vy+\alpha^2\vec n+\zeta\vn)-\vec F(\vx,\vy-\alpha^2\vec
  n+\zeta\vn)\big)\,
  \text{d}\zeta\qquad
  \vy\in\partial\Sigma~.
\end{align*}
By carrying out the leading-order asymptotic expansion for the singular
jump of~\eqref{eq:singular_jump}, we see that
$\jse{\,.\,}\;\to\;d\,\je{\,.\,}$. Consequently, the first condition of
\eqref{eq:compatibility_condition} is expanded to
\begin{align}
  \je{\vec n\times\vH^{(1)}}\;&=\;\lambda\vE^{(0)}
  \quad
  \text{on }\Omega\times\big(\partial\Sigma\setminus \partial Y\big)~.
  \label{eq:compatibility_rescaled1}
\end{align}
Furthermore, the expansion of the second one of
conditions~\eqref{eq:compatibility_condition} results in the following
boundary conditions to the two lowest orders in $d$:
\begin{align}
  \vec n \cdot\big(\sigma\vE^{(0)}\big)\;=\;
  \nabla_y\cdot\big(\lambda\vE^{(0)}\big)~,
  \qquad
  \vec n \cdot\big(\sigma\vE^{(1)}\big)\;=\;
  \nabla_x\cdot\big(\lambda\vE^{(0)}\big) +
  \nabla_y\cdot\big(\lambda\vE^{(1)}\big)~,
  \label{eq:compatibility_rescaled}
\end{align}
which hold on $\Omega\times\big(\partial\Sigma\setminus \partial Y\big)$.

\subsection{Characterization of $\vE^{(0)}$ and $\vH^{(0)}$}
\label{subsec:characterization}

We now use~\eqref{eq:orderdm1} to characterize $\vE^{(0)}(\vx,\vy)$ and
$\vH^{(0)}(\vx,\vy)$ in more detail. Since a conservative periodic vector
field is the sum of a constant vector and the gradient of a periodic
function (potential), we can write the general solution to the first
equation of \eqref{eq:orderdm1} as
\begin{equation}
  \label{eq:ansatz1}
  \vE^{(0)}(\vx,\vy) = \vME(\vx) + \nabla_y\varphi(\vx,\vy)~,
  \quad \varphi(\vx,\vy) = \sum_j\chi_j(\vx,\vy)\mathcal{E}_j(\vx)~.
\end{equation}
In this vein, a conservative and divergence-free periodic vector field must
be constant. Hence, the general solution to the second and fourth laws of
the first group of equations in~\eqref{eq:orderdm1} is given by
 \begin{equation}
  \label{eq:ansatz2}
  \vH^{(0)}(\vx,\vy) = \vMH(\vx)~.
\end{equation}
The functions $\vME(\vx)$, $\vMH(\vx)$ and $\chi_j(\vx,\vy)$ ($j=1,2,3$)
are further characterized below.

\subsection{Derivation of cell problem}
\label{subsec:cell_problem}

Next, we derive a closed set of equations that fully describe the functions
$\chi_j(\vx, \vy)$ introduced in~\eqref{eq:ansatz1}. These equations
comprise the cell problem accounting for the microstructure details.

First, we substitute \eqref{eq:ansatz1} into the respective, zeroth-order
expressions in~\eqref{eq:orderdm1} and~\eqref{eq:compatibility_rescaled}.
Specifically, we use the third law of the first group of equations
in~\eqref{eq:orderdm1}; the first and third jump conditions in
\eqref{eq:orderdm1}; and the first condition in
\eqref{eq:compatibility_rescaled}. Thus, we obtain the following equations:
\begin{align*}
  \begin{cases}
    \begin{aligned}
      \sum_{j}\nabla_y\cdot\Big(\e(\vx,\vy)\big(\vec
      e_j+\nabla_y\chi_j(\vx,\vy)\big)\Big)\mathcal{E}_j(\vx)
      \;&=\;0
      \quad\text{in }\Omega\times Y~,
      \\[0.3em]
      \sum_j\js{\vn\times\Big(\e(\vx,\vy)\big(\vec
      e_j+\nabla_y\chi_j(\vx,\vy)\big)\Big)}\:\mathcal{E}_j(\vx)
      \;&= \;0
      \quad\text{on }\Omega\times\Sigma~,
      \\[0.3em]
      \sum_j\js{\vn\cdot\Big(\e(\vx,\vy)\big(\vec
      e_j+\nabla_y\chi_j(\vx,\vy)\big)\Big)}\:\mathcal{E}_j(\vx)
      \;&=\;
      \\[0.3em]
      \frac{1}{i\omega}\sum_j\nabla_y\cdot\Big(\sigma(\vx,\vy)
      &\big(\vec e_j+\nabla_y\chi_j(\vx,\vy)\big)\Big)\mathcal{E}_j(\vx)
      \quad\text{on }\Omega\times\Sigma~,
      \\[0.3em]
      \sum_j\vec n\cdot\Big(\sigma(\vx,\vy)\big(\vec
      e_j+\nabla_y\chi_j(\vx,\vy)\big)\Big)\:\mathcal{E}_j(\vx)
      \;&= \;
      \\[0.3em]
      \sum_j\nabla_y\cdot\Big(\lambda(\vx,\vy) &\big(\vec
      e_j+\nabla_y\chi_j(\vx,\vy)\big)\Big)\mathcal{E}_j(\vx) \quad\text{on
      }\Omega\times\big(\partial\Sigma\setminus\partial Y\big)~.
    \end{aligned}
  \end{cases}
\end{align*}
Here, $\mathcal{E}_j$ and $\chi_j$ are coupled. To simplify this
description, we treat each term containing $\chi_j(\vx,\vy)$, which
accounts for the microscale behavior of $\vE^{(0)}(\vx,\vy)$, as
independent from $\mathcal{E}_j(\vx)$ ($j=1,2,3$). Thus, the above
equations decouple into three distinct problems, one for each $\chi_j$,
as displayed in~\eqref{eq:cell_problem}. Equations~\eqref{eq:cell_problem} along with
the condition that $\chi_j(\vx,\,.\,)$ be $Y$-periodic form the desired,
closed cell problem. Notice that $\vx$ plays the role of a parameter.
Hence, the above cell problem uniquely describes $\chi_j(\vx,\,.\,)$ for
any given (macroscopic) point $\vx$.

\subsection{Homogenized macroscale problem}
\label{subsec:homogenized_problem}

Our remaining task is to derive corresponding macroscale equations for the
functions $\vME(\vx)$ and $\vMH(\vx)$. We start by substituting
\eqref{eq:ansatz1} and \eqref{eq:ansatz2} into the first equation
of~\eqref{eq:order0a} and averaging (in cell $Y$) over the fast variable,
$\vy$. Hence, we obtain the following expression:
\begin{multline}
  \nabla_x\times\vMH(\vx) +
  \int_Y\nabla_y\times\vH^{(1)}(\vx,\vy)\dy =
  \\
  -i\omega\sum_j \int_Y\e(\vx,\vy)\big(\vec
  e_j+\nabla_y\chi_j(\vx,\vy)\big)\dy\:\mathcal{E}_j(\vx) + \vJa(\vx)~.
  \label{eq:macro_1}
\end{multline}
By use of the Gauss theorem and the $Y$-periodicity of $\vH^{(1)}$, the
second term on the left-hand side of the above equation is written as
\begin{align*}
  \int_Y\nabla_y\times\vH^{(1)}\dy
  &= -\int_\Sigma\js{\vn\times\vH^{(1)}}\doy
     -\int_{\partial\Sigma\setminus\partial Y} \jse{\vec
     n\times\vH^{(1)}}\text{d} s
  \\
  &= -\int_\Sigma\sigma(\vx,\vy) \vE^{(0)}\doy -\int_{\partial\Sigma\setminus\partial
     Y} \lambda(\vx,\vy) \vE^{(0)}\text{d} s
  \\
  &=
  -\sum_j
  \left\{
  \int_\Sigma \sigma(\vx,\vy) \big(\vec e_j+\nabla_y\chi_j\big)\doy
  + \int_{\partial\Sigma\setminus Y} \lambda(\vx,\vy) \big(\vec
    e_j+\nabla_y\chi_j\big)\text{d} s
  \right\}
  \:\mathcal{E}_j(\vx)~.
\end{align*}
In the above, the second equality comes from using the second jump
condition of \eqref{eq:order0b}, and \eqref{eq:compatibility_rescaled1}.
The third equality follows from~\eqref{eq:ansatz1}. Let us now define the
\emph{effective permittivity tensor} $\eff$ by
\begin{multline}
  \eff_{ij}(\vx):=
  \int_Y\e(\vx,\vy)\big(\vec e_j+\nabla_y\chi_j(\vx,\vy)\big)
  \cdot\vec e_i \dy
  -\frac1{i\omega}
  \int_\Sigma
  \sigma(\vx, \vy)\big(\vec e_j+\nabla_y\chi_j(\vx,\vy)\big) \cdot \vec
  e_i\doy
  \\
  -\frac1{i\omega}
  \int_{\partial\Sigma\setminus\partial Y}
  \lambda(\vx, \vy)\big(\vec e_j+\nabla_y\chi_j(\vx,\vy)\big) \cdot \vec
  e_i\,\text{d} s~.
  \label{eq:sm_effectiveperm}
\end{multline}
In view of this definition of $\eff$, \eqref{eq:macro_1} takes the form
\begin{align*}
  \nabla\times\vMH &= -i\omega\eff\vME+\vJa~,
\end{align*}
which describes the effective Amp\'ere-Maxwell law.

The last three equations of~\eqref{eq:order0a} can be manipulated in a
similar fashion. For example, consider the third equation. By
using~\eqref{eq:ansatz1} and~\eqref{eq:ansatz2} and averaging over the fast
variable, $\vy$, we obtain
\begin{multline}
  \nabla_x\cdot\Big(
  \sum_j \int_Y\e(\vx,\vy)\big(\vec
  e_j+\nabla_y\chi_j(\vx,\vy)\big)\dy\:\mathcal{E}_j(\vx)
  \Big)
  \\
  + \int_Y\nabla_y\cdot(\e(\vx,\vy)\vE^{(1)}(\vx,\vy))\dy =
  \frac{1}{i\omega}{\nabla_x\cdot\vJa(\vx)}~.
  \label{eq:temp2}
\end{multline}
Next, we manipulate the second term of the left-hand side by applying the
Gauss theorem and utilizing the second jump condition of
\eqref{eq:order0b}, as follows:
\begin{align*}
  \int_Y\nabla_y\cdot(\e\vE^{(1)})\dy
  &= -\int_\Sigma\big[\nu\cdot(\e\vE^{(1)})\big]_\Sigma\doy
  \\
  &= -\int_\Sigma\frac{1}{i\omega} \big(\nabla_x\cdot(\sigma\vE^{(0)})+
  \nabla_y\cdot(\sigma\vE^{(1)})\big)\doy
  \\
  &= -\frac{1}{i\omega} \nabla_x\cdot\Big(\sum_j \int_\Sigma \sigma(\vx,
  \vy)\big(\vec e_j+\nabla_y\chi_j(\vx,\vy)\big)
  \doy\,\mathcal{E}_j(\vx)\doy\Big)
  \\
  &\qquad\qquad \qquad\qquad \qquad\qquad \qquad\qquad \qquad\qquad
  -\frac{1}{i\omega} \int_\Sigma\nabla_y\cdot(\sigma\vE^{(1)})\doy~.
\end{align*}
By applying the Gauss theorem and utilizing the second boundary condition
in \eqref{eq:compatibility_rescaled}, we find
\begin{align*}
  \int_\Sigma\nabla_y\cdot(\sigma\vE^{(1)})\doy
  \;&=\;
  \int_{\partial \Sigma}\vec n\cdot\big(\sigma\vE^{(1)}\big)\,\text{d}s
  \;=\;
  \int_{\partial \Sigma\setminus\partial Y}\vec
  n\cdot\big(\sigma\vE^{(1)}\big)\,\text{d}s
  \\
  &=\;
  \int_{\partial \Sigma\setminus\partial Y}
  \nabla_x\cdot(\lambda\vE^{0})+
  \nabla_y\cdot(\lambda\vE^{1})
  \,\text{d}s
  \\
  &=\;
  \nabla_x\cdot\left(
  \sum_j
  \int_{\partial\Sigma\setminus\partial Y}
  \lambda(\vx, \vy)\big(\vec e_j+\nabla_y\chi_j(\vx,\vy)\big)
  \,\text{d} s
  \,
  \vME_j(\vx)
  \right)\;.
\end{align*}
The second equality (exclusion of $\partial Y$) exploits the fact that
$\sigma\vE^{(1)}(\vx,.)$ is $Y$-periodic and, thus, single valued on
$\partial\Sigma\cap\partial Y$. Notice, however, that the normal
vector $\vec n$ changes sign. The last integral in the second line vanishes
because $\partial\Sigma\setminus Y$ has no boundary. Substituting the
result of these manipulations into~\eqref{eq:temp2} and
utilizing~\eqref{eq:sm_effectiveperm}, we obtain
\begin{align*}
  \nabla\cdot(\eff\vME) = \frac{1}{i\omega}{\nabla\cdot\vJa}.
\end{align*}

Similar steps can be applied to the remaining equations
of~\eqref{eq:order0a}. The homogenized system \textcolor{black}{finally
reads}
\begin{align}
  \begin{cases}
    \begin{aligned}
      & \nabla\times\vME = i\omega\mu_0\vMH~,
      \quad
      \nabla\times\vMH = -i\omega\eff\vME+\vJa~;
      \\[0.3em]
      &\nabla\cdot(\eff\vME) = \frac{1}{i\omega}{\nabla\cdot\vJa}~,
      \quad
      \nabla\cdot\vMH = 0~.
    \end{aligned}
  \end{cases}
  \label{eq:sm_maxwellhomogenized}
\end{align}

\section{Derivation of compatibility conditions from Gauss' and Amp\`ere's law}
\label{sec:gauss_law}

In this section, we give a short derivation of the internal compatibility
condition \eqref{eq:compatibility_condition}.

For an arbitrary volume $V$ in $\Omega,$ the integral over the
\emph{electric displacement in normal direction} has to be equal to the
integral over the total charge density contained in $V,$
\begin{align*}
  i\omega\int_{V}\rho\dx
  \;=\;
  \int_{\partial V}\vec n_{\partial V}\cdot\big(\e^d\vE^d\big)\dox.
\end{align*}
We now choose $V$ to be an arbitrary rectangular box containing a part of
the edge $\partial\Sigma^d$. We extend the sheet over the edge parallel in
$\vec n$-direction, and assume $\sigma^d=0$ in the extension. The box shall
be of vanishing length and height, and with top and bottom faces parallel
to the extended sheet $\Sigma^d_\ast$; see Figure~\ref{fig:box}. Then,
\begin{align*}
  \lim_{\text{height}\to0}
  i\omega
  \int_{V}\rho\dx \;&=\;
  \lim_{\text{height}\to0}\;
  \int_{\partial V}\vec n_{\partial V}\cdot\big(\e^d\vE^d\big)\dox
  \\
  &=\;
  \lim_{\text{height}\to0}\;
  \int_{\text{top}}\vn\cdot\big(\e^d\vE^d\big)^{\text{above}}\dox
  -
  \int_{\text{bottom}}\vn\cdot\big(\e^d\vE^d\big)^{\text{below}}\dox
  \\
  &=\;
  +\int_{V\cap\Sigma^d_\ast}\vn\cdot\big(\e^d\vE^d\big)^{\text{above}}\dox
  -\int_{V\cap\Sigma^d_\ast}\vn\cdot\big(\e^d\vE^d\big)^{\text{below}}\dox.
\end{align*}
Here, $\vec n_{\partial V}$ denotes the outward pointing unit normal on
faces of the volume $V$ and $\vn$ is the normal field on $\Sigma^d_\ast$.
Now, utilizing the third jump condition in \eqref{eq:jump} we conclude that
\begin{align*}
  \lim_{\text{height}\to0}
  i\omega
  \int_{V}\rho\dx \;&=\;
  \int_{V\cap\Sigma^d_\ast}\jsd{\vn\cdot\big(\e^d\vE^d\big)}\dox
  \\
  &=\;\int_{V\cap\Sigma^d_\ast}\nabla\cdot\big(\sigma^d\vE^d\big)\dox
  \\
  &=\;\int_{\partial V\cap\Sigma^d_\ast}\vec
  n\cdot\big(\sigma^d\vE^d\big)\dox.
\end{align*}
Here, $\vec n$ is the outward pointing normal on the edge, see
Figure~\ref{fig:box}. By keeping the width (dimension parallel to the edge
$\partial\Sigma^d$) fixed and in the limit of vanishing height and length,
we conclude that the volume integral over the charge density $\rho$ reduces
to
\begin{align*}
  \lim_{\text{height}\to0}\;
  \lim_{\text{length}\to0}\;
  i\omega
  \int_{V}\rho\dx \;=\;
  \int_{V\cap\partial\Sigma^d}\nabla\cdot\big(\lambda^d(\vx)\vE^d(\vx)\big)
  \,\text{d} s.
\end{align*}
Consequently,
\begin{align*}
  \int_{V\cap\partial\Sigma^d}\nabla\cdot\big(\lambda^d(\vx)
  \vE^d(\vx)\big) \,\text{d} s
  =\int_{V\cap\partial\Sigma^d}
  \jsd{\vec n\cdot\big(\sigma^d\vE^d\big)}\dox.
\end{align*}
Due to the fact that $V$ was chosen arbitrarily, we conclude that
\begin{align*}
  \jsd{\vec n\cdot\big(\sigma^d\vE^d\big)} \;=\;
  \nabla\cdot\big(\lambda^d(\vx)\vE^d(\vx)\big)
\end{align*}
has to hold true pointwise on $\partial\Sigma^d$. But $\sigma^d$ vanishes
outside of $\Sigma^d$, thus
\begin{align*}
  \vec n\cdot\big(\sigma^d\vE^d\big) \;=\;
  \nabla\cdot\big(\lambda^d(\vx)\vE^d(\vx)\big)
  \quad\text{on }\partial\Sigma^d.
\end{align*}
\begin{figure}
  \begin{center}
    \begin{tikzpicture}[scale=0.75]
        \path[thick,draw,->] (-2.5,-0.6) -- (-2.0,-0.6);
        \path[thick,draw,->] (-2.5,-0.6) -- (-2.5,-0.1);
        \path[thick,draw,->] (-2.5,-0.6) -- (-2.8,-0.9);
        \node at (-2.5, +0.1) {\tiny height};
        \node at (-1.8, -0.9) {\tiny length};
        \node at (-3.0, -1.1) {\tiny width};
    \end{tikzpicture}
    \hspace{2em}
    \begin{tikzpicture}[scale=2.00]
        \path [color=black,thick, draw, dashed] (0.60, 0.15) -- (0.60,-0.10);
        \path [color=black,thick, draw, dashed] (1.65, 0.15) -- (1.65,-0.10);
        \path [color=black,thick, draw, dashed] (0.80, 0.35) -- (0.80, 0.10);
        \path [color=black,thick, draw, dashed] (1.85, 0.35) -- (1.85, 0.10);
        \path [color=black,thick, draw] (0.60, -0.10) -- (1.15, -0.10) -- (1.65, -0.10);
        \path [color=black,thick, draw, dashed] (0.80, 0.10) -- (1.35, 0.10) -- (1.85, 0.10);
        \path [color=black,thick, draw, dashed] (0.60, -0.10) to[out=35,in=210] (0.80, 0.10);
        \path [color=black,thick, draw] (1.65, -0.10) to[out=35,in=210] (1.80, 0.10);

        \path [fill=white]
          (1, 0) -- (3, 0) to[out=-30,in=170] (2.5, 0.5)
          -- (1.5, 0.5) to[out=170,in=-30] (1.0, 0.0);

        \path [fill=white]
          (-1, 0) to[out=-30,in=170] (1, 0) to[out=-30,in=170] (1.5, 0.5)
          to[out=170,in=-30] (-0.5, 0.5) to[out=170,in=-30] (-1.0, 0.0);
        \path [fill=black, fill opacity=0.2]
          (-1, 0) to[out=-30,in=170] (1, 0) to[out=-30,in=170] (1.5, 0.5)
          to[out=170,in=-30] (-0.5, 0.5) to[out=170,in=-30] (-1.0, 0.0);

        \path [thick, opacity=0.5, draw, dashed] (-1, 0) to[out=-30,in=170] (1, 0);
        \path [thick, draw] (1, 0) to[out=-30,in=170] (1.5, 0.5);
        \path [thick, opacity=0.5, draw, dashed] (-0.5, 0.5) to[out=-30,in=170] (1.5, 0.5);
        \path [thick, opacity=0.5, draw, dashed] (-1.0, 0.0) to[out=-30,in=170] (-0.5, 0.5);

        \path [thick, opacity=0.5, draw, dashed] (1, 0) -- (2, 0);
        \path [thick, opacity=0.5, draw, dashed] (2, 0) to[out=-30,in=170] (2.5, 0.5);
        \path [thick, opacity=0.5, draw, dashed] (1.5, 0.5) -- (2.5, 0.5);

        \path [color=black,thick, draw] (0.60, 0.40) -- (1.15, 0.40) -- (1.65, 0.40);
        \path [color=black,thick, draw] (0.80, 0.60) -- (1.35, 0.60) -- (1.85, 0.60);
        \path [color=black,thick, draw] (0.60, 0.40) to[out=35,in=210] (0.80, 0.60);
        \path [color=black,thick, draw] (1.65, 0.40) to[out=35,in=210] (1.80, 0.60);
        \path [color=black,thick, draw] (0.60, 0.40) -- (0.60, 0.15);
        \path [color=black,thick, draw] (1.65, 0.40) -- (1.65, 0.15);
        \path [color=black,thick, draw, dashed] (0.80, 0.60) -- (0.80, 0.35);
        \path [color=black,thick, draw] (1.85, 0.60) -- (1.85, 0.35);

        \path [color=black,thick, draw, dashed] (0.60, 0.15) -- (1.15, 0.15) -- (1.65, 0.15);
        \path [color=black,thick, draw, dashed] (0.80, 0.35) -- (1.35, 0.35) -- (1.85, 0.35);
        \path [color=black,thick, draw, dashed] (0.60, 0.15) to[out=35,in=210] (0.80, 0.35);
        \path [color=black,thick, draw, dashed] (1.65, 0.15) to[out=35,in=210] (1.80, 0.35);

        \node at (0.0, 0.2) {\small $\Sigma^d$};
        \path[thick,draw,->] (0.3,0.12) -- (0.3,0.24);
        \node at (0.3, 0.30) {\small $\vn$};

        \path [color=black,thick, draw] (0.55, 0.70) -- (0.75, 0.38);
        \node at (0.5, 0.75) {\small $V$};
        \path[thick,draw,->] (1.28,0.28) -- (1.4,0.28);
        \node at (1.5, 0.25) {\small $\vec n$};
    \end{tikzpicture}
  \end{center}
  \caption{
    Choice of rectangular box $V$; a curved, rectangular box containing a part of
    the edge $\partial\Sigma^d$. We extend the sheet over the edge parallel in
    $\vec n$-direction, and assume $\sigma^d=0$ in the extension. The box shall
    be of vanishing length and height, and with top and bottom faces parallel
    to the extended sheet $\Sigma^d$.
  }
  \label{fig:box}
\end{figure}
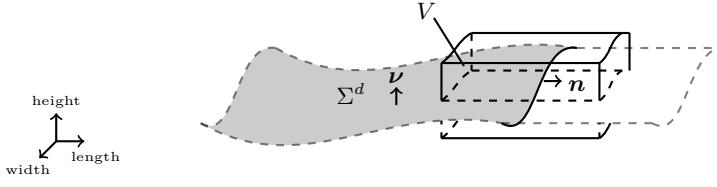
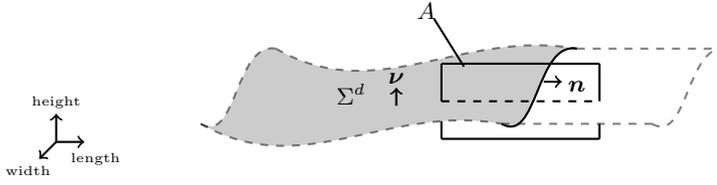
\begin{figure}
  \begin{center}
    \begin{tikzpicture}[scale=0.75]
        \path[thick,draw,->] (-2.5,-0.6) -- (-2.0,-0.6);
        \path[thick,draw,->] (-2.5,-0.6) -- (-2.5,-0.1);
        \path[thick,draw,->] (-2.5,-0.6) -- (-2.8,-0.9);
        \node at (-2.5, +0.1) {\tiny height};
        \node at (-1.8, -0.9) {\tiny length};
        \node at (-3.0, -1.1) {\tiny width};
    \end{tikzpicture}
    \hspace{2em}
    \begin{tikzpicture}[scale=2.00]
        \path [color=black,thick, draw] (0.60, 0.15) -- (0.60,-0.10);
        \path [color=black,thick, draw] (1.65, 0.15) -- (1.65,-0.10);
        \path [color=black,thick, draw] (0.60, -0.10) -- (1.15, -0.10) -- (1.65, -0.10);

        \path [fill=white]
          (1, 0) -- (3, 0) to[out=-30,in=170] (2.5, 0.5)
          -- (1.5, 0.5) to[out=170,in=-30] (1.0, 0.0);

        \path [fill=white]
          (-1, 0) to[out=-30,in=170] (1, 0) to[out=-30,in=170] (1.5, 0.5)
          to[out=170,in=-30] (-0.5, 0.5) to[out=170,in=-30] (-1.0, 0.0);
        \path [fill=black, fill opacity=0.2]
          (-1, 0) to[out=-30,in=170] (1, 0) to[out=-30,in=170] (1.5, 0.5)
          to[out=170,in=-30] (-0.5, 0.5) to[out=170,in=-30] (-1.0, 0.0);

        \path [thick, opacity=0.5, draw, dashed] (-1, 0) to[out=-30,in=170] (1, 0);
        \path [thick, draw] (1, 0) to[out=-30,in=170] (1.5, 0.5);
        \path [thick, opacity=0.5, draw, dashed] (-0.5, 0.5) to[out=-30,in=170] (1.5, 0.5);
        \path [thick, opacity=0.5, draw, dashed] (-1.0, 0.0) to[out=-30,in=170] (-0.5, 0.5);

        \path [thick, opacity=0.5, draw, dashed] (1, 0) -- (2, 0);
        \path [thick, opacity=0.5, draw, dashed] (2, 0) to[out=-30,in=170] (2.5, 0.5);
        \path [thick, opacity=0.5, draw, dashed] (1.5, 0.5) -- (2.5, 0.5);

        \path [color=black,thick, draw] (0.60, 0.40) -- (1.15, 0.40) -- (1.65, 0.40);
        \path [color=black,thick, draw] (0.60, 0.40) -- (0.60, 0.15);
        \path [color=black,thick, draw] (1.65, 0.40) -- (1.65, 0.15);

        \path [color=black,thick, draw, dashed] (0.60, 0.15) -- (1.15, 0.15) -- (1.65, 0.15);

        \node at (0.0, 0.2) {\small $\Sigma^d$};
        \path[thick,draw,->] (0.3,0.12) -- (0.3,0.24);
        \node at (0.3, 0.30) {\small $\vn$};

        \path [color=black,thick, draw] (0.55, 0.70) -- (0.75, 0.38);
        \node at (0.5, 0.75) {\small $A$};
        \path[thick,draw,->] (1.28,0.28) -- (1.4,0.28);
        \node at (1.5, 0.25) {\small $\vec n$};
    \end{tikzpicture}
  \end{center}
  \caption{
    Choice of area element $A$; a curved, rectangular rectangle enclosing a
    point of the edge $\partial\Sigma^d$. We extend the sheet over the edge
    parallel in $\vec n$-direction, and assume $\sigma^d=0$ in the
    extension. The area shall be of vanishing length and height, and with
    top and bottom lines parallel to the extended sheet $\Sigma^d$.}
  \label{fig:rectangle}
\end{figure}
In a similar vain, let $A$ be an arbitrary area element perpendicular to
the edge; see Figure~\ref{fig:rectangle}. By virtue of Amp\`ere's law we
have
\begin{align}
  \int_{\partial A}\vH^d\cdot\text{d}\vec s = \int_A\vJ\cdot\vec\tau\dox,
  \label{eq:ampere}
\end{align}
where $\vec\tau$ is the unit vector in edge direction, orthogonal to $\vec
n$ and $\vn$. In the limit of vanishing length, we can rewrite the
left-hand side:
\begin{align*}
  \lim_{\text{length}\;\to\;0}
  \int_{\partial A}\vH^d\cdot\text{d}\vec s &\;=\;
  \lim_{\text{length}\;\to\;0}\left\{
  \int_{\text{left,right}}(\vec n\times\vH^d)\cdot\vec\tau\,\text{d}s
  \;+\;
  \int_{\text{top,bottom}}(\pm\vec\nu\times\vH^d)\cdot\vec\tau\,\text{d}s
  \right\}
  \\
  &\;=\;
  \lim_{\text{length}\;\to\;0}
  \int_{\text{left,right}}(\vec n\times\vH^d)\cdot\vec\tau\,\text{d}s.
\end{align*}
Exploiting the fact that $\vJ = \vJ_a + \delta_{\Sigma^d}\sigma^d\vE^d +
\delta_{\partial\Sigma^d}\lambda^d\vE^d$ and by taking the limit of
vanishing length we conclude that:
\begin{align}
  \lambda^d\vE^d\Big|_{\partial\Sigma^d\cap A}\cdot\vec\tau
  \;=\;
  \lim_{\text{length}\;\to\;0}\,
  \int_{\text{left,right}}(\vec n\times\vH^d)\cdot\vec\tau\,\text{d}s.
  \label{eq:19ff}
\end{align}
This implies that the jump over $\vec n\times\vH^d$ must have a singular
point contribution:
\begin{align*}
  \jse{\vec n\times\vH^d}\cdot\vec\tau
  \;=\;
  \lambda^d\vE^d\Big|_{\partial\Sigma^d\cap A}\cdot\vec\tau.
\end{align*}
Here, we defined $\jse{.}$ rigorously as the corresponding limit in
\eqref{eq:19ff}. Note that the height of the area element $A$ was
arbitrarily chosen. Indeed, the actual value of $\jse{\vec n\times\vH^d}$
does not depend on the particular choice of the area element $A$ because it
corresponds directly to a residue of an analytic function $(\vec
n\times\vH^d)\cdot\vec\tau$. In thise sense, Definition
\eqref{eq:compatibility_condition} is an equivalent, sligthly less
technical definition.

\section{Derivation of interface and internal boundary condition from weak
formulation}
\label{sec:weak_formulation}

In this appendix we derive the strong formulation with all jump and
compatibility conditions starting from a varational formulation. The weak
formulation reads, find a vector field $E$ such that
\begin{multline}
  \int_\Omega\mu_0^{-1}\nabla\times\vE^d\,\cdot\nabla\times\overline{\vec\psi}\,\dx
  - \omega^2\int_\Omega\e\vE^d\cdot\overline{\vec\psi}\,\dx
  \\
  - i\omega\int_{\Sigma^d}\sigma^d\vE^d\cdot\overline{\vec\psi}\,\dx
  - i\omega\int_{\partial\Sigma^d}\lambda^d
    \vE^d\cdot\overline{\vec\psi}\,\text{d}s
  \;=\;
  \int_\Omega i\omega \vec J_a\cdot\overline{\vec\psi}\,\dx,
  \label{eq:weak_formulation}
\end{multline}
for all smooth, vector-valued test functions $\vec\psi$ with compact
support in $\Omega$. Let us now define
\begin{align}
  i\omega\mu_0\int_\Omega\vH^d\cdot\overline{\vec\psi}\,\dx
  :=
  \int_\Omega\vE^d\cdot(\nabla\times\overline{\vec\psi})\,\dx.
  \label{eq:weak_h}
\end{align}
Integrating \eqref{eq:weak_h} by parts yields
\begin{align*}
  i\omega\mu_0\int_\Omega\vH^d\cdot\overline{\vec\psi}\,\dx
  =
  \int_\Omega\big(\nabla\times\vE^d\big)\cdot\overline{\vec\psi}\,\dx
  -
  \int_{\Sigma^d}\jsd{\nu\times\vE^d}\cdot\overline{\vec\psi}\,\dox.
\end{align*}
Thus, testing with (a) a smooth, vector-valued test function $\vec\psi$
with $\vec\psi=0$ on $\Sigma^d$, and (b) a sequence $\vec\psi_h$ of test
functions with vanishing support outside $\Sigma^d$ gives
\begin{align*}
  i\omega\mu_0\vH^d = \nabla\times\vE^d
  \quad\text{in }\Omega\setminus\Sigma^d,
  \qquad
  \jsd{\nu\times\vE^d} = 0
  \quad\text{on }\Sigma^d.
\end{align*}
Similarly, integration by parts of \eqref{eq:weak_formulation} and
substituting $\vH$:
\begin{multline*}
  i\omega\int_\Omega\big(\nabla\times\vH^d\big)\cdot\overline{\vec\psi}\dx
  -\omega^2\int_\Omega\e\vE^d\cdot\overline{\vec\psi}\dx
  - i\omega\int_\Omega\vec J_a\cdot\overline{\vec\psi}\dx
  \\
  =
  +i\omega\int_{\Sigma^d}\jsd{\vn\times\vH^d}\cdot
  \overline{\vec\psi}\,\dox
  +i\omega\int_{\partial\Sigma^d} \jse{\vec
  n\times\vH^d}\,\overline{\vec\psi}\, \text{d}s
  \\
  -i\omega\int_{\Sigma^d} \sigma^d \vE^d\cdot\overline{\vec\psi}\,\dox
  -i\omega\int_{\partial\Sigma^d}\lambda^d\vE^d
  \cdot\overline{\vec\psi}\,\text{d}s.
\end{multline*}
The occurence of the jump term over $\partial\Sigma^d$ after the
integration by parts has to be justified more precisely. Similarly, to the
discussion in Appendix~\ref{sec:gauss_law} we assume that the function
space for $\vH$ admits singular distributions on the edge. More precisely,
we define
\begin{multline}
  \int_{\partial\Sigma^d}\jse{\vec n\times\vH}\cdot\vec\psi\,\text{d}s
  \;:=\;
  \\
  \int_\Omega(\nabla\times\vH)\cdot\vec\psi\dx -
  \int_\Omega\vH\cdot(\nabla\times\vec\psi)\dx -
  \int_{\Sigma^d}\jsd{\vn\times\vH}\cdot\vec\psi\dox.
  \label{eq:integration_by_parts}
\end{multline}
Utilizing the same sequences (a) and (b) of test functions yields a similar
result:
\begin{align*}
  \begin{aligned}
  \nabla\times\big(i\omega\vH^d\big) -\omega^2\e\vE^d -i\omega\vec J_a &= 0
  &\quad&\text{in }\Omega\setminus\Sigma^d,
  \\
  i\omega\jsd{\vn\times\vH^d} &= i\omega\sigma^d\vE^d
  &\quad&\text{on }\Sigma^d\setminus{\partial\Sigma^d},
  \\
  i\omega\,\jse{\vec n\times\vH^d} &= i\omega\lambda^d\vE^d
  &\quad&\text{on }\partial\Sigma^d.
  \end{aligned}
\end{align*}

Now, let $\varphi$ be an arbitrary scalar-valued test function with compact
support and set $\vec\psi=\nabla\varphi$. And choose again (a) $\varphi=0$
on $\Sigma^d$, and (b) a sequence $\varphi_h$ of test functions with
vanishing support outside $\Sigma^d$. Testing \eqref{eq:weak_h} and
subsequent integration by parts results in
\begin{align*}
  \nabla\cdot\vH^d = 0
  \quad\text{in }\Omega\setminus\Sigma^d,
  \qquad
  \jsd{\vn\cdot\vH^d} = 0
  \quad\text{on }\Sigma^d.
\end{align*}
In case of the first equation we start again at
\eqref{eq:weak_formulation}. Utilizing the vector identity
$\nabla\times(\nabla\varphi)=0$:
\begin{align*}
  - \omega^2\int_\Omega\e\vE^d\cdot\nabla{\overline\varphi}\,\dx
  - i\omega\int_{\Sigma^d}\sigma^d\vE^d\cdot\nabla{\overline\varphi}\,\dox
  - i\omega\int_{\partial\Sigma^d}\lambda^d
  \vE^d\cdot\nabla{\overline\varphi}\,\text{d}s
  =
   i\omega \int_\Omega \vec J_a\cdot\nabla{\overline\varphi}\,\dx.
\end{align*}
Integration by parts:
\begin{multline*}
    \omega^2\int_\Omega\nabla\cdot\big(\e\vE^d\big)\overline\varphi\,\dx
  + i\omega \int_\Omega \nabla\cdot\vec J_a{\overline\varphi}\,\dx
  \\
  =
  - \omega^2\int_{\Sigma^d}\jsd{\nu\cdot
  \big(\e\vE^d)}\overline\varphi\,\dox
  - i\omega\int_{\Sigma^d}\nabla\cdot
  \big(\sigma^d\vE^d\big)\overline\varphi\,\dox
  \\
  + i\omega\int_{\partial\Sigma^d}\vec n \cdot \big(\sigma^d\vE^d\big)
  \overline\varphi\,\text{d}s
  - i\omega\int_{\partial\Sigma^d}\nabla\cdot
    \big(\lambda^d\vE^d\big)\overline\varphi\,\text{d}s.
\end{multline*}
Here, $\vec n$ denotes the outward-pointing unit vector tangential to
$\Sigma^d$ and normal to $\partial\Sigma^d$. We point out that for the
integration by parts of the interface term $\int_{\Sigma^d}\sigma^d\vE^d$
it is essential that $\sigma^d$ projects onto the tangential space of
$\Sigma^d$. We thus recover
\begin{align*}
  \nabla\cdot\big(\e\vE^d\big) = \frac1{i\omega}\nabla\cdot\vec J_a
  \quad\text{in }\Omega\setminus\Sigma^d,
  \qquad
  \jsd{\nu\cdot\big(\e\vE^d)} =
  \frac1{i\omega}\nabla\cdot\big(\sigma^d\vE^d\big)
  \quad\text{on }\Sigma^d,
\end{align*}
and
\begin{align*}
  \vec n \cdot \big(\sigma^d\vE^d\big) =
  \nabla\cdot \big(\lambda^d\vE^d\big)
  \quad\text{on }\partial\Sigma^d.
\end{align*}


\section*{Acknowledgments}

We acknowledge support by ARO MURI Award W911NF-14-0247 (MMai, MMat,
EK, ML, DM); EFRI 2-DARE NSF Grant 1542807 (MMat); and NSF DMS-1412769
(DM).



\end{document}

%% file: nano-ribbons.tex
\begin{tikzpicture}[gnuplot, xscale=0.45, yscale=0.5]
\gpmonochromelines
\path (0.000,0.000) rectangle (12.500,8.750);
\gpcolor{color=gp lt color axes}
\gpsetlinetype{gp lt axes}
\gpsetdashtype{gp dt axes}
\gpsetlinewidth{0.50}
\draw[gp path] (1.012,0.985)--(11.947,0.985);
\gpcolor{color=gp lt color border}
\gpsetlinetype{gp lt border}
\gpsetdashtype{gp dt solid}
\gpsetlinewidth{1.00}
\draw[gp path] (1.012,0.985)--(1.192,0.985);
\draw[gp path] (11.947,0.985)--(11.767,0.985);
\gpcolor{color=gp lt color axes}
\gpsetlinetype{gp lt axes}
\gpsetdashtype{gp dt axes}
\gpsetlinewidth{0.50}
\draw[gp path] (1.012,1.813)--(11.947,1.813);
\gpcolor{color=gp lt color border}
\gpsetlinetype{gp lt border}
\gpsetdashtype{gp dt solid}
\gpsetlinewidth{1.00}
\draw[gp path] (1.012,1.813)--(1.192,1.813);
\draw[gp path] (11.947,1.813)--(11.767,1.813);
\node[gp node right] at (0.828,1.813) {$-20$};
\gpcolor{color=gp lt color axes}
\gpsetlinetype{gp lt axes}
\gpsetdashtype{gp dt axes}
\gpsetlinewidth{0.50}
\draw[gp path] (1.012,2.642)--(11.947,2.642);
\gpcolor{color=gp lt color border}
\gpsetlinetype{gp lt border}
\gpsetdashtype{gp dt solid}
\gpsetlinewidth{1.00}
\draw[gp path] (1.012,2.642)--(1.192,2.642);
\draw[gp path] (11.947,2.642)--(11.767,2.642);
\gpcolor{color=gp lt color axes}
\gpsetlinetype{gp lt axes}
\gpsetdashtype{gp dt axes}
\gpsetlinewidth{0.50}
\draw[gp path] (1.012,3.470)--(11.947,3.470);
\gpcolor{color=gp lt color border}
\gpsetlinetype{gp lt border}
\gpsetdashtype{gp dt solid}
\gpsetlinewidth{1.00}
\draw[gp path] (1.012,3.470)--(1.192,3.470);
\draw[gp path] (11.947,3.470)--(11.767,3.470);
\node[gp node right] at (0.828,3.470) {$0$};
\gpcolor{color=gp lt color axes}
\gpsetlinetype{gp lt axes}
\gpsetdashtype{gp dt axes}
\gpsetlinewidth{0.50}
\draw[gp path] (1.012,4.299)--(11.947,4.299);
\gpcolor{color=gp lt color border}
\gpsetlinetype{gp lt border}
\gpsetdashtype{gp dt solid}
\gpsetlinewidth{1.00}
\draw[gp path] (1.012,4.299)--(1.192,4.299);
\draw[gp path] (11.947,4.299)--(11.767,4.299);
\gpcolor{color=gp lt color axes}
\gpsetlinetype{gp lt axes}
\gpsetdashtype{gp dt axes}
\gpsetlinewidth{0.50}
\draw[gp path] (1.012,5.127)--(11.947,5.127);
\gpcolor{color=gp lt color border}
\gpsetlinetype{gp lt border}
\gpsetdashtype{gp dt solid}
\gpsetlinewidth{1.00}
\draw[gp path] (1.012,5.127)--(1.192,5.127);
\draw[gp path] (11.947,5.127)--(11.767,5.127);
\node[gp node right] at (0.828,5.127) {$20$};
\gpcolor{color=gp lt color axes}
\gpsetlinetype{gp lt axes}
\gpsetdashtype{gp dt axes}
\gpsetlinewidth{0.50}
\draw[gp path] (1.012,5.956)--(11.947,5.956);
\gpcolor{color=gp lt color border}
\gpsetlinetype{gp lt border}
\gpsetdashtype{gp dt solid}
\gpsetlinewidth{1.00}
\draw[gp path] (1.012,5.956)--(1.192,5.956);
\draw[gp path] (11.947,5.956)--(11.767,5.956);
\gpcolor{color=gp lt color axes}
\gpsetlinetype{gp lt axes}
\gpsetdashtype{gp dt axes}
\gpsetlinewidth{0.50}
\draw[gp path] (1.012,6.784)--(11.947,6.784);
\gpcolor{color=gp lt color border}
\gpsetlinetype{gp lt border}
\gpsetdashtype{gp dt solid}
\gpsetlinewidth{1.00}
\draw[gp path] (1.012,6.784)--(1.192,6.784);
\draw[gp path] (11.947,6.784)--(11.767,6.784);
\node[gp node right] at (0.828,6.784) {$40$};
\gpcolor{color=gp lt color axes}
\gpsetlinetype{gp lt axes}
\gpsetdashtype{gp dt axes}
\gpsetlinewidth{0.50}
\draw[gp path] (1.012,7.613)--(11.947,7.613);
\gpcolor{color=gp lt color border}
\gpsetlinetype{gp lt border}
\gpsetdashtype{gp dt solid}
\gpsetlinewidth{1.00}
\draw[gp path] (1.012,7.613)--(1.192,7.613);
\draw[gp path] (11.947,7.613)--(11.767,7.613);
\gpcolor{color=gp lt color axes}
\gpsetlinetype{gp lt axes}
\gpsetdashtype{gp dt axes}
\gpsetlinewidth{0.50}
\draw[gp path] (1.012,8.441)--(11.947,8.441);
\gpcolor{color=gp lt color border}
\gpsetlinetype{gp lt border}
\gpsetdashtype{gp dt solid}
\gpsetlinewidth{1.00}
\draw[gp path] (1.012,8.441)--(1.192,8.441);
\draw[gp path] (11.947,8.441)--(11.767,8.441);
\node[gp node right] at (0.828,8.441) {$60$};
\gpcolor{color=gp lt color axes}
\gpsetlinetype{gp lt axes}
\gpsetdashtype{gp dt axes}
\gpsetlinewidth{0.50}
\draw[gp path] (1.012,0.985)--(1.012,8.441);
\gpcolor{color=gp lt color border}
\gpsetlinetype{gp lt border}
\gpsetdashtype{gp dt solid}
\gpsetlinewidth{1.00}
\draw[gp path] (1.012,0.985)--(1.012,1.165);
\draw[gp path] (1.012,8.441)--(1.012,8.261);
\node[gp node center] at (1.012,0.3) {$1.5$};
\gpcolor{color=gp lt color axes}
\gpsetlinetype{gp lt axes}
\gpsetdashtype{gp dt axes}
\gpsetlinewidth{0.50}
\draw[gp path] (3.746,0.985)--(3.746,8.441);
\gpcolor{color=gp lt color border}
\gpsetlinetype{gp lt border}
\gpsetdashtype{gp dt solid}
\gpsetlinewidth{1.00}
\draw[gp path] (3.746,0.985)--(3.746,1.165);
\draw[gp path] (3.746,8.441)--(3.746,8.261);
\node[gp node center] at (3.746,0.3) {$2$};
\gpcolor{color=gp lt color axes}
\gpsetlinetype{gp lt axes}
\gpsetdashtype{gp dt axes}
\gpsetlinewidth{0.50}
\draw[gp path] (6.480,0.985)--(6.480,8.441);
\gpcolor{color=gp lt color border}
\gpsetlinetype{gp lt border}
\gpsetdashtype{gp dt solid}
\gpsetlinewidth{1.00}
\draw[gp path] (6.480,0.985)--(6.480,1.165);
\draw[gp path] (6.480,8.441)--(6.480,8.261);
\node[gp node center] at (6.480,0.3) {$2.5$};
\gpcolor{color=gp lt color axes}
\gpsetlinetype{gp lt axes}
\gpsetdashtype{gp dt axes}
\gpsetlinewidth{0.50}
\draw[gp path] (9.213,0.985)--(9.213,7.645);
\draw[gp path] (9.213,8.261)--(9.213,8.441);
\gpcolor{color=gp lt color border}
\gpsetlinetype{gp lt border}
\gpsetdashtype{gp dt solid}
\gpsetlinewidth{1.00}
\draw[gp path] (9.213,0.985)--(9.213,1.165);
\draw[gp path] (9.213,8.441)--(9.213,8.261);
\node[gp node center] at (9.213,0.3) {$3$};
\gpcolor{color=gp lt color axes}
\gpsetlinetype{gp lt axes}
\gpsetdashtype{gp dt axes}
\gpsetlinewidth{0.50}
\draw[gp path] (11.947,0.985)--(11.947,8.441);
\gpcolor{color=gp lt color border}
\gpsetlinetype{gp lt border}
\gpsetdashtype{gp dt solid}
\gpsetlinewidth{1.00}
\draw[gp path] (11.947,0.985)--(11.947,1.165);
\draw[gp path] (11.947,8.441)--(11.947,8.261);
\node[gp node center] at (11.947,0.3) {$3.5$};
\draw[gp path] (1.012,8.441)--(1.012,0.985)--(11.947,0.985)--(11.947,8.441)--cycle;
\node[gp node center] at (6.479,-0.415) {$\tilde\omega$};
\node[gp node right] at (10.479,7.107) {Re $\varepsilon_{11}^{\text{R}}$};
\gpsetlinewidth{3.00}
\draw[gp path] (10.663,7.107)--(11.779,7.107);
\path[fill=red, opacity=0.1]
  (5.113,3.470)--(5.118,3.120)--(5.124,2.672)%
  --(5.129,2.268)--(5.134,1.924)--(5.140,1.648)--(5.145,1.440)--(5.151,1.293)--(5.156,1.201)%
  --(5.162,1.151)--(5.167,1.136)--(5.195,1.325)--(5.222,1.632)--(5.249,1.905)--(5.277,2.126)%
  --(5.304,2.301)--(5.331,2.441)--(5.359,2.555)--(5.386,2.649)--(5.413,2.728)--(5.441,2.795)%
  --(5.468,2.852)--(5.495,2.901)--(5.523,2.944)--(5.550,2.982)--(5.577,3.015)--(5.605,3.045)%
  --(5.632,3.072)--(5.659,3.096)--(5.687,3.118)--(5.714,3.138)--(5.741,3.157)--(5.769,3.174)%
  --(5.796,3.189)--(5.823,3.203)--(5.851,3.217)--(5.878,3.229)--(5.905,3.241)--(5.933,3.251)%
  --(5.960,3.261)--(5.987,3.271)--(6.015,3.280)--(6.042,3.288)--(6.069,3.296)--(6.097,3.303)%
  --(6.124,3.310)--(6.151,3.317)--(6.179,3.323)--(6.206,3.329)--(6.233,3.335)--(6.261,3.340)%
  --(6.288,3.346)--(6.315,3.350)--(6.343,3.355)--(6.370,3.360)--(6.397,3.364)--(6.425,3.368)%
  --(6.452,3.372)--(6.480,3.376)--(6.753,3.407)--(7.026,3.429)--(7.300,3.446)--(7.573,3.459)%
  --(7.846,3.469)--(8.120,3.470);
\draw[gp path] (1.012,3.629)--(1.285,3.633)--(1.559,3.638)--(1.832,3.644)--(2.105,3.651)%
  --(2.379,3.659)--(2.652,3.669)--(2.926,3.682)--(3.199,3.699)--(3.472,3.721)--(3.746,3.751)%
  --(3.773,3.755)--(3.800,3.759)--(3.828,3.763)--(3.855,3.768)--(3.882,3.772)--(3.910,3.777)%
  --(3.937,3.782)--(3.964,3.787)--(3.992,3.792)--(4.019,3.798)--(4.046,3.804)--(4.074,3.810)%
  --(4.101,3.816)--(4.128,3.823)--(4.156,3.831)--(4.183,3.838)--(4.210,3.846)--(4.238,3.855)%
  --(4.265,3.864)--(4.293,3.874)--(4.320,3.885)--(4.347,3.896)--(4.375,3.908)--(4.402,3.921)%
  --(4.429,3.935)--(4.457,3.950)--(4.484,3.966)--(4.511,3.984)--(4.539,4.004)--(4.566,4.025)%
  --(4.593,4.049)--(4.621,4.075)--(4.648,4.104)--(4.675,4.136)--(4.703,4.173)--(4.730,4.214)%
  --(4.757,4.261)--(4.785,4.315)--(4.812,4.379)--(4.839,4.453)--(4.867,4.541)--(4.894,4.648)%
  --(4.921,4.778)--(4.949,4.939)--(4.976,5.140)--(5.003,5.385)--(5.031,5.656)--(5.058,5.826)%
  --(5.063,5.815)--(5.069,5.775)--(5.074,5.698)--(5.080,5.575)--(5.085,5.396)--(5.091,5.154)%
  --(5.096,4.845)--(5.102,4.472)--(5.107,4.046)--(5.113,3.587)--(5.118,3.120)--(5.124,2.672)%
  --(5.129,2.268)--(5.134,1.924)--(5.140,1.648)--(5.145,1.440)--(5.151,1.293)--(5.156,1.201)%
  --(5.162,1.151)--(5.167,1.136)--(5.195,1.325)--(5.222,1.632)--(5.249,1.905)--(5.277,2.126)%
  --(5.304,2.301)--(5.331,2.441)--(5.359,2.555)--(5.386,2.649)--(5.413,2.728)--(5.441,2.795)%
  --(5.468,2.852)--(5.495,2.901)--(5.523,2.944)--(5.550,2.982)--(5.577,3.015)--(5.605,3.045)%
  --(5.632,3.072)--(5.659,3.096)--(5.687,3.118)--(5.714,3.138)--(5.741,3.157)--(5.769,3.174)%
  --(5.796,3.189)--(5.823,3.203)--(5.851,3.217)--(5.878,3.229)--(5.905,3.241)--(5.933,3.251)%
  --(5.960,3.261)--(5.987,3.271)--(6.015,3.280)--(6.042,3.288)--(6.069,3.296)--(6.097,3.303)%
  --(6.124,3.310)--(6.151,3.317)--(6.179,3.323)--(6.206,3.329)--(6.233,3.335)--(6.261,3.340)%
  --(6.288,3.346)--(6.315,3.350)--(6.343,3.355)--(6.370,3.360)--(6.397,3.364)--(6.425,3.368)%
  --(6.452,3.372)--(6.480,3.376)--(6.753,3.407)--(7.026,3.429)--(7.300,3.446)--(7.573,3.459)%
  --(7.846,3.469)--(8.120,3.478)--(8.393,3.485)--(8.667,3.491)--(8.940,3.496)--(9.213,3.500)%
  --(9.487,3.504)--(9.760,3.508)--(10.033,3.511)--(10.307,3.513)--(10.580,3.516)--(10.853,3.518)%
  --(11.127,3.520)--(11.400,3.522)--(11.674,3.523)--(11.947,3.525);
\node[gp node right] at (10.479,6.099) {Im $\varepsilon_{11}^{\text{R}}$};
\draw[gp path, dashed] (10.663,6.099)--(11.779,6.099);
\draw[gp path, dashed] (1.012,3.471)--(1.285,3.471)--(1.559,3.471)--(1.832,3.472)--(2.105,3.472)%
  --(2.379,3.472)--(2.652,3.473)--(2.926,3.473)--(3.199,3.474)--(3.472,3.476)--(3.746,3.478)%
  --(3.773,3.478)--(3.800,3.478)--(3.828,3.479)--(3.855,3.479)--(3.882,3.480)--(3.910,3.480)%
  --(3.937,3.480)--(3.964,3.481)--(3.992,3.481)--(4.019,3.482)--(4.046,3.483)--(4.074,3.483)%
  --(4.101,3.484)--(4.128,3.485)--(4.156,3.486)--(4.183,3.486)--(4.210,3.487)--(4.238,3.488)%
  --(4.265,3.490)--(4.293,3.491)--(4.320,3.492)--(4.347,3.494)--(4.375,3.496)--(4.402,3.498)%
  --(4.429,3.500)--(4.457,3.502)--(4.484,3.505)--(4.511,3.509)--(4.539,3.512)--(4.566,3.516)%
  --(4.593,3.521)--(4.621,3.527)--(4.648,3.534)--(4.675,3.542)--(4.703,3.551)--(4.730,3.563)%
  --(4.757,3.578)--(4.785,3.596)--(4.812,3.618)--(4.839,3.648)--(4.867,3.687)--(4.894,3.741)%
  --(4.921,3.817)--(4.949,3.928)--(4.976,4.100)--(5.003,4.379)--(5.031,4.862)--(5.058,5.724)%
  --(5.063,5.961)--(5.069,6.221)--(5.074,6.502)--(5.080,6.799)--(5.085,7.103)--(5.091,7.401)%
  --(5.096,7.676)--(5.102,7.907)--(5.107,8.073)--(5.113,8.156)--(5.118,8.147)--(5.124,8.047)%
  --(5.129,7.867)--(5.134,7.626)--(5.140,7.345)--(5.145,7.044)--(5.151,6.740)--(5.156,6.446)%
  --(5.162,6.169)--(5.167,5.913)--(5.195,4.969)--(5.222,4.440)--(5.249,4.136)--(5.277,3.951)%
  --(5.304,3.832)--(5.331,3.752)--(5.359,3.695)--(5.386,3.654)--(5.413,3.623)--(5.441,3.599)%
  --(5.468,3.580)--(5.495,3.565)--(5.523,3.553)--(5.550,3.543)--(5.577,3.535)--(5.605,3.528)%
  --(5.632,3.522)--(5.659,3.517)--(5.687,3.513)--(5.714,3.509)--(5.741,3.506)--(5.769,3.503)%
  --(5.796,3.500)--(5.823,3.498)--(5.851,3.496)--(5.878,3.494)--(5.905,3.493)--(5.933,3.491)%
  --(5.960,3.490)--(5.987,3.489)--(6.015,3.488)--(6.042,3.487)--(6.069,3.486)--(6.097,3.485)%
  --(6.124,3.484)--(6.151,3.483)--(6.179,3.483)--(6.206,3.482)--(6.233,3.481)--(6.261,3.481)%
  --(6.288,3.480)--(6.315,3.480)--(6.343,3.480)--(6.370,3.479)--(6.397,3.479)--(6.425,3.478)%
  --(6.452,3.478)--(6.480,3.478)--(6.753,3.476)--(7.026,3.474)--(7.300,3.473)--(7.573,3.473)%
  --(7.846,3.472)--(8.120,3.472)--(8.393,3.472)--(8.667,3.471)--(8.940,3.471)--(9.213,3.471)%
  --(9.487,3.471)--(9.760,3.471)--(10.033,3.471)--(10.307,3.471)--(10.580,3.471)--(10.853,3.471)%
  --(11.127,3.471)--(11.400,3.471)--(11.674,3.471)--(11.947,3.471);
\gpsetdashtype{gp dt solid}
\gpsetlinewidth{1.00}
\draw[gp path] (1.012,8.441)--(1.012,0.985)--(11.947,0.985)--(11.947,8.441)--cycle;
\gpdefrectangularnode{gp plot 1}{\pgfpoint{1.012cm}{0.985cm}}{\pgfpoint{11.947cm}{8.441cm}}
\end{tikzpicture}

%% file: nano-tubes.tex
\begin{tikzpicture}[gnuplot, xscale=0.45, yscale=0.5]
\gpmonochromelines
\path (0.000,0.000) rectangle (12.500,8.750);
\gpcolor{color=gp lt color axes}
\gpsetlinetype{gp lt axes}
\gpsetdashtype{gp dt axes}
\gpsetlinewidth{0.50}
\draw[gp path] (1.196,0.985)--(11.947,0.985);
\gpcolor{color=gp lt color border}
\gpsetlinetype{gp lt border}
\gpsetdashtype{gp dt solid}
\gpsetlinewidth{1.00}
\draw[gp path] (1.196,0.985)--(1.376,0.985);
\draw[gp path] (11.947,0.985)--(11.767,0.985);
\node[gp node right] at (1.012,0.985) {$-100$};
\gpcolor{color=gp lt color axes}
\gpsetlinetype{gp lt axes}
\gpsetdashtype{gp dt axes}
\gpsetlinewidth{0.50}
\draw[gp path] (1.196,2.476)--(11.947,2.476);
\gpcolor{color=gp lt color border}
\gpsetlinetype{gp lt border}
\gpsetdashtype{gp dt solid}
\gpsetlinewidth{1.00}
\draw[gp path] (1.196,2.476)--(1.376,2.476);
\draw[gp path] (11.947,2.476)--(11.767,2.476);
\node[gp node right] at (1.012,2.476) {$-50$};
\gpcolor{color=gp lt color axes}
\gpsetlinetype{gp lt axes}
\gpsetdashtype{gp dt axes}
\gpsetlinewidth{0.50}
\draw[gp path] (1.196,3.967)--(11.947,3.967);
\gpcolor{color=gp lt color border}
\gpsetlinetype{gp lt border}
\gpsetdashtype{gp dt solid}
\gpsetlinewidth{1.00}
\draw[gp path] (1.196,3.967)--(1.376,3.967);
\draw[gp path] (11.947,3.967)--(11.767,3.967);
\node[gp node right] at (1.012,3.967) {$0$};
\gpcolor{color=gp lt color axes}
\gpsetlinetype{gp lt axes}
\gpsetdashtype{gp dt axes}
\gpsetlinewidth{0.50}
\draw[gp path] (1.196,5.459)--(11.947,5.459);
\gpcolor{color=gp lt color border}
\gpsetlinetype{gp lt border}
\gpsetdashtype{gp dt solid}
\gpsetlinewidth{1.00}
\draw[gp path] (1.196,5.459)--(1.376,5.459);
\draw[gp path] (11.947,5.459)--(11.767,5.459);
\node[gp node right] at (1.012,5.459) {$50$};
\gpcolor{color=gp lt color axes}
\gpsetlinetype{gp lt axes}
\gpsetdashtype{gp dt axes}
\gpsetlinewidth{0.50}
\draw[gp path] (1.196,6.950)--(11.947,6.950);
\gpcolor{color=gp lt color border}
\gpsetlinetype{gp lt border}
\gpsetdashtype{gp dt solid}
\gpsetlinewidth{1.00}
\draw[gp path] (1.196,6.950)--(1.376,6.950);
\draw[gp path] (11.947,6.950)--(11.767,6.950);
\node[gp node right] at (1.012,6.950) {$100$};
\gpcolor{color=gp lt color axes}
\gpsetlinetype{gp lt axes}
\gpsetdashtype{gp dt axes}
\gpsetlinewidth{0.50}
\draw[gp path] (1.196,8.441)--(11.947,8.441);
\gpcolor{color=gp lt color border}
\gpsetlinetype{gp lt border}
\gpsetdashtype{gp dt solid}
\gpsetlinewidth{1.00}
\draw[gp path] (1.196,8.441)--(1.376,8.441);
\draw[gp path] (11.947,8.441)--(11.767,8.441);
\node[gp node right] at (1.012,8.441) {$150$};
\gpcolor{color=gp lt color axes}
\gpsetlinetype{gp lt axes}
\gpsetdashtype{gp dt axes}
\gpsetlinewidth{0.50}
\draw[gp path] (1.196,0.985)--(1.196,8.441);
\gpcolor{color=gp lt color border}
\gpsetlinetype{gp lt border}
\gpsetdashtype{gp dt solid}
\gpsetlinewidth{1.00}
\draw[gp path] (1.196,0.985)--(1.196,1.165);
\draw[gp path] (1.196,8.441)--(1.196,8.261);
\node[gp node center] at (1.196,0.3) {$0.5$};
\gpcolor{color=gp lt color axes}
\gpsetlinetype{gp lt axes}
\gpsetdashtype{gp dt axes}
\gpsetlinewidth{0.50}
\draw[gp path] (3.884,0.985)--(3.884,7.645);
\draw[gp path] (3.884,8.261)--(3.884,8.441);
\gpcolor{color=gp lt color border}
\gpsetlinetype{gp lt border}
\gpsetdashtype{gp dt solid}
\gpsetlinewidth{1.00}
\draw[gp path] (3.884,0.985)--(3.884,1.165);
\draw[gp path] (3.884,8.441)--(3.884,8.261);
\node[gp node center] at (3.884,0.3) {$1$};
\gpcolor{color=gp lt color axes}
\gpsetlinetype{gp lt axes}
\gpsetdashtype{gp dt axes}
\gpsetlinewidth{0.50}
\draw[gp path] (6.572,0.985)--(6.572,8.441);
\gpcolor{color=gp lt color border}
\gpsetlinetype{gp lt border}
\gpsetdashtype{gp dt solid}
\gpsetlinewidth{1.00}
\draw[gp path] (6.572,0.985)--(6.572,1.165);
\draw[gp path] (6.572,8.441)--(6.572,8.261);
\node[gp node center] at (6.572,0.3) {$1.5$};
\gpcolor{color=gp lt color axes}
\gpsetlinetype{gp lt axes}
\gpsetdashtype{gp dt axes}
\gpsetlinewidth{0.50}
\draw[gp path] (9.259,0.985)--(9.259,8.441);
\gpcolor{color=gp lt color border}
\gpsetlinetype{gp lt border}
\gpsetdashtype{gp dt solid}
\gpsetlinewidth{1.00}
\draw[gp path] (9.259,0.985)--(9.259,1.165);
\draw[gp path] (9.259,8.441)--(9.259,8.261);
\node[gp node center] at (9.259,0.3) {$2$};
\gpcolor{color=gp lt color axes}
\gpsetlinetype{gp lt axes}
\gpsetdashtype{gp dt axes}
\gpsetlinewidth{0.50}
\draw[gp path] (11.947,0.985)--(11.947,8.441);
\gpcolor{color=gp lt color border}
\gpsetlinetype{gp lt border}
\gpsetdashtype{gp dt solid}
\gpsetlinewidth{1.00}
\draw[gp path] (11.947,0.985)--(11.947,1.165);
\draw[gp path] (11.947,8.441)--(11.947,8.261);
\node[gp node center] at (11.947,0.3) {$2.5$};
\draw[gp path] (1.196,8.441)--(1.196,0.985)--(11.947,0.985)--(11.947,8.441)--cycle;
\node[gp node center] at (6.571,-0.415) {$\tilde\omega$};
\node[gp node right] at (3.852,7.107) {Re $\varepsilon_{11}^{\text{T}}$};
\gpsetlinewidth{3.00}
\draw[gp path] (4.036,7.107)--(5.152,7.107);
\path[fill=red, opacity=0.1]
    (7.663,3.967)--(7.668,3.526)%
  --(7.673,3.139)--(7.679,2.794)--(7.684,2.502)--(7.690,2.271)--(7.695,2.099)--(7.700,1.981)%
  --(7.706,1.907)--(7.711,1.871)--(7.716,1.863)--(7.722,1.876)--(7.727,1.905)--(7.733,1.945)%
  --(7.738,1.992)--(7.743,2.044)--(7.749,2.098)--(7.754,2.154)--(7.759,2.209)--(7.765,2.264)%
  --(7.770,2.317)--(7.776,2.369)--(7.781,2.419)--(7.786,2.467)--(7.792,2.514)--(7.797,2.558)%
  --(7.802,2.600)--(7.808,2.641)--(7.813,2.679)--(7.819,2.716)--(7.824,2.751)--(7.829,2.785)%
  --(7.835,2.817)--(7.840,2.848)--(7.845,2.877)--(7.851,2.905)--(7.856,2.932)--(7.862,2.957)%
  --(7.867,2.982)--(7.872,3.005)--(7.878,3.028)--(7.883,3.049)--(7.888,3.070)--(7.894,3.090)%
  --(7.899,3.109)--(7.905,3.127)--(7.910,3.145)--(7.915,3.162)--(7.921,3.179)--(7.926,3.194)%
  --(7.932,3.210)--(7.937,3.224)--(7.942,3.239)--(7.948,3.252)--(7.953,3.266)--(7.958,3.278)%
  --(7.964,3.291)--(7.969,3.303)--(7.975,3.314)--(7.980,3.326)--(7.985,3.337)--(7.991,3.347)%
  --(7.996,3.357)--(8.001,3.367)--(8.007,3.377)--(8.012,3.386)--(8.018,3.395)--(8.023,3.404)%
  --(8.028,3.413)--(8.034,3.421)--(8.039,3.429)--(8.044,3.437)--(8.050,3.445)--(8.055,3.453)%
  --(8.061,3.460)--(8.066,3.467)--(8.071,3.474)--(8.077,3.481)--(8.082,3.487)--(8.087,3.494)%
  --(8.093,3.500)--(8.098,3.506)--(8.104,3.512)--(8.109,3.518)--(8.114,3.524)--(8.120,3.529)%
  --(8.125,3.535)--(8.130,3.540)--(8.141,3.551)--(8.147,3.556)--(8.152,3.560)--(8.157,3.565)%
  --(8.163,3.570)--(8.168,3.575)--(8.179,3.583)--(8.184,3.588)--(8.211,3.608)--(8.238,3.627)%
  --(8.265,3.644)--(8.292,3.659)--(8.319,3.673)--(8.345,3.686)--(8.372,3.698)--(8.399,3.710)%
  --(8.426,3.720)--(8.453,3.730)--(8.480,3.739)--(8.507,3.748)--(8.534,3.756)--(8.560,3.763)%
  --(8.587,3.770)--(8.614,3.777)--(8.641,3.783)--(8.668,3.789)--(8.695,3.795)--(8.722,3.800)%
  --(8.749,3.806)--(8.775,3.810)--(8.802,3.815)--(8.829,3.820)--(8.856,3.824)--(8.883,3.828)%
  --(8.910,3.832)--(8.937,3.835)--(8.964,3.839)--(8.990,3.842)--(9.017,3.846)--(9.044,3.849)%
  --(9.071,3.852)--(9.098,3.855)--(9.125,3.858)--(9.152,3.860)--(9.179,3.863)--(9.205,3.865)%
  --(9.232,3.868)--(9.259,3.870)--(9.528,3.890)--(9.797,3.905)--(10.066,3.916)--(10.334,3.925)%
  --(10.603,3.933)--(10.872,3.939)--(11.141,3.944)--(11.409,3.948)--(11.678,3.952)--(11.947,3.956);
\draw[gp path] (1.196,4.046)--(1.465,4.048)--(1.734,4.050)--(2.002,4.052)--(2.271,4.054)%
  --(2.540,4.056)--(2.809,4.058)--(3.077,4.061)--(3.346,4.063)--(3.615,4.067)--(3.884,4.070)%
  --(4.153,4.074)--(4.421,4.079)--(4.690,4.085)--(4.959,4.092)--(5.228,4.101)--(5.496,4.112)%
  --(5.765,4.126)--(6.034,4.145)--(6.303,4.171)--(6.572,4.211)--(6.598,4.216)--(6.625,4.221)%
  --(6.652,4.226)--(6.679,4.232)--(6.706,4.238)--(6.733,4.245)--(6.760,4.252)--(6.787,4.259)%
  --(6.813,4.267)--(6.840,4.275)--(6.867,4.284)--(6.894,4.293)--(6.921,4.303)--(6.948,4.314)%
  --(6.975,4.326)--(7.002,4.339)--(7.028,4.352)--(7.055,4.367)--(7.082,4.383)--(7.109,4.401)%
  --(7.114,4.404)--(7.120,4.408)--(7.125,4.412)--(7.131,4.416)--(7.136,4.420)--(7.141,4.424)%
  --(7.147,4.428)--(7.152,4.433)--(7.157,4.437)--(7.163,4.442)--(7.168,4.446)--(7.174,4.451)%
  --(7.179,4.455)--(7.184,4.460)--(7.190,4.465)--(7.195,4.470)--(7.200,4.475)--(7.206,4.481)%
  --(7.211,4.486)--(7.217,4.492)--(7.222,4.497)--(7.227,4.503)--(7.233,4.509)--(7.238,4.515)%
  --(7.243,4.521)--(7.249,4.528)--(7.254,4.534)--(7.260,4.541)--(7.265,4.548)--(7.270,4.555)%
  --(7.276,4.562)--(7.281,4.569)--(7.286,4.577)--(7.292,4.585)--(7.297,4.593)--(7.303,4.601)%
  --(7.308,4.609)--(7.313,4.618)--(7.319,4.627)--(7.324,4.636)--(7.329,4.646)--(7.340,4.666)%
  --(7.346,4.676)--(7.351,4.687)--(7.356,4.698)--(7.367,4.721)--(7.372,4.733)--(7.378,4.745)%
  --(7.383,4.758)--(7.389,4.772)--(7.394,4.786)--(7.399,4.800)--(7.405,4.815)--(7.410,4.830)%
  --(7.415,4.846)--(7.421,4.863)--(7.426,4.880)--(7.432,4.898)--(7.437,4.917)--(7.442,4.936)%
  --(7.448,4.956)--(7.453,4.977)--(7.458,4.999)--(7.464,5.021)--(7.469,5.045)--(7.475,5.070)%
  --(7.480,5.095)--(7.485,5.122)--(7.491,5.150)--(7.496,5.180)--(7.501,5.210)--(7.507,5.242)%
  --(7.512,5.275)--(7.518,5.310)--(7.523,5.346)--(7.528,5.384)--(7.534,5.424)--(7.539,5.465)%
  --(7.544,5.507)--(7.550,5.551)--(7.555,5.596)--(7.561,5.643)--(7.566,5.690)--(7.571,5.737)%
  --(7.577,5.783)--(7.582,5.828)--(7.587,5.871)--(7.593,5.908)--(7.598,5.938)--(7.604,5.958)%
  --(7.609,5.963)--(7.614,5.948)--(7.620,5.908)--(7.625,5.833)--(7.630,5.717)--(7.636,5.551)%
  --(7.641,5.329)--(7.647,5.050)--(7.652,4.715)--(7.657,4.337)--(7.663,3.933)--(7.668,3.526)%
  --(7.673,3.139)--(7.679,2.794)--(7.684,2.502)--(7.690,2.271)--(7.695,2.099)--(7.700,1.981)%
  --(7.706,1.907)--(7.711,1.871)--(7.716,1.863)--(7.722,1.876)--(7.727,1.905)--(7.733,1.945)%
  --(7.738,1.992)--(7.743,2.044)--(7.749,2.098)--(7.754,2.154)--(7.759,2.209)--(7.765,2.264)%
  --(7.770,2.317)--(7.776,2.369)--(7.781,2.419)--(7.786,2.467)--(7.792,2.514)--(7.797,2.558)%
  --(7.802,2.600)--(7.808,2.641)--(7.813,2.679)--(7.819,2.716)--(7.824,2.751)--(7.829,2.785)%
  --(7.835,2.817)--(7.840,2.848)--(7.845,2.877)--(7.851,2.905)--(7.856,2.932)--(7.862,2.957)%
  --(7.867,2.982)--(7.872,3.005)--(7.878,3.028)--(7.883,3.049)--(7.888,3.070)--(7.894,3.090)%
  --(7.899,3.109)--(7.905,3.127)--(7.910,3.145)--(7.915,3.162)--(7.921,3.179)--(7.926,3.194)%
  --(7.932,3.210)--(7.937,3.224)--(7.942,3.239)--(7.948,3.252)--(7.953,3.266)--(7.958,3.278)%
  --(7.964,3.291)--(7.969,3.303)--(7.975,3.314)--(7.980,3.326)--(7.985,3.337)--(7.991,3.347)%
  --(7.996,3.357)--(8.001,3.367)--(8.007,3.377)--(8.012,3.386)--(8.018,3.395)--(8.023,3.404)%
  --(8.028,3.413)--(8.034,3.421)--(8.039,3.429)--(8.044,3.437)--(8.050,3.445)--(8.055,3.453)%
  --(8.061,3.460)--(8.066,3.467)--(8.071,3.474)--(8.077,3.481)--(8.082,3.487)--(8.087,3.494)%
  --(8.093,3.500)--(8.098,3.506)--(8.104,3.512)--(8.109,3.518)--(8.114,3.524)--(8.120,3.529)%
  --(8.125,3.535)--(8.130,3.540)--(8.141,3.551)--(8.147,3.556)--(8.152,3.560)--(8.157,3.565)%
  --(8.163,3.570)--(8.168,3.575)--(8.179,3.583)--(8.184,3.588)--(8.211,3.608)--(8.238,3.627)%
  --(8.265,3.644)--(8.292,3.659)--(8.319,3.673)--(8.345,3.686)--(8.372,3.698)--(8.399,3.710)%
  --(8.426,3.720)--(8.453,3.730)--(8.480,3.739)--(8.507,3.748)--(8.534,3.756)--(8.560,3.763)%
  --(8.587,3.770)--(8.614,3.777)--(8.641,3.783)--(8.668,3.789)--(8.695,3.795)--(8.722,3.800)%
  --(8.749,3.806)--(8.775,3.810)--(8.802,3.815)--(8.829,3.820)--(8.856,3.824)--(8.883,3.828)%
  --(8.910,3.832)--(8.937,3.835)--(8.964,3.839)--(8.990,3.842)--(9.017,3.846)--(9.044,3.849)%
  --(9.071,3.852)--(9.098,3.855)--(9.125,3.858)--(9.152,3.860)--(9.179,3.863)--(9.205,3.865)%
  --(9.232,3.868)--(9.259,3.870)--(9.528,3.890)--(9.797,3.905)--(10.066,3.916)--(10.334,3.925)%
  --(10.603,3.933)--(10.872,3.939)--(11.141,3.944)--(11.409,3.948)--(11.678,3.952)--(11.947,3.956);
\node[gp node right] at (3.852,6.099) {Im $\varepsilon_{11}^{\text{T}}$};
\draw[gp path, dashed] (4.036,6.099)--(5.152,6.099);
\draw[gp path, dashed] (1.196,3.968)--(1.465,3.968)--(1.734,3.968)--(2.002,3.968)--(2.271,3.968)%
  --(2.540,3.968)--(2.809,3.968)--(3.077,3.968)--(3.346,3.968)--(3.615,3.968)--(3.884,3.968)%
  --(4.153,3.968)--(4.421,3.968)--(4.690,3.969)--(4.959,3.969)--(5.228,3.969)--(5.496,3.970)%
  --(5.765,3.971)--(6.034,3.972)--(6.303,3.974)--(6.572,3.977)--(6.598,3.978)--(6.625,3.978)%
  --(6.652,3.979)--(6.679,3.980)--(6.706,3.980)--(6.733,3.981)--(6.760,3.982)--(6.787,3.983)%
  --(6.813,3.984)--(6.840,3.985)--(6.867,3.986)--(6.894,3.987)--(6.921,3.989)--(6.948,3.990)%
  --(6.975,3.992)--(7.002,3.994)--(7.028,3.996)--(7.055,3.999)--(7.082,4.002)--(7.109,4.005)%
  --(7.114,4.006)--(7.120,4.007)--(7.125,4.008)--(7.131,4.009)--(7.136,4.009)--(7.141,4.010)%
  --(7.147,4.011)--(7.152,4.012)--(7.157,4.013)--(7.163,4.014)--(7.168,4.015)--(7.174,4.016)%
  --(7.179,4.017)--(7.184,4.018)--(7.190,4.019)--(7.195,4.021)--(7.200,4.022)--(7.206,4.023)%
  --(7.211,4.024)--(7.217,4.026)--(7.222,4.027)--(7.227,4.029)--(7.233,4.030)--(7.238,4.032)%
  --(7.243,4.033)--(7.249,4.035)--(7.254,4.037)--(7.260,4.039)--(7.265,4.040)--(7.270,4.042)%
  --(7.276,4.044)--(7.281,4.047)--(7.286,4.049)--(7.292,4.051)--(7.297,4.054)--(7.303,4.056)%
  --(7.308,4.059)--(7.313,4.061)--(7.319,4.064)--(7.324,4.067)--(7.329,4.071)--(7.340,4.077)%
  --(7.346,4.081)--(7.351,4.085)--(7.356,4.089)--(7.367,4.097)--(7.372,4.102)--(7.378,4.107)%
  --(7.383,4.112)--(7.389,4.118)--(7.394,4.124)--(7.399,4.130)--(7.405,4.136)--(7.410,4.143)%
  --(7.415,4.150)--(7.421,4.158)--(7.426,4.167)--(7.432,4.175)--(7.437,4.185)--(7.442,4.195)%
  --(7.448,4.206)--(7.453,4.217)--(7.458,4.230)--(7.464,4.243)--(7.469,4.257)--(7.475,4.273)%
  --(7.480,4.289)--(7.485,4.307)--(7.491,4.327)--(7.496,4.348)--(7.501,4.371)--(7.507,4.396)%
  --(7.512,4.423)--(7.518,4.453)--(7.523,4.486)--(7.528,4.523)--(7.534,4.562)--(7.539,4.607)%
  --(7.544,4.655)--(7.550,4.710)--(7.555,4.770)--(7.561,4.837)--(7.566,4.912)--(7.571,4.997)%
  --(7.577,5.091)--(7.582,5.198)--(7.587,5.317)--(7.593,5.452)--(7.598,5.603)--(7.604,5.773)%
  --(7.609,5.962)--(7.614,6.172)--(7.620,6.401)--(7.625,6.649)--(7.630,6.910)--(7.636,7.176)%
  --(7.641,7.435)--(7.647,7.672)--(7.652,7.868)--(7.657,8.005)--(7.663,8.068)--(7.668,8.051)%
  --(7.673,7.954)--(7.679,7.790)--(7.684,7.574)--(7.690,7.325)--(7.695,7.061)--(7.700,6.796)%
  --(7.706,6.540)--(7.711,6.300)--(7.716,6.079)--(7.722,5.878)--(7.727,5.697)--(7.733,5.536)%
  --(7.738,5.392)--(7.743,5.264)--(7.749,5.150)--(7.754,5.049)--(7.759,4.959)--(7.765,4.879)%
  --(7.770,4.807)--(7.776,4.743)--(7.781,4.685)--(7.786,4.634)--(7.792,4.587)--(7.797,4.545)%
  --(7.802,4.506)--(7.808,4.472)--(7.813,4.440)--(7.819,4.411)--(7.824,4.385)--(7.829,4.361)%
  --(7.835,4.338)--(7.840,4.318)--(7.845,4.299)--(7.851,4.282)--(7.856,4.266)--(7.862,4.251)%
  --(7.867,4.237)--(7.872,4.224)--(7.878,4.212)--(7.883,4.201)--(7.888,4.190)--(7.894,4.181)%
  --(7.899,4.171)--(7.905,4.163)--(7.910,4.155)--(7.915,4.147)--(7.921,4.140)--(7.926,4.133)%
  --(7.932,4.127)--(7.937,4.121)--(7.942,4.115)--(7.948,4.110)--(7.953,4.105)--(7.958,4.100)%
  --(7.964,4.095)--(7.969,4.091)--(7.975,4.087)--(7.980,4.083)--(7.985,4.079)--(7.991,4.076)%
  --(7.996,4.072)--(8.001,4.069)--(8.007,4.066)--(8.012,4.063)--(8.018,4.060)--(8.023,4.058)%
  --(8.028,4.055)--(8.034,4.052)--(8.039,4.050)--(8.044,4.048)--(8.050,4.046)--(8.055,4.044)%
  --(8.061,4.042)--(8.066,4.040)--(8.071,4.038)--(8.077,4.036)--(8.082,4.034)--(8.087,4.033)%
  --(8.093,4.031)--(8.098,4.029)--(8.104,4.028)--(8.109,4.026)--(8.114,4.025)--(8.120,4.024)%
  --(8.125,4.022)--(8.130,4.021)--(8.141,4.019)--(8.147,4.018)--(8.152,4.017)--(8.157,4.016)%
  --(8.163,4.015)--(8.168,4.014)--(8.179,4.012)--(8.184,4.011)--(8.211,4.007)--(8.238,4.003)%
  --(8.265,4.000)--(8.292,3.997)--(8.319,3.995)--(8.345,3.993)--(8.372,3.991)--(8.399,3.989)%
  --(8.426,3.988)--(8.453,3.986)--(8.480,3.985)--(8.507,3.984)--(8.534,3.983)--(8.560,3.982)%
  --(8.587,3.981)--(8.614,3.980)--(8.641,3.980)--(8.668,3.979)--(8.695,3.979)--(8.722,3.978)%
  --(8.749,3.977)--(8.775,3.977)--(8.802,3.977)--(8.829,3.976)--(8.856,3.976)--(8.883,3.975)%
  --(8.910,3.975)--(8.937,3.975)--(8.964,3.974)--(8.990,3.974)--(9.017,3.974)--(9.044,3.974)%
  --(9.071,3.973)--(9.098,3.973)--(9.125,3.973)--(9.152,3.973)--(9.179,3.973)--(9.205,3.972)%
  --(9.232,3.972)--(9.259,3.972)--(9.528,3.971)--(9.797,3.970)--(10.066,3.969)--(10.334,3.969)%
  --(10.603,3.969)--(10.872,3.969)--(11.141,3.968)--(11.409,3.968)--(11.678,3.968)--(11.947,3.968);
\gpsetdashtype{gp dt solid}
\gpsetlinewidth{1.00}
\draw[gp path] (1.196,8.441)--(1.196,0.985)--(11.947,0.985)--(11.947,8.441)--cycle;
\gpdefrectangularnode{gp plot 1}{\pgfpoint{1.196cm}{0.985cm}}{\pgfpoint{11.947cm}{8.441cm}}
\end{tikzpicture}